\documentclass[11pt,reqno]{amsart}
\usepackage{amsmath,amssymb,amsthm,mathtools,mathrsfs}
\usepackage{needspace}
\usepackage[a4paper,margin=30mm]{geometry}
\usepackage[hidelinks]{hyperref}
\newcommand{\doi}[1]{\href{https://doi.org/#1}{doi:#1}}
\numberwithin{equation}{section}
\newcommand{\HH}{\mathbb H}
\newcommand{\RR}{\mathbb R}
\newcommand{\Proj}{\mathcal P_{\HH}}
\newcommand{\ProjR}{\mathcal P_{\RR}}
\newcommand{\Real}{\mathcal R}

\newcommand{\trH}{\operatorname{tr}_{\mathbb H}}
\newcommand{\trR}{\operatorname{tr}_{\mathbb R}}
\newcommand{\Ree}{\operatorname{Re}}
\newcommand{\pair}[2]{\langle #1,\,#2\rangle}
\newcommand{\closure}[1]{\overline{#1}}
\newcommand{\dd}{\,\mathrm d}
\newcommand{\cL}{\mathcal L}
\newcommand{\bd}{\partial\Omega}
\newcommand{\ubar}{\underline u}
\newtheorem{theorem}{Theorem}[section]
\newtheorem{proposition}[theorem]{Proposition}
\newtheorem{lemma}[theorem]{Lemma}
\newtheorem{corollary}[theorem]{Corollary}

\theoremstyle{definition}
\newtheorem{definition}[theorem]{Definition}
\theoremstyle{plain}

\begin{document}

\title[Quaternionic $k$-Hessian equations]{The Dirichlet Problem for Quaternionic $k$-Hessian Equations}
\author[Xi-Nan Ma]{Xi-Nan Ma}
\address{School of Mathematical Sciences,
University of Science and Technology of China, Hefei, China}
\email{xinan@ustc.edu.cn}

\author{Yang Zhou}
\address{School of Mathematical Sciences,
University of Science and Technology of China, Hefei, China}
\email{zy19700816@mail.ustc.edu.cn}

\author[Jingyong Zhu]{Jingyong Zhu}
\address{School of Mathematical Sciences, Sichuan University,
Chengdu 610065, China}
\email{jzhu@scu.edu.cn}

\begin{abstract}
We study the Dirichlet problem for quaternionic $k$-Hessian equations on
smooth bounded domains. For a smooth positive density $f(q,z)$ satisfying
$f_z\ge0$, we prove existence and uniqueness of a smooth admissible
solution whenever a smooth admissible subsolution with the prescribed
boundary values exists. The proof extends the quaternionic Monge--Amp\`ere
subsolution method to $2\le k<n$ through boundary second derivative
estimates. We also construct subsolutions on strictly quaternionic
$(k-1)$-convex domains and characterize smooth solvability for constant
boundary values.
\end{abstract}
\maketitle

\section{Introduction}

Hessian equations form a basic class of fully nonlinear elliptic equations
in analysis and geometry. The Monge--Amp\`ere equation arises in the
Minkowski problem and the Calabi conjecture
\cite{ChengYau1976,Pogorelov1978,Yau1978}, while the $k$-Hessian equations
interpolate between the Laplace and Monge--Amp\`ere equations.
Caffarelli, Nirenberg and Spruck developed the classical Dirichlet theory
for real Hessian equations \cite{CNS1984-I,CNS1985-III}.
The complex theory was developed by Caffarelli, Kohn, Nirenberg and Spruck
\cite{CNS1985-II} and, for complex Hessian equations, by Li
\cite{Li2004}. Further results, including degenerate equations, were
obtained by Trudinger and Ivochkina--Trudinger--Wang
\cite{Trudinger1995,IvochkinaTrudingerWang2005}.

The subsolution method allows the boundary geometry required in these
theories to be encoded in a global admissible subsolution. The work of
Guan and Guan--Spruck, and its subsequent extensions, established
Dirichlet results under such hypotheses in real and complex settings
\cite{GuanSpruck1993,Guan1994,Guan1998,Guan2009,Guan2014,Guan2023}.
The subsolution supplies comparison functions and boundary barriers,
which are essential for estimates on general domains.

Alesker established the quaternionic linear algebra and pluripotential
framework for quaternionic Monge--Amp\`ere equations
\cite{AleskerLinearAlgebra,Alesker2003}.

For the Dirichlet problem, Zhu \cite{Zhu2016} proved that a smooth
admissible subsolution guarantees smooth solvability of the quaternionic
Monge--Amp\`ere equation on any smooth bounded domain.
The subsolution method also applies to flat quaternionic quotients Ma-Zhu
\cite{MaZhu2026}, and Wan \cite{Wan2020} established a subsolution
theorem in the weak solution setting.
Harvey--Lawson's general Dirichlet theory includes inhomogeneous
quaternionic $k$-Hessian equations and gives unique continuous viscosity
solutions under the corresponding strict boundary convexity assumptions
\cite[Theorem~2.11 and Example~6.4]{HarveyLawson2019}.
Liu--Wang \cite{LiuWang2023} developed the pluripotential theory of
quaternionic $m$-subharmonic functions.
Amal, Asserda and Barloub studied the homogeneous Dirichlet problem on the
unit ball \cite{AmalAsserdaBarloub2023} and finite-energy solutions with
positive Radon measure data \cite{AmalAsserdaBarloub2025}.

We study smooth solutions up to the boundary for quaternionic
$k$-Hessian equations on general smooth bounded domains admitting an
admissible subsolution. The main extension beyond the Monge--Amp\`ere
case is the range $2\le k<n$. In this range the admissible quaternionic
Hessian need not be positive definite, and the boundary second derivative
estimates require an argument adapted to the $k$-Hessian cone.

Let \(\Omega\subset\!\subset\mathbb H^n\) be a bounded domain with smooth
boundary. We seek a \(k\)-admissible solution of
\begin{equation}\label{eq:dirichlet}
 \sigma_k(H(u))=f(q,u)\quad\text{in }\Omega,
 \qquad u=\varphi\quad\text{on }\bd,
 \qquad 2\le k\le n.
\end{equation}
Here \(H(u)=\frac14(\partial_{\bar q_a}\partial_{q_b}u)_{a,b=1}^n\)
is the normalized quaternionic Hessian of a real-valued function. The
quaternionic linear algebra and the notion of \(k\)-admissibility are
given in Sections~\ref{sec:prelim} and~\ref{sec:hessian-tools}, respectively. Our main result is the following.

\begin{theorem}\label{thm:main}
Let \(2\le k\le n\) and let \(\Omega\subset\!\subset\HH^n\) be a bounded domain
with smooth boundary. Suppose that \(\varphi\in C^\infty(\bd)\) and
\[
 f(q,z)\in C^\infty(\closure\Omega\times\RR),\qquad f>0,
 \qquad f_z:=\frac{\partial f}{\partial z}\ge0.
\]
Assume that a real-valued function \(\ubar\in C^\infty(\closure\Omega)\) satisfies
\begin{equation}\label{eq:subsolution}
 H(\ubar)\in\Gamma_k\quad\text{in }\Omega,
 \qquad \sigma_k(H(\ubar))\ge f(q,\ubar)\quad\text{in }\Omega,
 \qquad \ubar=\varphi\quad\text{on }\bd.
\end{equation}
Then \eqref{eq:dirichlet} has a unique \(k\)-admissible solution
\(u\in C^\infty(\closure\Omega)\). Moreover, \(u\ge\ubar\).
\end{theorem}

 For $k=n$, the theorem belongs
the subsolution result of Zhu \cite{Zhu2016} and  Ma-Zhu
\cite{MaZhu2026}.
The same conclusion holds on smooth relatively compact domains in flat
quaternionic quotients $\HH^n/\Gamma$, where $\Gamma$ is any
discrete additive subgroup, as stated in Theorem~\ref{thm:quotient}.

We next record two geometric consequences. Let \(d\) be the interior distance to
\(\bd\), \(\varrho=-d\), and let \(\nu=D\varrho\) on \(\bd\). At a boundary point \(q\),
define the \emph{quaternionic tangent space} and the \emph{quaternionic Levi form} by
\[
 E_q=\{\xi\in\HH^n:e_n(q)^*\xi=0\},
 \qquad R_q=H(\varrho)|_{E_q}.
\]
Here $e_n(q):=\iota^{-1}(\nu(q))\in\HH^n$, where $\iota$ is the real
coordinate identification in \eqref{eq:real-coordinate-identification}. See
Subsection~\ref{sec:boundary-tools}.
The condition \(e_n(q)^*\xi=0\) is quaternionic orthogonality, so \(E_q\) is a right
quaternionic hyperplane. A \emph{defining function} for \(\Omega\)
is a smooth function which is negative in \(\Omega\), zero on \(\bd\), and has nonzero gradient on
\(\bd\). Analogously to the real and complex boundary conditions in
\cite{CNS1985-III,Li2004}, we call \(\bd\)
\emph{strictly quaternionic \((k-1)\)-convex} if
\begin{equation}\label{eq:levi-condition}
 \sigma_j(R_q)>0\quad(1\le j\le k-1,\ q\in\bd).
\end{equation}
Here $\sigma_j(R_q)$ denotes the $j$-th elementary symmetric polynomial
of the $n-1$ real eigenvalues of $R_q$.

\begin{corollary}[Geometric sufficient condition]\label{cor:construction}
Let \(\Omega\Subset\HH^n\) be a smooth bounded domain satisfying
\eqref{eq:levi-condition}. Then there exists a smooth defining function
\(\rho\) for \(\Omega\) such that \(H(\rho)\in\Gamma_k\) on
\(\closure\Omega\). For every smooth boundary value \(\varphi\) and every
smooth positive \(f(q,z)\) with \(f_z\ge0\), a subsolution satisfying
\eqref{eq:subsolution} exists. Consequently, the Dirichlet problem
\eqref{eq:dirichlet} has a unique smooth admissible solution.
\end{corollary}

\begin{corollary}[Constant boundary values]\label{cor:necessity}
Let \(\Omega\Subset\HH^n\) be a smooth bounded domain, let \(\varphi\) be
constant, and let \(f\) satisfy the assumptions of Theorem~\ref{thm:main}.
Then \eqref{eq:dirichlet} has a smooth admissible solution if and only if
\eqref{eq:levi-condition} holds. If \(\bd\) is connected, this condition
is equivalent to
\begin{equation}\label{eq:scalar-levi-condition}
 \sigma_{k-1}(R_q)>0\quad\text{for every }q\in\bd.
\end{equation}
\end{corollary}

The proof combines the subsolution method
\cite{Guan1994,Guan2014,Zhu2016,MaZhu2026} with the algebraic properties of the
quaternionic $k$-Hessian operator. Comparison and differentiation reduce
the global second derivative estimate to the boundary. A barrier
constructed from the subsolution controls the mixed derivatives.
The remaining step is a uniform positive lower bound for the
$(k-1)$-st symmetric function of the tangential quaternionic Hessian.
We obtain this bound using the dual cone and a local auxiliary function,
and then use the equation to control the double normal derivative.
The resulting estimates yield smooth solvability by the continuity method.

Related Calabi--Yau problems on compact quaternionic manifolds have been
studied in \cite{AleskerVerbitsky2010,Alesker2013,AleskerShelukhin2013,BGV2022, DinewSroka2023}.
Gentili--Zhang \cite{GentiliZhang2022} obtained a priori estimates under
a $\mathcal C$-subsolution hypothesis and solved quaternionic Hessian
equations on compact flat hyperk\"ahler manifolds.
Gentili--Vezzoni \cite{GentiliVezzoni2026} established $C^{2,\alpha}$
estimates for a class of fully nonlinear equations on compact
hyperk\"ahler manifolds that includes quaternionic Hessian equations.

The paper is organized as follows.  Section~\ref{sec:prelim} fixes the quaternionic conventions and Hessian
calculus. Section~\ref{sec:hessian-tools} collects the cone properties and
boundary identities needed in the estimates.
Section~\ref{sec:lower} establishes the lower-order, tangential, and mixed
derivative estimates, and Section~\ref{sec:boundary} proves the double
normal estimate. Section~\ref{sec:existence} completes the proof of the
main theorem and extends it to flat quaternionic quotients.
Section~\ref{sec:geometry} proves the geometric corollaries
and gives examples.

\section{Quaternionic linear algebra and Hessian calculus}\label{sec:prelim}
We fix the quaternionic notation and collect the algebraic and
differential identities needed for the subsequent a priori estimates.
Unless otherwise indicated, the quaternionic conventions and standard
linear algebra follow \cite{AleskerLinearAlgebra,Alesker2003} and
\cite[Section~2]{Zhu2016,MaZhu2026}.
For real representations and the four-term averaging of real quadratic
forms, see \cite[Section~1.1]{AleskerShelukhin2013}.
We use the Hessian normalization in \eqref{eq:hessian-definition}.
The ambient space is $\HH^n$. Algebraic definitions and identities are
also used for matrices of size $m$, in particular $m=n-1$ for tangential blocks.

\subsection{Quaternionic linear algebra}\label{sec:linear-algebra}

\begin{definition}[Quaternions]\label{def:quaternions}
The quaternion algebra is
\begin{equation}\label{eq:quaternion-algebra}
 \HH=\{x+y\mathbf i+z\mathbf j+w\mathbf k:x,y,z,w\in\RR\},
 \qquad \mathbf i^2=\mathbf j^2=\mathbf k^2=-1,
 \quad \mathbf i\mathbf j=\mathbf k=-\mathbf j\mathbf i.
\end{equation}
Conjugation is given by
$\overline{x+y\mathbf i+z\mathbf j+w\mathbf k}
=x-y\mathbf i-z\mathbf j-w\mathbf k$.
\end{definition}

For two quaternions
\[
 q_a=x_a+y_a\mathbf i+z_a\mathbf j+w_a\mathbf k,\qquad
 q_b=x_b+y_b\mathbf i+z_b\mathbf j+w_b\mathbf k,
\]
with real coefficients, we have
\begin{equation}\label{eq:quaternion-product-identities}
 \begin{gathered}
 \overline{q_aq_b}=\bar q_b\bar q_a,\qquad |q_aq_b|=|q_a||q_b|,\\
 \Ree(q_aq_b)=x_ax_b-y_ay_b-z_az_b-w_aw_b=\Ree(q_bq_a).
 \end{gathered}
\end{equation}
In particular,
\[
 (1+\mathbf i)(1+\mathbf j)=1+\mathbf i+\mathbf j+\mathbf k,
 \qquad
 (1+\mathbf j)(1+\mathbf i)=1+\mathbf i+\mathbf j-\mathbf k.
\]
For a square quaternionic matrix $C$, write $\trH C=\sum_a c_{aa}$.
For $C\in\HH^{p\times m}$ and $D\in\HH^{m\times p}$,
the products $CD$ and $DC$ are square matrices of sizes $p$ and $m$, respectively, and
\begin{equation}\label{eq:cyclic-real-trace}
 \Ree\trH(CD)
 =\sum_{ a, b}\Ree(c_{ ab}d_{ ba})
 =\Ree\trH(DC).
\end{equation}
In particular, for $C_1,\ldots,C_j\in\HH^{m\times m}$ with $j\ge2$,
\begin{equation}\label{eq:cyclic-real-trace-many}
 \Ree\trH(C_1C_2\cdots C_j)
 =\Ree\trH(C_2\cdots C_jC_1)
 =\Ree\trH(C_jC_1\cdots C_{j-1}).
\end{equation}
\begin{definition}[Real coordinates and inner product]\label{def:real-coordinates}
We regard $\HH^n$ as a right quaternionic vector space of columns
$q=(q_1,\ldots,q_n)^T$ in its standard basis
$\epsilon_1,\ldots,\epsilon_n$.
Unless another basis is specified, quaternionic vectors and matrices
are expressed in this standard basis.
Write $q_a=x_a+y_a\mathbf i+z_a\mathbf j+w_a\mathbf k$.
We identify $\HH^n$ with $\RR^{4n}$ via the real-linear isomorphism
\begin{equation}\label{eq:real-coordinate-identification}
 \iota:\HH^n\longrightarrow\RR^{4n},\qquad
 \iota(q)=(x_1,y_1,z_1,w_1,\ldots,x_n,y_n,z_n,w_n)^T.
\end{equation}
For real column vectors $\xi,\zeta\in\RR^{4n}$, write
$\xi\cdot\zeta:=\xi^T\zeta$.
For $q,\zeta\in\HH^n$, set $q^*:=(\bar q_1,\ldots,\bar q_n)$ and define
\begin{equation}\label{eq:real-inner-product}
 \begin{aligned}
 \langle q,\zeta\rangle_{\RR}
 &:=\Ree(q^*\zeta)=\sum_{a=1}^n\Ree(\bar q_a \zeta_a)\\
 &=\iota(q)^T\iota(\zeta)=\iota(q)\cdot\iota(\zeta),\\
 |q|^2&:=\sum_{a=1}^n|q_a|^2=|\iota(q)|^2.
 \end{aligned}
\end{equation}
\end{definition}

Applying \eqref{eq:quaternion-product-identities} to
$\bar q_a\zeta_a$ and summing over $a$ gives \eqref{eq:real-inner-product}.

The operators $D,D^2,\Delta$ are the real gradient, Hessian and Laplacian.
 Write the components of \eqref{eq:real-coordinate-identification} as
\[
 \iota(q)=(x^1, \ldots,x^{4n})^T,\qquad
 (x^{4a-3},x^{4a-2},x^{4a-1},x^{4a})=(x_a,y_a,z_a,w_a).
 \]
We retain this ordering in every quaternionic orthonormal coordinate system.
Real coordinate indices $r,s,i,j$ range from $1$ to $4n$ unless a
different range is stated. In particular,
\begin{equation}\label{eq:real-coordinate-derivatives}
 \partial_r:=\frac{\partial}{\partial x^r},\qquad
 u_r:=\partial_ru=\frac{\partial u}{\partial x^r},\qquad
 u_{rs}:=\partial_r\partial_su
 =\frac{\partial^2u}{\partial x^r\partial x^s}.
\end{equation}
For real vectors $\xi,\zeta$, set
\begin{equation}\label{eq:directional-derivatives}
 u_\xi:=Du\cdot\xi,\qquad
 u_{\xi\zeta}:=D^2u(\xi,\zeta)
 =\sum_{r,s=1}^{4n}u_{rs}\xi^r\zeta^s.
\end{equation}
For real vector fields $X=\sum_rX^r\partial_r$ and $Y=\sum_sY^s\partial_s$,
\begin{equation}\label{eq:euclidean-covariant-derivative}
 DY=(\partial_sY^r)_{r,s=1}^{4n},\qquad
 D_XY:=(DY)X
 =\sum_{r,s=1}^{4n}X^r(\partial_rY^s)\partial_s.
\end{equation}
Thus
\begin{equation}\label{eq:iterated-directional-derivative}
  D_X(D _Yu)
  =\sum_{r,s=1}^{4n}X^r Y^s u_{rs}
   +\sum_{r,s=1}^{4n}X^r(\partial_r Y^s)u_s
  =D^2u(X, Y)+Du\cdot(D_X Y).
 \end{equation}

\begin{definition}[Matrix norms]\label{def:matrix-norms}
For a real or quaternionic matrix $M\in\mathbb F^{p\times m}$,
where $\mathbb F=\RR$ or $\HH$, and a column vector
$q=(q_1,\ldots,q_m)^T\in\mathbb F^m$, we use
\begin{equation}\label{eq:matrix-norms}
 |M|:=\left(\sum_{i=1}^p\sum_{j=1}^m|M_{ij}|^2\right)^{1/2},\qquad
 \|M\|_{\mathrm{op}}:=\sup_{\substack{q\in\mathbb F^m\\|q|=1}}|Mq|,
 \qquad |q|^2=\sum_{a=1}^m|q_a|^2.
\end{equation}
Unsubscripted matrix norms are Frobenius norms $|M|$.
\end{definition}

\begin{definition}[Quaternionic structures]\label{def:quaternionic-structures}
Let $I,J,K\in\RR^{4n\times4n}$ be determined by
\begin{equation}\label{eq:quaternionic-structures}
 I\iota(q)=\iota(-q\mathbf i),\qquad
 J\iota(q)=\iota(-q\mathbf j),\qquad
 K\iota(q)=\iota(-q\mathbf k)
 \quad(q\in\HH^n).
\end{equation}
In the coordinates \eqref{eq:real-coordinate-identification},
\begin{equation}\label{eq:quaternionic-structure-matrices}
 I=\operatorname{diag}(\underbrace{I_0,\ldots,I_0}_{n}),\qquad
 J=\operatorname{diag}(\underbrace{J_0,\ldots,J_0}_{n}),\qquad
 K=\operatorname{diag}(\underbrace{K_0,\ldots,K_0}_{n}),
\end{equation}
where
\begin{equation}\label{eq:quaternionic-structure-blocks}
 I_0=\begin{pmatrix}
 0&1&0&0\\
 -1&0&0&0\\
 0&0&0&-1\\
 0&0&1&0
 \end{pmatrix},\quad
 J_0=\begin{pmatrix}
 0&0&1&0\\
 0&0&0&1\\
 -1&0&0&0\\
 0&-1&0&0
 \end{pmatrix},\quad
 K_0=\begin{pmatrix}
 0&0&0&1\\
 0&0&-1&0\\
 0&1&0&0\\
 -1&0&0&0
 \end{pmatrix}.
\end{equation}
\end{definition}

Direct matrix multiplication gives
\begin{equation}\label{eq:quaternionic-structure-orthogonality}
 Q^TQ=I_{4n},\qquad Q^T=-Q,\qquad Q^2=-I_{4n}
 \quad(Q\in\{I,J,K\}),
\end{equation}
\begin{equation}\label{eq:quaternionic-structure-relations}
 IJ=K=-JI,\qquad JK=I=-KJ,\qquad KI=J=-IK.
\end{equation}

\begin{definition}[Hyperhermitian and unitary matrices]\label{def:hyperhermitian-unitary}
For a quaternionic matrix $C=(c_{ab})\in\HH^{n\times n}$, set $C^*=(\bar c_{ba})$.
A matrix $A\in\HH^{n\times n}$ is \emph{hyperhermitian} if $A^*=A$. Set
\begin{equation}\label{eq:hyperhermitian-space}
 \mathcal H_n:=\{A\in\HH^{n\times n}:A^*=A\}.
\end{equation}
In particular, $q^*Aq$ is real. We write $A\ge0$ if $q^*Aq\ge0$ for all
$q\in\HH^{n}$, and $A>0$ if the inequality is strict for $q\ne0$.
A matrix $U\in\HH^{n\times n}$ is \emph{quaternionic unitary} if $U^*U=UU^*=I_n$.
\end{definition}

For a quaternionic unitary matrix $U$, write its columns in the standard basis as
\begin{equation}\label{eq:unitary-columns}
 U=(e_1,\ldots,e_n),\qquad e_a=U\epsilon_a,\qquad
 e_a^*e_b=\delta_{ab}\quad(1\le a,b\le n).
\end{equation}
Thus $e_1,\ldots,e_n$ is a quaternionic orthonormal basis.
Every $A\in\mathcal H_n$ has a unitary diagonalization
\begin{equation}\label{eq:unitary-diagonalization}
 U^*AU=\operatorname{diag}(\lambda_1,\ldots,\lambda_n),
 \qquad \lambda _a\in\RR\quad(1\le a\le n).
\end{equation}
We denote these real numbers, counted with multiplicity, by
$\lambda_1(A),\ldots,\lambda_n(A)$ and call them the eigenvalues of $A$.
 For the columns in \eqref{eq:unitary-columns}, this is equivalent to
$Ae_a=e_a\lambda_a(A)$ for $1\le a\le n$.
For a hyperhermitian $m\times m$ matrix $B$, unitary diagonalization gives
\begin{equation}\label{eq:spectral-operator-norm}
 \|B\|_{\mathrm{op}}=\max_{1\le a\le m}|\lambda _a(B)|.
\end{equation}

\begin{definition}[Real representation of a quaternionic matrix]\label{def:real-representation}
For $C=(c_{ab})\in\HH^{n\times n}$, its \emph{real representation}
$\Real(C)$ is determined by
\begin{equation}\label{eq:real-representation-definition}
 \Real:\HH^{n\times n}\longrightarrow\RR^{4n\times4n},\qquad
 \Real(C)\iota(q)=\iota(Cq)\quad(q\in\HH^n).
\end{equation}
\end{definition}

For $c=x+y\mathbf i+z\mathbf j+w\mathbf k$, quaternionic multiplication gives
\[
 \begin{aligned}
 c\mathbf i&=-y+x\mathbf i+w\mathbf j-z\mathbf k,\\
 c\mathbf j&=-z-w\mathbf i+x\mathbf j+y\mathbf k,\\
 c\mathbf k&=-w+z\mathbf i-y\mathbf j+x\mathbf k.
 \end{aligned}
\]
For a single quaternion we use the same coordinate notation:
\[
 q_a=x_a+y_a\mathbf i+z_a\mathbf j+w_a\mathbf k,\qquad
 \iota(q_a)=\begin{pmatrix}x_a\\y_a\\z_a\\w_a\end{pmatrix},
 \qquad x_a,y_a,z_a,w_a\in\RR.
\]
Then
\[
 cq_a=x_ac+y_a(c\mathbf i)+z_a(c\mathbf j)+w_a(c\mathbf k),
\]
and
\begin{equation}\label{eq:left-multiplication-columns}
 \iota(cq_a)=
 \underbrace{\begin{pmatrix}
 \iota(c)&\iota(c\mathbf i)&\iota(c\mathbf j)&\iota(c\mathbf k)
 \end{pmatrix}}_{L_c}
 \begin{pmatrix}x_a\\y_a\\z_a\\w_a\end{pmatrix}
 =L_c\iota(q_a),
\end{equation}
with
\begin{equation}\label{eq:left-multiplication-matrix}
 L_c=\begin{pmatrix}
 x&-y&-z&-w\\ y&x&-w&z\\ z&w&x&-y\\ w&-z&y&x
 \end{pmatrix}.
\end{equation}
Consequently,
\begin{equation}\label{eq:real-representation-blocks}
 \Real(C)=
 \begin{pmatrix}
 L_{c_{11}}&L_{c_{12}}&\cdots&L_{c_{1n}}\\
 L_{c_{21}}&L_{c_{22}}&\cdots&L_{c_{2n}}\\
 \vdots&\vdots&\ddots&\vdots\\
 L_{c_{n1}}&L_{c_{n2}}&\cdots&L_{c_{nn}}
 \end{pmatrix}.
\end{equation}
Thus its $(a,b)$ block is the $4\times4$ real matrix $L_{c_{ab}}$.
We write $\trR$ for the ordinary trace of a real matrix.

\begin{lemma}[Properties of the real representation]\label{lem:real-representation}
The map $\Real:\HH^{n\times n}\to\RR^{4n\times4n}$ is real-linear and
injective. For $C,D\in\HH^{n\times n}$,
\begin{equation}\label{eq:real-representation-algebra}
 \Real(CD)=\Real(C)\Real(D),\qquad
 \Real(C^*)=\Real(C)^T,\qquad
 \Real(I_n)=I_{4n}.
\end{equation}
Its image is characterized by
\begin{equation}\label{eq:real-representation-image}
 \begin{aligned}
 \Real(\HH^{n\times n})
 &=\{M\in\RR^{4n\times4n}:MQ=QM\ \text{for }Q=I,J,K\},\\
 \Real(\mathcal H_n)
 &=\{S\in\operatorname{Sym}(4n,\RR):SQ=QS\ \text{for }Q=I,J,K\}.
 \end{aligned}
\end{equation}
For $A\in\mathcal H_n$ and $q\in\HH^n$,
\begin{equation}\label{eq:real-representation-positivity}
 \begin{gathered}
 \iota(q)^T\Real(A)\iota(q)=q^*Aq,\\
 A\ge0\iff\Real(A)\ge0,\qquad
 A>0\iff\Real(A)>0.
 \end{gathered}
\end{equation}
Moreover,
\begin{equation}\label{eq:real-representation-trace-norm}
 \trR\Real(C)=4\Ree\trH C,\qquad |\Real(C)|^2=4|C|^2.
\end{equation}
If $U$ is quaternionic unitary, then $\Real(U)$ is a real orthogonal matrix, that is,
\begin{equation}\label{eq:unitary-real-orthogonality}
 \Real(U)^T\Real(U)=\Real(U)\Real(U)^T=I_{4n}.
\end{equation}
For $A\in\mathcal H_n$, each eigenvalue, counted with multiplicity,
occurs four times in the spectrum of $\Real(A)$.
\end{lemma}

\begin{proof}
We first prove \eqref{eq:real-representation-algebra}.
Real-linearity follows from \eqref{eq:real-representation-definition}, and
\[
 \Real(CD)\iota(q)=\iota(CDq)
 =\Real(C)\Real(D)\iota(q).
\]
Also $L_{\bar c}=L_c^T$, so transposing the block matrix gives
$\Real(C^*)=\Real(C)^T$. The identity matrix is preserved.
If $\Real(C)=0$, then $Cq=0$ for every $q$, hence $C=0$.

Next, we prove \eqref{eq:real-representation-image}.
For $(\varepsilon,Q)\in\{(\mathbf i,I),(\mathbf j,J),(\mathbf k,K)\}$,
\eqref{eq:quaternionic-structures} and the real-linearity of $\iota$ give
\[
 Q\iota(q)=\iota(-q\varepsilon)=-\iota(q\varepsilon),\qquad
 \iota(q\varepsilon)=-Q\iota(q).
\]
Associativity gives $C(q\varepsilon)=(Cq)\varepsilon$, hence
\[
 \begin{aligned}
 \Real(C)Q\iota(q)
 &=\Real(C)\iota(-q\varepsilon)
 =\iota\bigl(C(-q\varepsilon)\bigr)\\
 &=\iota\bigl(-(Cq)\varepsilon\bigr)
 =Q\iota(Cq)=Q\Real(C)\iota(q).
 \end{aligned}
\]
Since $\iota(\HH^n)=\RR^{4n}$,
\[
 \Real(C)Q=Q\Real(C)\qquad(Q\in\{I,J,K\}).
\]
Conversely, let $M\in\RR^{4n\times4n}$ commute with $I,J,K$ and define
\[
 T:=\iota^{-1}\circ M\circ\iota:\HH^n\longrightarrow\HH^n,
 \qquad T(q)=\iota^{-1}\bigl(M\iota(q)\bigr).
\]
 The map $T$ is real-linear. For the same pairs $(\varepsilon,Q)$,
$MQ=QM$ gives
\[
  T(q\varepsilon) =\iota^{-1} ( -MQ\iota(q) )
 =\iota^{-1} (-QM\iota(q) ) =T(q)\varepsilon.
 \]
 By real-linearity, $T(qc)=T(q)c$ for every $c\in\HH$.
 Using the standard basis, set $C=(T(\epsilon_1),\ldots,T(\epsilon_n))$.
For $q=(q_1,\ldots,q_n)^T=\sum_{a=1}^n \epsilon_aq_a$,
\[
 T(q)=\sum_{a=1}^nT(\epsilon_aq_a)=\sum_{a=1}^nT(\epsilon_a)q_a=Cq.
\]
Since $M\iota(q)=\iota(T(q))=\iota(Cq)$ for every $q\in\HH^n$,
\eqref{eq:real-representation-definition} gives $M=\Real(C)$.
Since $\Real(C^*)=\Real(C)^T$ by \eqref{eq:real-representation-algebra},
the injectivity of $\Real$ gives
\[
 \Real(C)^T=\Real(C)
 \iff\Real(C^*)=\Real(C)
 \iff C^*=C,
\]
which proves the assertion for hyperhermitian matrices.

We now prove \eqref{eq:real-representation-positivity}.
By \eqref{eq:real-representation-definition} and \eqref{eq:real-inner-product},
\[
 \begin{aligned}
 \iota(q)^T\Real(A)\iota(q)
 &=\iota(q)^T\iota(Aq)\\
 &=\iota(q)\cdot\iota(Aq)=\Ree(q^*Aq)=q^*Aq.
 \end{aligned}
\]
Here $q^*Aq\in\RR$ because $A^*=A$. Both positivity equivalences follow.

To prove \eqref{eq:real-representation-trace-norm}, first use
\eqref{eq:left-multiplication-matrix} to obtain
\[
 \trR L_c=4x=4\Ree c,\qquad
 |L_c|^2=4(x^2+y^2+z^2+w^2)=4|c|^2.
\]
The block decomposition \eqref{eq:real-representation-blocks} then gives
\[
 \begin{aligned}
 \trR\Real(C)
 &=\sum_{ a=1}^n\trR L_{c _{aa}}
 =4\sum_{ a=1}^n\Ree c_{ aa}=4\Ree\trH C,\\
 |\Real(C)|^2
 &=\sum_{ a, b=1}^n|L_{c _{ab}}|^2
 =4\sum_{ a, b=1}^n|c_{ ab}|^2=4|C|^2.
 \end{aligned}
\]

Equation \eqref{eq:unitary-real-orthogonality} follows from
\[
 \Real(U)^T\Real(U)=\Real(U^*U)=I_{4n},\qquad
 \Real(U)\Real(U)^T=\Real(UU^*)=I_{4n}.
\]

Finally, choose $U$ as in \eqref{eq:unitary-diagonalization}.
By \eqref{eq:real-representation-algebra} and
\eqref{eq:left-multiplication-matrix}--\eqref{eq:real-representation-blocks},
\[
  \Real(U)^T\Real(A)\Real(U)
  =\Real(U^* AU)
  = \operatorname{diag}(\lambda_1(A) I_4,\ldots,\lambda_n(A)I_4).
 \]
 Since $\Real(U)$ is orthogonal by \eqref{eq:unitary-real-orthogonality},
 each eigenvalue of $A$ occurs in $\Real(A)$ with four times its
quaternionic multiplicity.
\end{proof}

\begin{definition}[Real averaging map]\label{def:real-projection}
Put $\mathscr Q:=\{I_{4n},I,J,K\}$. For a real matrix $S\in\RR^{4n\times 4n}$, define
\begin{equation}\label{eq:projection}
 \begin{aligned}
 \ProjR&:\RR^{4n\times4n}\longrightarrow\RR^{4n\times4n},\\
 \ProjR(S)&:=\frac14\sum_{Q\in\mathscr Q}Q^TSQ=
 \frac14(S+I^TSI+J^TSJ+K^TSK).
 \end{aligned}
\end{equation}
\end{definition}

The following projection identities will be used repeatedly below.

\begin{lemma}
\label{lem:realification-projection}
The map $\ProjR$ is the orthogonal projection of
$\RR^{4n\times4n}$ onto $\Real(\HH^{n\times n})$ with respect to
the Frobenius inner product
\begin{equation}\label{eq:frobenius-inner-product}
 \langle M,N\rangle_{\mathrm F}:=\trR(M^TN),\qquad
 M,N\in\RR^{4n\times4n}.
\end{equation}
For every $S\in\RR^{4n\times4n}$,
\begin{gather}
 \ProjR(S)Q=Q\ProjR(S)\qquad(Q=I,J,K),
 \label{eq:projection-commutation}\\
 \ProjR(S^T)=\ProjR(S)^T,\qquad
 \trR\ProjR(S)=\trR S,\notag\\
 S=S^T\quad\Longrightarrow\quad\ProjR(S)\in\Real(\mathcal H_n),\notag\\
 S=S^T,\quad v^TSv\ge0\ (\forall v\in\RR^{4n})
 \quad\Longrightarrow\quad\ProjR(S)\ge0.\notag
\end{gather}
\end{lemma}

\begin{proof}
By \eqref{eq:quaternionic-structure-orthogonality} and
\eqref{eq:projection},
\[
 4\ProjR(S)=S-ISI-JSJ-KSK.
\]
Using \eqref{eq:quaternionic-structure-relations}, we obtain
\[
 \begin{aligned}
 4\ProjR(S)I&=SI+IS+JSK-KSJ=4I\ProjR(S),\\
 4\ProjR(S)J&=SJ-ISK+JS+KSI=4J\ProjR(S),\\
 4\ProjR(S)K&=SK+ISJ-JSI+KS=4K\ProjR(S).
 \end{aligned}
\]
Transposing \eqref{eq:projection} gives
\begin{equation}\label{eq:projection-transpose}
 \ProjR(S)^T
 =\frac14\bigl(S^T+I^TS^TI+J^TS^TJ+K^TS^TK\bigr)
 =\ProjR(S^T).
\end{equation}
Since $QQ^T=I_{4n}$ for every $Q\in\mathscr Q$,
\[
 \trR\ProjR(S)
 =\frac14\sum_{Q\in\mathscr Q}\trR(Q^TSQ)
 =\frac14\sum_{Q\in\mathscr Q}\trR(SQQ^T)
 =\trR S.
\]

By \eqref{eq:real-representation-image},
\eqref{eq:projection-commutation}, and \eqref{eq:projection-transpose},
\[
 S=S^T\quad\Longrightarrow\quad
 \ProjR(S)=\ProjR(S)^T\in\Real(\mathcal H_n).
\]

If $S=S^T\ge0$ and $v\in\RR^{4n}$, then
\[
 \begin{aligned}
 v^T\ProjR(S)v
 &=\frac14\sum_{Q\in\mathscr Q}v^TQ^TSQv\\
 &=\frac14\sum_{Q\in\mathscr Q}(Qv)^TS(Qv)\ge0.
 \end{aligned}
\]

Finally, \eqref{eq:real-representation-image} and
\eqref{eq:projection-commutation} give
$\ProjR(S)\in\Real(\HH^{n\times n})$.
For $N\in\Real(\HH^{n\times n})$, we have $QNQ^T=N$ for every
$Q\in\mathscr Q$ , so
\[
  \langle\ProjR(S),N\rangle_{\mathrm F}
  =\frac14\sum_{Q\in\mathscr Q}\trR ( S^T QNQ^T)
  = \trR(S^T N)=\langle S,N\rangle_{\mathrm F}.
 \]
 Thus $S-\ProjR(S)$ is orthogonal to $\Real(\HH^{n\times n})$,
 which proves the projection property.
\end{proof}

\begin{definition}[Quaternionic projection]\label{def:quaternionic-projection}
For a symmetric matrix $S$, define
\begin{equation}\label{eq:quaternionic-projection}
 \Proj:\operatorname{Sym}(4n,\RR)\longrightarrow\mathcal H_n,\qquad
 \Proj(S):=\Real^{-1}(\ProjR(S)).
\end{equation}
Here $\Real^{-1}$ is defined on $\Real(\mathcal H_n)$.
By Lemma~\ref{lem:realification-projection},
$\ProjR(S)\in\Real(\mathcal H_n)$, and $\Real$ is injective by
Lemma~\ref{lem:real-representation}. Thus $\Proj(S)$ is the unique
$A\in\mathcal H_n$ satisfying $\Real(A)=\ProjR(S)$.
\end{definition}

The map $\Proj$ is surjective: for every $A\in\mathcal H_n$,
Lemma~\ref{lem:realification-projection} gives
\begin{equation}\label{eq:quaternionic-projection-right-inverse}
 \Proj(\Real(A))
 =\Real^{-1}\bigl(\ProjR(\Real(A))\bigr)
 =\Real^{-1}(\Real(A))=A.
\end{equation}
 \subsection{Quaternionic derivatives and elementary Hessian calculus}
\begin{definition}[Quaternionic differential operators]\label{def:quaternionic-derivatives}
Let $\Omega\subset\HH^n$ be open. For $u\in C^1(\Omega,\HH)$, define
\begin{equation}\label{eq:quaternionic-differential-operators}
 \begin{aligned}
 u_{\bar a}:=\frac{\partial u}{\partial\bar q_a}
 &=u_{x_a}+\mathbf i u_{y_a}+\mathbf j u_{z_a}+\mathbf k u_{w_a},\\
 u_a:=\frac{\partial u}{\partial q_a}
 &=\overline{\frac{\partial\bar u}{\partial\bar q_a}}
 =u_{x_a}-u_{y_a}\mathbf i-u_{z_a}\mathbf j-u_{w_a}\mathbf k,
 \end{aligned}
\end{equation}
where
\[
 u_{x_a}:=\frac{\partial u}{\partial x_a},\qquad
 u_{y_a}:=\frac{\partial u}{\partial y_a},\qquad
 u_{z_a}:=\frac{\partial u}{\partial z_a},\qquad
 u_{w_a}:=\frac{\partial u}{\partial w_a}.
\]
\end{definition}
Subscripts $a,\bar a$ abbreviate quaternionic derivatives. Real coordinate
names are retained, as in $u_{x_a}$ and $u_{x_ay_b}$.
We also use $r,s$ for real coordinate directions as in
\eqref{eq:real-coordinate-derivatives}.
For real $u$, we have $u_a=\overline{u_{\bar a}}$.
For $u\in C^2(\Omega,\HH)$, commutation of real partial derivatives and
of left and right multiplication gives
\begin{equation}\label{eq:quaternionic-derivative-commutation}
 \frac{\partial}{\partial\bar q_a}\left(\frac{\partial u}{\partial q_b}\right)
 =\frac{\partial}{\partial q_b}\left(\frac{\partial u}{\partial\bar q_a}\right).
\end{equation}

\begin{definition}[Normalized quaternionic Hessian]\label{def:quaternionic-hessian}
For a real function $u\in C^2(\Omega)$, define the quaternionic Hessian by
\begin{equation}\label{eq:hessian-definition}
 H(u):=\frac14(u_{\bar a b})_{a,b=1}^n,\qquad
 u_{\bar a b}:=\frac{\partial^2u}{\partial\bar q_a\,\partial q_b}
 =\frac{\partial}{\partial\bar q_a}\left(\frac{\partial u}{\partial q_b}\right).
\end{equation}
We write $H(u)(q)$ for its value at $q$. In real coordinates,
\begin{equation}\label{eq:hessian-component-expansion}
 \begin{aligned}
 u_{\bar a b}
 &=u_{x_ax_b}+u_{y_ay_b}+u_{z_az_b}+u_{w_aw_b}\\
 &\quad+\mathbf i(u_{y_ax_b}-u_{x_ay_b}-u_{z_aw_b}+u_{w_az_b})\\
 &\quad+\mathbf j(u_{z_ax_b}-u_{x_az_b}+u_{y_aw_b}-u_{w_ay_b})\\
 &\quad+\mathbf k(u_{w_ax_b}-u_{x_aw_b}-u_{y_az_b}+u_{z_ay_b}).
 \end{aligned}
\end{equation}
\end{definition}

 \begin{lemma}[Hyperhermitian symmetry, trace, and real Hessians]\label{lem:hessian-projection}
For every real function $u\in C^2(\Omega)$, the matrix $H(u)$ is
hyperhermitian and
\[
 H(u)_{aa}=\frac14(u_{x_ax_a}+u_{y_ay_a}+u_{z_az_a}+u_{w_aw_a}),
 \qquad \trH(H(u))=\frac14\Delta u.
\]
Moreover,
\begin{equation}\label{eq:hessian-projection}
 H(u)=\Proj(D^2u),\qquad
 \Real(H(u))=\ProjR(D^2u).
\end{equation}
\end{lemma}

\begin{proof}
Since $u$ is real, \eqref{eq:quaternionic-differential-operators} and
\eqref{eq:quaternionic-derivative-commutation} give
\[
 \begin{aligned}
 \overline{H(u)_{ab}}
 &=\frac14\overline{\frac{\partial}{\partial\bar q_a}
       \left(\frac{\partial u}{\partial q_b}\right)}
 =\frac14\frac{\partial}{\partial q_a}
       \left(\frac{\partial u}{\partial\bar q_b}\right)\\
 &=\frac14\frac{\partial}{\partial\bar q_b}
       \left(\frac{\partial u}{\partial q_a}\right)
 =H(u)_{ba}.
 \end{aligned}
\]
Setting $b=a$ in \eqref{eq:hessian-component-expansion} gives the diagonal
formula, and summing over $a$ gives the trace identity.

For fixed $a,b$, write the corresponding $4\times4$ block of $D^2u$ as
\[
 S_{ab}=
 \begin{pmatrix}
 u_{x_ax_b}&u_{x_ay_b}&u_{x_az_b}&u_{x_aw_b}\\
 u_{y_ax_b}&u_{y_ay_b}&u_{y_az_b}&u_{y_aw_b}\\
 u_{z_ax_b}&u_{z_ay_b}&u_{z_az_b}&u_{z_aw_b}\\
 u_{w_ax_b}&u_{w_ay_b}&u_{w_az_b}&u_{w_aw_b}
 \end{pmatrix}
 =(s_{ij})_{i,j=1}^4.
\]
Using \eqref{eq:quaternionic-structure-blocks}, we obtain
\[
 I_0^TS_{ab}I_0=
 \begin{pmatrix}
 s_{22}&-s_{21}&-s_{24}&s_{23}\\
 -s_{12}&s_{11}&s_{14}&-s_{13}\\
 -s_{42}&s_{41}&s_{44}&-s_{43}\\
 s_{32}&-s_{31}&-s_{34}&s_{33}
 \end{pmatrix}.
\]
The products with $J_0$ and $K_0$ are computed in the same way. Set
\[
 \begin{aligned}
 h_{ab}^0&=s_{11}+s_{22}+s_{33}+s_{44},\\
 h_{ab}^1&=s_{21}-s_{12}-s_{34}+s_{43},\\
 h_{ab}^2&=s_{31}-s_{13}+s_{24}-s_{42},\\
 h_{ab}^3&=s_{41}-s_{14}-s_{23}+s_{32}.
 \end{aligned}
\]
By \eqref{eq:hessian-component-expansion},
\[
 u_{\bar a b}=h_{ab}^0+h_{ab}^1\mathbf i+h_{ab}^2\mathbf j+h_{ab}^3\mathbf k.
\]
Adding the four matrices and using \eqref{eq:left-multiplication-matrix} gives
\[
 \begin{aligned}
 &S_{ab}+I_0^TS_{ab}I_0+J_0^TS_{ab}J_0+K_0^TS_{ab}K_0\\
 &\qquad=
 \begin{pmatrix}
 h_{ab}^0&-h_{ab}^1&-h_{ab}^2&-h_{ab}^3\\[3pt]
 h_{ab}^1&h_{ab}^0&-h_{ab}^3&h_{ab}^2\\[3pt]
 h_{ab}^2&h_{ab}^3&h_{ab}^0&-h_{ab}^1\\[3pt]
 h_{ab}^3&-h_{ab}^2&h_{ab}^1&h_{ab}^0
 \end{pmatrix}
 =L_{u_{\bar a b}}.
 \end{aligned}
\]
Consequently,
\[
 \ProjR(D^2u)
 =\Real\left(\frac14(u_{\bar a b})_{a,b=1}^n\right)
 =\Real(H(u)).
\]
The injectivity of $\Real$ proves \eqref{eq:hessian-projection}.
\end{proof}

In particular,
\begin{equation}\label{eq:normalization}
 H\left(\frac{|q|^2}{2}\right)=I_n,
 \qquad \Real\left(H\left(\frac{|q|^2}{2}\right)\right)=I_{4n}.
\end{equation}

\begin{lemma}[Averages, rank-one matrices, and scalar calculus]
\label{lem:hessian-calculus}
For $S\in\operatorname{Sym}(4n,\RR)$ and $\xi\in\HH^n$,
\begin{equation}\label{eq:quadratic-average}
 \xi^*\Proj(S)\xi
 =\frac14\sum_{Q\in\mathscr Q}
       \bigl(Q\iota(\xi)\bigr)^TS\bigl(Q\iota(\xi)\bigr).
\end{equation}
For real column vectors $a,b\in\RR^{4n}$, write
$a\otimes b:=ab^T=(a_rb_s)_{r,s=1}^{4n}$.
Set $\widehat a:=\iota^{-1}(a)\in\HH^n$. Then
\begin{equation}\label{eq:projected-rank-one}
 \Proj(a\otimes a)=\frac14\widehat a\,\widehat a^*.
\end{equation}
For real smooth functions $g,h$ and $\chi:\RR\to\RR$,
\begin{align}
 H(gh)&=gH(h)+hH(g)
       +\Proj(Dg\otimes Dh+Dh\otimes Dg),\label{eq:hessian-product}\\
 H(\chi(g))&=\chi'(g)H(g)
       +\chi''(g)\Proj(Dg\otimes Dg).\label{eq:hessian-chain}
\end{align}
Here $Dg$ and $Dh$ are real gradient columns.
\end{lemma}
\begin{proof}
By Lemma~\ref{lem:real-representation} and \eqref{eq:projection},
\[
 \begin{aligned}
 \xi^*\Proj(S)\xi
 &=\iota(\xi)^T\Real(\Proj(S))\iota(\xi)\\
 &=\iota(\xi)^T\ProjR(S)\iota(\xi)\\
 &=\frac14\sum_{Q\in\mathscr Q}
       \bigl(Q\iota(\xi)\bigr)^TS\bigl(Q\iota(\xi)\bigr),
 \end{aligned}
\]
which proves \eqref{eq:quadratic-average}.
For $S=a\otimes a$, write
\[
 c:=\widehat a^*\xi=x+y\mathbf i+z\mathbf j+w\mathbf k,
 \qquad x,y,z,w\in\RR.
\]
Since $a=\iota(\widehat a)$, \eqref{eq:quaternionic-structures} and
\eqref{eq:real-inner-product} give
\[
 \begin{aligned}
 a^T\iota(\xi)
 &=\Ree(\widehat a^*\xi)=\Ree c=x,\\
 a^TI\iota(\xi)
 &=a^T\iota(-\xi\mathbf i)
 =\Ree\bigl(\widehat a^*(-\xi\mathbf i)\bigr)
 =-\Ree(c\mathbf i)=y,\\
 a^TJ\iota(\xi)
 &=a^T\iota(-\xi\mathbf j)
 =\Ree\bigl(\widehat a^*(-\xi\mathbf j)\bigr)
 =-\Ree(c\mathbf j)=z,\\
 a^TK\iota(\xi)
 &=a^T\iota(-\xi\mathbf k)
 =\Ree\bigl(\widehat a^*(-\xi\mathbf k)\bigr)
 =-\Ree(c\mathbf k)=w.
 \end{aligned}
\]
Moreover, $a\otimes a=aa^T$, so for each $Q\in\mathscr Q$,
\[
 \begin{aligned}
 \bigl(Q\iota(\xi)\bigr)^T(a\otimes a)\bigl(Q\iota(\xi)\bigr)
 &=\bigl((Q\iota(\xi))^Ta\bigr)\bigl(a^TQ\iota(\xi)\bigr)\\
 &=\bigl(a^TQ\iota(\xi)\bigr)^2,
 \end{aligned}
\]
since $(Q\iota(\xi))^Ta=a^TQ\iota(\xi)\in\RR$.
Applying \eqref{eq:quadratic-average} therefore yields
\[
 \begin{aligned}
 \xi^*\Proj(a\otimes a)\xi
 &=\frac14\sum_{Q\in\mathscr Q}
       \bigl(Q\iota(\xi)\bigr)^T(a\otimes a)\bigl(Q\iota(\xi)\bigr)\\
 &=\frac14\sum_{Q\in\mathscr Q}
       \bigl(a^TQ\iota(\xi)\bigr)^2\\
 &=\frac14(x^2+y^2+z^2+w^2)=\frac14\bar c c\\
 &=\frac14(\xi^*\widehat a)(\widehat a^*\xi)
 =\frac14\xi^*\widehat a\,\widehat a^*\xi.
 \end{aligned}
\]
Since this holds for every $\xi\in\HH^n$, the two hyperhermitian
matrices coincide, proving \eqref{eq:projected-rank-one}.
Finally,
\[
 (gh)_{rs}=g h_{rs}+h g_{rs}+g_rh_s+h_rg_s,\qquad
 (\chi(g))_{rs}=\chi'(g)g_{rs}+\chi''(g)g_rg_s.
\]
Applying the real-linear map $\Proj$ gives the last two identities.
\end{proof}

\begin{definition}[Quaternionic plurisubharmonicity]
\label{def:quaternionic-plurisubharmonicity}
For $p,\xi\in\HH^n$ with $\xi\ne0$, a right quaternionic affine line is
\begin{equation}\label{eq:quaternionic-affine-line}
 p+\xi\HH:=\{p+\xi a:a\in\HH\}.
\end{equation}
A real upper semicontinuous function $u$ on $\Omega$ is
\emph{quaternionic plurisubharmonic} if its restriction to the
intersection of $\Omega$ with every such line is subharmonic.
A real $C^2$ function is \emph{strictly quaternionic plurisubharmonic}
if these restrictions have positive Laplacian.
\end{definition}

\begin{proposition}[Hessian criterion]
\label{prop:psh-hessian-criterion}
For real $u\in C^2(\Omega)$,
\begin{equation}\label{eq:psh-hessian-equivalence}
 \begin{aligned}
 u\text{ is quaternionic plurisubharmonic}
 &\quad\Longleftrightarrow\quad H(u)\ge0,\\
 u\text{ is strictly quaternionic plurisubharmonic}
 &\quad\Longleftrightarrow\quad H(u)>0.
 \end{aligned}
\end{equation}
\end{proposition}

\begin{lemma}[Affine covariance]\label{lem:unitary-change}
Let $U\in\HH^{n\times n}$ satisfy $U^*U=I_n$, and let $q_0\in\HH^n$.
For a real $C^2$ function $u$, under the affine change
$q=q_0+U\widetilde q$, define
$\widetilde u(\widetilde q)=u(q_0+U\widetilde q)$. Then
\begin{equation}\label{eq:unitary-covariance}
 \begin{aligned}
 \widetilde u_{\bar a b}(\widetilde q)
 &=\sum_{c,d=1}^n\overline{U_{ca}}\,
       u_{\bar c d}(q_0+U\widetilde q)U_{db},\\
 H(\widetilde u)(\widetilde q)
 &=U^*H(u)(q_0+U\widetilde q)U.
 \end{aligned}
\end{equation}
The coordinate change preserves Hessian positivity and the Moore
determinant at corresponding points.
\end{lemma}

\begin{definition}[Real trace pairing]\label{def:real-trace-pairing}
We use the positive definite real pairing on $\mathcal H_n$
\begin{equation}\label{eq:real-pairing}
 \pair{A}{B}:=\Ree\trH(AB)
 =\frac14\trR(\Real(A)\Real(B)).
\end{equation}
\end{definition}

The equality follows from the product and trace identities in
Lemma~\ref{lem:real-representation}.

\begin{lemma}[Orthonormal bases and matrix representations]
\label{lem:basis-change}
Let $U=(e_1,\ldots,e_n)$ be quaternionic unitary as in
 \eqref{eq:unitary-columns}, and put $e=(e_1,\ldots,e_n)$.
Let $\widehat q=(\widehat q_1,\ldots,\widehat q_n)^T$ be the coordinate
column of $q\in\HH^n$ in the basis $e$, and let $[A]_e$ be the matrix
of $A\in\mathcal H_n$ in this basis. Then
\begin{equation}\label{eq:orthonormal-basis-change}
 q=U\widehat q,\qquad \widehat q=U^*q,\qquad
 [A]_e=U^*AU=(e_a^*Ae_b)_{a,b=1}^n.
\end{equation}
The corresponding real transformations are
\begin{equation}\label{eq:real-basis-change}
 \iota(\widehat q)=\Real(U)^T\iota(q),\qquad
 \Real([A]_e)=\Real(U)^T\Real(A)\Real(U).
\end{equation}
For $1\le m\le n$, let
$E=\operatorname{span}_{\HH}\{e_1,\ldots,e_m\}\subset\HH^n$ and
$V=(e_1,\ldots,e_m)\in\HH^{n\times m}$.
We write $A|_E$ for the restriction of the hyperhermitian form
$(q,\zeta)\mapsto q^*A\zeta$ to $E$.
Its matrix in the basis $e_1,\ldots,e_m$ is
\begin{equation}\label{eq:subspace-matrix}
 [A|_E]=V^*AV=(e_a^*Ae_b)_{a,b=1}^m,\qquad V^*V=I_m.
\end{equation}
Thus $V:\HH^m\to E$ is an isometry. For $q\in\HH^n$, $V^*q$ is
the coordinate column of its orthogonal projection $VV^*q$ onto $E$.

For $A,B\in\mathcal H_n$ and $q\in\HH^n$,
\begin{equation}\label{eq:pairing-basis-invariance}
 \pair{U^*AU}{U^*BU}=\pair{A}{B},\qquad
 \pair{A}{qq^*}=q^*Aq.
\end{equation}
The pairing on $E$ is therefore defined using any orthonormal basis.
If $\Theta:E\to E$ is hyperhermitian and
$\Theta e_a=e_a\mu_a$ with $\mu_a\in\RR$, then
\begin{equation}\label{eq:pairing-in-eigenbasis}
 \pair{\Theta}{A|_E}
 =\Ree\trH\bigl(\operatorname{diag}(\mu_1,\ldots,\mu_m)V^*AV\bigr)
 =\sum_{a=1}^m\mu_a e_a^*Ae_a.
\end{equation}
\end{lemma}
\begin{proof}
Since the $a$-th column of $U$ is $e_a$, expansion in the standard basis gives
\[
 e_a=\sum_{b=1}^n\epsilon_b U_{ba},\qquad
 q=\sum_{a=1}^n e_a\widehat q_a
   =\sum_{b=1}^n\epsilon_b
      \left(\sum_{a=1}^n U_{ba}\widehat q_a\right).
\]
Comparing with $q=\sum_{b=1}^n\epsilon_b q_b$ gives
$q_b=\sum_{a=1}^nU_{ba}\widehat q_a$, hence $q=U\widehat q$.
Since $U^*U=I_n$, we obtain $\widehat q=U^*q$.
For $\zeta=U\widehat\zeta$ and $\eta,\xi\in\HH^m$,
\[
 q^*A\zeta=\widehat q^*(U^*AU)\widehat\zeta,\qquad
 (V\eta)^*A(V\xi)=\eta^*(V^*AV)\xi.
\]
The first identity gives the matrix formula in
\eqref{eq:orthonormal-basis-change}, and the second gives
\eqref{eq:subspace-matrix}. Equation \eqref{eq:real-basis-change}
follows from \eqref{eq:real-representation-definition} and
\eqref{eq:real-representation-algebra}.
Equations \eqref{eq:cyclic-real-trace} and \eqref{eq:real-pairing}
give \eqref{eq:pairing-basis-invariance}.
Finally, the diagonal entries of $V^*AV$ are the real numbers
$e_a^*Ae_a$, proving \eqref{eq:pairing-in-eigenbasis}.
\end{proof}

\subsection{The Moore determinant}\label{sec:moore-determinant}
\begin{proposition}[Polynomial characterization]\label{prop:moore-polynomial}
There is a unique real polynomial $P$ on $\mathcal H_n$ satisfying
\begin{equation}\label{eq:moore-polynomial-characterization}
 P(A)^4=\det_{\RR}\Real(A)\quad(A\in\mathcal H_n),
 \qquad P(I_n)=1.
\end{equation}
It is homogeneous of degree $n$ and has integer coefficients.
\end{proposition}

\begin{definition}[Moore determinant]\label{def:moore-determinant}
The Moore determinant of $A\in\mathcal H_n$ is
\begin{equation}\label{eq:moore-determinant-definition}
 \det A:=P(A),
\end{equation}
where $P$ is the polynomial in Proposition~\ref{prop:moore-polynomial}.
\end{definition}
Throughout the paper, determinants of hyperhermitian matrices are Moore
determinants.

\begin{lemma}[Moore determinant identities]\label{lem:moore-identities}
The Moore determinant agrees with the usual determinant on real symmetric
and complex Hermitian matrices. For $A\in\mathcal H_n$ and
$C\in\HH^{n\times n}$,
\begin{equation}\label{eq:moore-congruence}
 \det(C^*AC)=\det A\,\det(C^*C).
\end{equation}
In particular,
\begin{equation}\label{eq:moore-unitary-eigenvalues}
 \det(U^*AU)=\det A\quad(U^*U=I_n),\qquad
 \det A=\prod_{a=1}^n\lambda_a(A).
\end{equation}
Thus
\begin{equation}\label{eq:moore-positivity}
 A\ge0\ \Longrightarrow\ \det A\ge0,\qquad
 A>0\ \Longrightarrow\ \det A>0.
\end{equation}
\end{lemma}

\subsection{Mixed discriminants and matrix inequalities}\label{sec:mixed-discriminants}
\begin{definition}[Mixed discriminant]\label{def:mixed-discriminant}
Let $A_1,\ldots,A_n\in\mathcal H_n$. The polynomial
\begin{equation}\label{eq:mixed-discriminant-definition}
 \det(\vartheta_1A_1+\cdots+\vartheta_nA_n),
 \qquad \vartheta_1,\ldots,\vartheta_n\in\RR,
\end{equation}
is homogeneous of degree $n$. The mixed discriminant
$\det(A_1,\ldots,A_n)$ is defined as $1/n!$ times the coefficient of
$\vartheta_1\cdots\vartheta_n$ in \eqref{eq:mixed-discriminant-definition}.
This normalization gives $\det(A,\ldots,A)=\det A$.
For $0\le\ell\le n$, write $A[\ell]$ for $\ell$ copies of $A$, omitting the argument
when $\ell=0$. The total number of matrix arguments is always $n$.
\end{definition}

\begin{lemma}[Properties of mixed discriminants]\label{lem:mixed-discriminants}
The mixed discriminant is symmetric and real-multilinear.
For $A,B,A_1,\ldots,A_n\in\mathcal H_n$, $\vartheta_1,\vartheta_2\in\RR$,
and $1\le k\le n$,
\begin{equation}\label{eq:mixed-binomial}
 \begin{aligned}
 &\det((\vartheta_1A+\vartheta_2B)[k],A_{k+1},\ldots,A_n)\\
 &\qquad=\sum_{ \ell=0}^k \binom{k}{\ell}\vartheta_1^{ k-\ell}\vartheta_2 ^\ell
        \det(B[ \ell],A[ k-\ell],A_{k+1},\ldots,A_n).
 \end{aligned}
\end{equation}
Here $\binom{k}{\ell}:=k!/[\ell!(k-\ell)!]$ is the binomial coefficient
( $0\le \ell\le k$, $0!:=1$), counting the choices of $\ell$ positions
among the first $k$ arguments that contribute $B$.
For every $C\in\HH^{n\times n}$,
\begin{equation}\label{eq:mixed-congruence}
 \det(C^*A_1C,\ldots,C^*A_nC)
 =\det(C^*C)\det(A_1,\ldots,A_n).
\end{equation}
It is nonnegative on nonnegative arguments, and positive on positive
arguments:
\begin{equation}\label{eq:mixed-positivity}
 \begin{aligned}
 A_j\ge0\ (1\le j\le n)&\ \Longrightarrow\
   \det(A_1,\ldots,A_n)\ge0,\\
 A_j>0\ (1\le j\le n)&\ \Longrightarrow\
   \det(A_1,\ldots,A_n)>0.
 \end{aligned}
\end{equation}
It is monotone in each argument when the others are nonnegative. In particular,
\begin{equation}\label{eq:mixed-monotonicity}
 \begin{gathered}
 A\le B,\qquad A_2,\ldots,A_n\ge0\\
 \Longrightarrow\quad
 \det(A,A_2,\ldots,A_n)\le\det(B,A_2,\ldots,A_n).
 \end{gathered}
\end{equation}
\end{lemma}
 \section{The quaternionic \texorpdfstring{$k$}{k}-Hessian operator and boundary tools}
\label{sec:hessian-tools}
Throughout this section, fix integers $n\ge2$ and $2\le k\le n$.
We recall the cone and linearization properties of $\sigma_k$ in
\cite{LiuWang2023,GentiliZhang2022}, and establish the matrix and boundary
identities needed in the a priori estimates.

\subsection{Elementary symmetric functions and admissibility}\label{sec:admissibility}
\leavevmode\par
\begin{definition}[Elementary symmetric functions]\label{def:elementary-symmetric-functions}
For a vector $\boldsymbol\lambda=(\lambda_1,\ldots,\lambda_m)\in\RR^m$, set
\[
 s_j(\boldsymbol\lambda)
 :=\sum_{1\le a_1<\cdots< a_j\le m}\lambda_{ a_1}\cdots\lambda_{ a_j}
 \quad(1\le j\le m),\qquad s_0(\boldsymbol\lambda):=1,
\]
and $s_j(\boldsymbol\lambda):=0$ for $j<0$ or $j>m$.
For $A\in\mathcal H_m$, let
$\boldsymbol\lambda(A)=(\lambda_1(A),\ldots,\lambda_m(A))\in\RR^m$
be its eigenvalue vector and define
\begin{equation}\label{eq:elementary-symmetric-functions}
 \sigma_j(A):=s_j(\boldsymbol\lambda(A)),\qquad 0\le j\le m.
\end{equation}
We set $\sigma_j(A)=0$ for $j>m$.
\end{definition}

For $A\in\mathcal H_n$, \eqref{eq:moore-unitary-eigenvalues} gives
\begin{equation}\label{eq:moore-characteristic-polynomial}
 \det(tI_n+A)=\prod_{ a=1}^n(t+\lambda _a(A))
 =\sum_{j=0}^n\sigma_j(A)t^{n-j},\qquad t\in\RR.
\end{equation}
Applying \eqref{eq:mixed-binomial} to $\det(A+tI_n)$ and comparing
coefficients gives
\begin{equation}\label{eq:sigma-mixed}
 \sigma_j(A)=\binom nj\det(A[j],I_n[n-j]),\qquad 0\le j\le n.
\end{equation}
Thus $\sigma_j$ is a homogeneous polynomial of degree $j$ in the real
coordinates of $A$. Expanding the products in
\eqref{eq:moore-characteristic-polynomial} also gives
\begin{equation}\label{eq:sigma-identity-direction}
 \sigma_j(A+tI_n)
 =\sum_{\ell=0}^j\binom{n-\ell}{j-\ell}\sigma_\ell(A)t^{j-\ell},
 \qquad 0\le j\le n.
\end{equation}

\begin{definition}[G\aa rding cones and admissibility]\label{def:garding-cones}
For $1\le j\le m$, define
\begin{equation}\label{eq:cone-definitions}
 \begin{aligned}
 \gamma_j^{(m)}
   &:=\{\boldsymbol\lambda\in\RR^m:s_\ell(\boldsymbol\lambda)>0,
                         \ 1\le\ell\le j\},\\
 \Gamma_j^{(m)}
   &:=\{A\in\mathcal H_m:\boldsymbol\lambda(A)\in\gamma_j^{(m)}\},
 \qquad \Gamma_j:=\Gamma_j^{(n)}.
 \end{aligned}
\end{equation}
By {\cite[Proposition~2.5 and Corollary~2.1]{LiuWang2023}}, the polynomial
$t\mapsto\sigma_k(A+tI_n)$ has $k$ real zeros, counted with multiplicity,
for every $A\in\mathcal H_n$. All these zeros are negative exactly when
$A\in\Gamma_k$. Thus $\sigma_k$ is hyperbolic in the direction $I_n$,
with hyperbolicity cone $\Gamma_k$.
A real $C^2$ function $u$ is \emph{$k$-admissible} on $\Omega$ if
\begin{equation}\label{eq:admissibility-definition}
 H(u)(q)\in\Gamma_k\qquad(q\in\Omega).
\end{equation}
\end{definition}

Moreover,
\[
 \Gamma_1\supset\cdots\supset\Gamma_n=\{A>0\},\qquad
 \sigma_1(H(u))=\tfrac14\Delta u,\qquad
 \sigma_n(H(u))=\det H(u).
\]

\begin{definition}[Normalized Hessian operator]\label{def:normalized-hessian-operator}
Define $F:\Gamma_k\to(0,\infty)$ and $g$ by
\begin{equation}\label{eq:normalized-hessian-operator}
 F(A)=\binom nk^{-1/k}\sigma_k(A)^{1/k},
 \qquad g(q,z)=\binom nk^{-1/k}f(q,z)^{1/k},
\end{equation}
where $f$ is the prescribed right-hand side of the Dirichlet problem
\eqref{eq:dirichlet}.
\end{definition}

\begin{definition}[Admissible subsolution]\label{def:admissible-subsolution}
For the data $f$ and $\varphi$ in \eqref{eq:dirichlet}, an
\emph{admissible subsolution} is a real function
$\ubar\in C^2(\overline\Omega)$ satisfying
\begin{equation}\label{eq:admissible-subsolution-normalized}
 H(\ubar)\in\Gamma_k\text{ in }\Omega,\qquad
 F(H(\ubar))\ge g(q,\ubar)\text{ in }\Omega,\qquad
 \ubar|_{\bd}=\varphi.
\end{equation}
\end{definition}

\subsection{Linearization and concavity}\label{sec:linearization}
For $A\in\Gamma_k$ and $Y\in\mathcal H_n$, define
\begin{equation}\label{eq:matrix-differential}
 DF(A)[Y]:=\left.\frac{\mathrm d}{\mathrm dt}\right|_{t=0}F(A+tY),
 \qquad t\in\RR.
\end{equation}
This differential is real-linear in $Y$.

\begin{definition}[Newton tensors]\label{def:newton-tensors}
For $0\le j\le m-1$, define
$T_j:\mathcal H_m\to\mathcal H_m$ by
\begin{equation}\label{eq:newton-tensor-definition}
 T_j(A):=\sum_{\ell=0}^j(-1)^\ell\sigma_{j-\ell}(A)A^\ell,
 \qquad T_0(A)=I_m.
\end{equation}
We call $T_j(A)$ the $j$-th Newton tensor (or Newton transformation) of $A$.
\end{definition}

Let $A=U\Lambda U^*$, where $U$ is unitary and
$\Lambda=\operatorname{diag}(\boldsymbol\lambda)$ with
$\boldsymbol\lambda=\boldsymbol\lambda(A)\in\RR^m$.
Write $\boldsymbol\lambda|a\in\RR^{m-1}$ for the vector obtained by
omitting $\lambda_a$. Separating the terms containing $\lambda_a$ gives
\[
 s _\ell(\boldsymbol\lambda)
 =s _\ell(\boldsymbol\lambda|a)+\lambda_a s_{ \ell-1}(\boldsymbol\lambda|a)
 \qquad( \ell\ge0).
\]
Consequently,
\[
 \sum_{\ell=0}^j(-\lambda_a)^\ell s_{j-\ell}(\boldsymbol\lambda)
 =\sum_{\ell=0}^j(-\lambda_a)^\ell s_{j-\ell}(\boldsymbol\lambda|a)
   -\sum_{\ell=1}^j(-\lambda_a)^\ell s_{j-\ell}(\boldsymbol\lambda|a)
 =s_j(\boldsymbol\lambda|a).
\]
Since $A^\ell=U\Lambda^\ell U^*$, it follows that
\begin{equation}\label{eq:newton-spectral-formula}
 T_j(A)=U\operatorname{diag}\bigl(s_j(\boldsymbol\lambda|1),\ldots,
                              s_j(\boldsymbol\lambda|m)\bigr)U^*.
\end{equation}
For $A,Y\in\mathcal H_n$, take $m=n$ above and set
$\widehat Y=U^*YU$ as in \eqref{eq:orthonormal-basis-change}.
The first-variation formula for quaternionic spectral
functions \cite[Section~3]{GentiliZhang2022},
\eqref{eq:newton-spectral-formula}, and \eqref{eq:pairing-in-eigenbasis} give
\[
 D\sigma_k(A)[Y]
 =\sum_{a=1}^n s_{k-1}(\boldsymbol\lambda|a)\widehat Y_{aa}
 =\Ree\trH\bigl(T_{k-1}(\Lambda)\widehat Y\bigr)
 =\pair{T_{k-1}(A)}{Y}.
\]
At repeated eigenvalues, the formula follows by continuity.
The same differential has a mixed-discriminant expression. By
\eqref{eq:sigma-mixed} and \eqref{eq:mixed-binomial},
\[
 \sigma_k(A+tY)
 =\binom nk\det((A+tY)[k],I_n[n-k])
 =\binom nk\sum_{ \ell=0}^k \binom{k}{\ell}t ^\ell
       \det(Y[ \ell],A[ k-\ell],I_n[n-k]).
\]
Taking the coefficient of $t$ yields
\begin{equation}\label{eq:newton-mixed}
 D\sigma_k(A)[Y]
 =\pair{T_{k-1}(A)}{Y}
 =k\binom nk\det(Y,A[k-1],I_n[n-k]).
\end{equation}
 For $A\in\Gamma_k$ and $Y\in\mathcal H_n$, the chain rule and \eqref{eq:real-pairing} give
\begin{equation}\label{eq:diff}
 \begin{aligned}
 DF(A)[Y]
 &=\frac1k\binom nk^{-1/k}\sigma_k(A)^{1/k-1}
       \Ree\trH(T_{k-1}(A)Y)\\
 &=\frac1{4k}\binom nk^{-1/k}\sigma_k(A)^{1/k-1}
       \trR\bigl(\Real(T_{k-1}(A))\Real(Y)\bigr)\\
 &=\frac1{4k}\binom nk^{-1/k}\sigma_k(A)^{1/k-1}
       \sum_{r,s=1}^{4n}\bigl[\Real(T_{k-1}(A))\bigr]_{rs}
                        \bigl[\Real(Y)\bigr]_{rs}.
 \end{aligned}
\end{equation}

\begin{definition}[Real linearization coefficients]\label{def:real-linearization}
For $A\in\Gamma_k$ and $1\le r,s\le4n$, put
\begin{equation}\label{eq:differential-representation}
 F^{rs}(A)=\frac1{4k}\binom nk^{-1/k}\sigma_k(A)^{1/k-1}
                  \bigl[\Real(T_{k-1}(A))\bigr]_{rs}.
\end{equation}
\end{definition}

The coefficient matrix $(F^{rs}(A))$ belongs to $\Real(\mathcal H_n)$
and is therefore fixed by $\ProjR$.
Using $\Real(\Proj(S))=\ProjR(S)$, \eqref{eq:diff}, and the orthogonal
projection property in Lemma~\ref{lem:realification-projection}, we obtain,
for every $S\in\operatorname{Sym}(4n,\RR)$,
\begin{equation}\label{eq:tool}
 \begin{aligned}
 DF(A)[\Proj(S)]
 &=\sum_{r,s=1}^{4n}F^{rs}(A)\bigl[\ProjR(S)\bigr]_{rs}
  =\sum_{r,s=1}^{4n}F^{rs}(A)S_{rs},\\
 \sum_{r=1}^{4n}F^{rr}(A)&=DF(A)[I_n].
 \end{aligned}
\end{equation}

\begin{lemma}[Concavity and ellipticity]\label{lem:concavity-ellipticity}
The cone $\Gamma_k$ is open and convex. Its closure and boundary satisfy
\begin{equation}\label{eq:cone-closure-boundary}
 \overline{\Gamma}_k=\{A\in\mathcal H_n:\sigma_j(A)\ge0,\ 1\le j\le k\},
 \qquad \partial\Gamma_k=\{A\in\overline{\Gamma}_k:\sigma_k(A)=0\}.
\end{equation}
On $\Gamma_k$, $F$ is positive, concave, and homogeneous of degree one. In particular,
\begin{equation}\label{eq:concavity-second-variation}
 D^2F(A)[Y,Y]\le0
 \qquad(A\in\Gamma_k,\ Y\in\mathcal H_n).
\end{equation}
For $A,Y\in\Gamma_k$,
\begin{equation}\label{eq:supporting-homogeneity}
 DF(A)[Y]\ge F(Y),\qquad DF(A)[A]=F(A).
\end{equation}
In particular,
\begin{equation}\label{eq:concavity-tools}
 DF(A)[I_n]\ge1,\qquad F(X+Y)\ge F(X)+F(Y)
 \quad(X\in\Gamma_k,\ Y\in\overline{\Gamma}_k).
\end{equation}
Here $X+Y\in\Gamma_k$, and $F$ is extended continuously to the closed cone.
Consequently,
\begin{equation}\label{eq:positive-perturbation}
 A\in\Gamma_k,\quad C\ge0
 \quad\Longrightarrow\quad A+C\in\Gamma_k,\quad F(A+C)\ge F(A).
\end{equation}
More generally, for $1\le\ell\le m$,
\begin{equation}\label{eq:newton-positivity}
 A\in\Gamma_\ell^{(m)}\quad\Longrightarrow\quad
 T_j(A)>0\quad(0\le j\le\ell-1).
\end{equation}
The real matrix $(F^{rs}(A))$ is positive definite and commutes with
$I,J,K$ and $\Real(A)$.
\end{lemma}
\begin{proof}
The cone properties follow from {\cite[Proposition~2.2(2),
Proposition~2.5, and Corollary~2.1]{LiuWang2023}}.
 The mixed inequality in {\cite[Corollary~2.2]{LiuWang2023}}, together with
\eqref{eq:sigma-mixed} and \eqref{eq:mixed-binomial}, gives the concavity
of $\sigma_k^{1/k}$ and hence \eqref{eq:concavity-second-variation}.
 Homogeneity and concavity give, respectively,
\[
 DF(A)[A]= F( A) , \qquad
 F( Y) \le F( A) +DF(A)[Y-A]=DF(A)[Y] .
\]
This proves \eqref{eq:supporting-homogeneity}. Taking $Y=I_n$ and using
$F(I_n)=1$ gives the first assertion of \eqref{eq:concavity-tools}.
The second follows from concavity at $(X+Y)/2$ and homogeneity, and
extends to $Y\in\overline\Gamma_k$ by continuity.
Since $\Gamma_k+\overline\Gamma_k\subset\Gamma_k$ and every $C\ge0$
belongs to $\overline\Gamma_k$, \eqref{eq:positive-perturbation} follows.

For $\boldsymbol\lambda\in\gamma_\ell^{(m)}$,
\cite[Proposition~1.1]{CNS1985-III} gives
\begin{equation}\label{eq:symmetric-coordinate-positivity}
 s_j(\boldsymbol\lambda| a)
 = \frac{\partial s_{j+1}}{\partial\lambda_a}(\boldsymbol\lambda)>0
 \qquad(0\le j\le\ell-1,\ 1\le a\le m).
\end{equation}
The spectral formula \eqref{eq:newton-spectral-formula} therefore gives
\eqref{eq:newton-positivity}. By \eqref{eq:differential-representation},
$(F^{rs}(A))$ is a positive multiple of $\Real(T_{k-1}(A))$.
Lemma~\ref{lem:real-representation} gives its positivity and commutation
with $I,J,K$. Commutation with $\Real(A)$ follows because $T_{k-1}(A)$
is a polynomial in $A$.
\end{proof}

\subsection{Hyperplane restrictions and rank-one perturbations}
\label{sec:block-identities}
\begin{lemma}\label{lem:block-identity}
For a hyperhermitian block matrix
\[
 A=\begin{pmatrix}B&\beta\\\beta^*&a_0\end{pmatrix},
 \qquad B\in\mathcal H_{n-1},\quad\beta\in\HH^{n-1},\quad a_0\in\RR,
\]
one has
\begin{equation}\label{eq:block}
 \sigma_k(A)=\sigma_k(B)+a_0\sigma_{k-1}(B)-\beta^*T_{k-2}(B)\beta,
\end{equation}
and, for $s\in\RR$,
\begin{equation}\label{eq:rank-one-sigma}
 \sigma_j(A+sE_{nn})=\sigma_j(A)+s\sigma_{j-1}(B),
 \qquad 1\le j\le n,
\end{equation}
where $E_{nn}$ has its only nonzero entry equal to $1$ in position $(n,n)$.
If $A\in\Gamma_k$, then $B\in\Gamma_{k-1}^{(n-1)}$ and
$A+sE_{nn}\in\Gamma_k$ for $s\ge0$.
Whenever $B\in\Gamma_{k-1}^{(n-1)}$, the last term in \eqref{eq:block}
is nonpositive.
\end{lemma}
\begin{proof}
Choose a unitary $U$ such that
\[
 U^*BU=D:=\operatorname{diag}(b_1,\ldots,b_{n-1}),\qquad
 \mathbf b=(b_1,\ldots,b_{n-1}),\qquad \widehat\beta:=U^*\beta.
\]
By \eqref{eq:orthonormal-basis-change}, with $V=\operatorname{diag}(U,1)$,
\[
 V^*AV=\begin{pmatrix}D& \widehat\beta\\ \widehat\beta^*&a_0\end{pmatrix}.
\]
To compute $\det(tI_n+A)$, take $t$ sufficiently large and set
\[
 B_t:=tI_{n-1}+B>0,\qquad
 B_t^{-1}=U(tI_{n-1}+D)^{-1}U^*.
\]
Thus
\[
 \beta^*B_t^{-1}\beta
 =\widehat\beta^*(tI_{n-1}+D)^{-1}\widehat\beta
 =\sum\limits_{a=1}^{n-1}\frac{|\widehat\beta_a|^2}{t+b_a}.
\]
Since $B_t(B_t^{-1}\beta)=\beta$, the following congruence eliminates
the off-diagonal blocks:
\[
 C=\begin{pmatrix}I_{n-1}&-B_t^{-1}\beta\\0&1\end{pmatrix},
 \qquad
 C^*(tI_n+A)C
 =\begin{pmatrix}B_t&0\\0&t+a_0-\beta^*B_t^{-1}\beta\end{pmatrix}.
\]
Since $\Real(C)$ is block triangular with
identity diagonal blocks,
\[
 \det(C^*C)^4=\bigl(\det_{\RR}\Real(C)\bigr)^2=1,
 \qquad C^*C>0\quad\Longrightarrow\quad\det(C^*C)=1.
\]
Thus \eqref{eq:moore-congruence} gives
\[
 \begin{aligned}
 \det(tI_n+A)
 &=\det B_t\,(t+a_0-\beta^*B_t^{-1}\beta)\\
 &=(t+a_0)\prod\limits_{a}(t+b _a)
       -\sum\limits_{a}| \widehat\beta_a|^2\prod_{\ell\ne a}(t+b_\ell).
 \end{aligned}
\]
This polynomial identity holds for all $t$. By
\eqref{eq:moore-characteristic-polynomial}, the coefficients of $t^{n-k}$ give
\[
 \sigma_k(A)
 =s_k(\mathbf b)+a_0s_{k-1}(\mathbf b)
  -\sum_{ a=1}^{n-1}| \widehat\beta_a|^2s_{k-2}(\mathbf b| a),
\]
where $\mathbf b|a$ omits $b_a$.
By \eqref{eq:newton-spectral-formula},
\[
 \sum_{ a=1}^{n-1}| \widehat\beta_a|^2s_{k-2}(\mathbf b| a)
 = \widehat\beta^*U^*T_{k-2}(B)U \widehat\beta
 =\beta^*T_{k-2}(B)\beta.
\]
This proves \eqref{eq:block}.
Replacing $a_0$ by $a_0+s$ in the determinant identity above gives
\[
 \det(tI_n+A+sE_{nn})-\det(tI_n+A)
 =s\det(tI_{n-1}+B).
\]
Comparing coefficients of $t^{n-j}$ using
\eqref{eq:moore-characteristic-polynomial} proves \eqref{eq:rank-one-sigma}.

If $A\in\Gamma_k$ and $2\le j\le k$, differentiating
\eqref{eq:rank-one-sigma} at $s=0$ and then applying
\eqref{eq:newton-mixed} and \eqref{eq:pairing-basis-invariance} gives
\[
 \sigma_{j-1}(B)=D\sigma_j(A)[E_{nn}]
 =\pair{T_{j-1}(A)}{E_{nn}}
 =(T_{j-1}(A))_{nn}>0,
\]
where the strict inequality follows from \eqref{eq:newton-positivity}.
Hence $B\in\Gamma_{k-1}^{(n-1)}$, and
$\beta^*T_{k-2}(B)\beta\ge0$ by Lemma~\ref{lem:concavity-ellipticity}.
The same lemma gives $A+sE_{nn}\in\Gamma_k$ for $s\ge0$.
\end{proof}

\subsection{Dual cones and quantitative estimates}\label{sec:tangential-cone}
Set
\begin{align}
 \mathscr C&:=\overline{\Gamma}_{k-1}^{(n-1)},
 \label{eq:tangential-cone}\\
 \mathscr C^*&:=\{\Theta\in\mathcal H_{n-1}:\pair{\Theta}{X}\ge0
                                      \text{ for all }X\in\mathscr C\},
 \label{eq:tangential-dual-cone}\\
 \mathscr S&:=\{\Theta\in\mathscr C^*:\trH \Theta=1\}.
 \label{eq:normalized-dual-cone}
\end{align}

Since $\mathscr C$ contains every nonnegative hyperhermitian matrix,
\eqref{eq:pairing-basis-invariance} gives
\[
 \Theta\in\mathscr C^*\ \Longrightarrow\
 \xi^*\Theta\xi=\pair{\Theta}{\xi\xi^*}\ge0,
 \qquad
 \Theta\ne0\ \Longrightarrow\ \trH \Theta>0.
\]
Hence $\mathscr S$ is a closed subset of the compact set
$\{\Theta\ge0:\trH\Theta=1\}$, and is nonempty since
$I_{n-1}/(n-1)\in\mathscr S$.

For $B\in\Gamma_{k-1}^{(n-1)}$, set
\begin{equation}\label{eq:cone-margin-definition}
 m(B):=\min_{\Theta\in\mathscr S}\pair{\Theta}{B}.
\end{equation}

\begin{lemma}\label{lem:4.1}
For $B\in\Gamma_{k-1}^{(n-1)}$,
\begin{equation}\label{eq:cone-margin-shift}
 m(B)=\max\{s\in\RR:B-sI_{n-1}\in\mathscr C\}>0.
\end{equation}
Moreover,
\begin{equation}\label{eq:boundary-cone}
 B-m(B)I_{n-1}\in\partial\Gamma_{k-1}^{(n-1)},\qquad
 \sigma_{k-1}(B-m(B)I_{n-1})=0.
\end{equation}
For $B_1,B_2\in\Gamma_{k-1}^{(n-1)}$ and unitary $U$,
\begin{equation}\label{eq:margin-continuity}
 |m(B_1)-m(B_2)|\le\|B_1-B_2\|_{\mathrm{op}},\qquad
 m(U^*BU)=m(B).
\end{equation}
\end{lemma}
\begin{proof}
The separation theorem for the closed convex cone $\mathscr C$, followed
by normalization by $\trH\Theta>0$ for nonzero $\Theta\in\mathscr C^*$, gives
\[
 X\in\mathscr C
 \iff\pair{\Theta}{X}\ge0\quad\text{for all }\Theta\in\mathscr S.
\]
Since $\pair{\Theta}{I_{n-1}}=1$ on $\mathscr S$,
\[
 B-sI_{n-1}\in\mathscr C
 \iff s\le\pair{\Theta}{B}\quad\text{for all }\Theta\in\mathscr S
 \iff s\le m(B).
\]
This proves the maximum formula in \eqref{eq:cone-margin-shift}.
Openness gives $B-\epsilon I_{n-1}\in\mathscr C$ for some $\epsilon>0$,
so $m(B)>0$. Maximality places $B-m(B)I_{n-1}$ on the boundary.
Equation~\eqref{eq:cone-closure-boundary} then gives \eqref{eq:boundary-cone}.

Choose $\Theta_i\in\mathscr S$ with $m(B_i)=\pair{\Theta_i}{B_i}$.
Since $\Theta_i\ge0$ and $\trH\Theta_i=1$,
\eqref{eq:pairing-in-eigenbasis} gives
\[
 m(B_i)-m(B_j)
 \le\pair{\Theta_j}{B_i-B_j}
 \le\|B_i-B_j\|_{\mathrm{op}},
 \qquad(i,j)=(1,2),(2,1).
\]
Finally, unitary invariance of $\mathscr C$ and
\eqref{eq:pairing-basis-invariance}, together with
\eqref{eq:cyclic-real-trace}, imply
\[
 \pair{U\Theta U^*}{X}=\pair{\Theta}{U^*XU},\qquad
 \trH(U\Theta U^*)=\trH\Theta,
 \qquad U^*\mathscr C U=\mathscr C.
\]
Consequently $\Theta\in\mathscr S\iff U\Theta U^*\in\mathscr S$, and
\[
 m(U^*BU)
 =\min_{\Theta\in\mathscr S}\pair{U\Theta U^*}{B}
 =m(B).
\]
\end{proof}

In particular, when $k=2$ or $k=n$, the corresponding closed cones are
$\{\trH B\ge0\}$ and $\{B\ge0\}$, respectively. Thus
\eqref{eq:boundary-cone} gives, respectively,
\[
 \trH(B-m(B)I_{n-1})=0,\qquad
 \lambda_{\min}(B-m(B)I_{n-1})=0,
\]
and hence
\begin{equation}\label{eq:margin-endpoints}
 k=2:\quad m(B)=\frac{\trH B}{n-1},\qquad
 k=n:\quad m(B)=\lambda_{\min}(B).
\end{equation}

\begin{lemma}[Quantitative cone estimates]\label{lem:margin-bounds}
If $B\in\Gamma_{k-1}^{(n-1)}$ and $\|B\|_{\mathrm{op}}\le M$, then
\begin{equation}\label{eq:margin-bounds}
 \binom{n-1}{k-1}m(B)^{k-1}
 \le \sigma_{k-1}(B)
 \le (k-1)\binom{n-1}{k-1}(2M)^{k-2}m(B).
\end{equation}
Moreover,
\begin{equation}\label{eq:newton-boundary}
 \sigma_k(B)\le
 \frac{n-k}{k(n-1)}\sigma_1(B)\sigma_{k-1}(B).
\end{equation}
\end{lemma}
\begin{proof}
For $k=2$, both bounds in \eqref{eq:margin-bounds} are equalities by
\eqref{eq:margin-endpoints}. For $k\ge3$, we have
$m(B)I_{n-1}\in\Gamma_{k-1}^{(n-1)}$ and
$B-m(B)I_{n-1}\in\partial\Gamma_{k-1}^{(n-1)}$ by
\eqref{eq:boundary-cone}. Applying \eqref{eq:concavity-tools} in dimension
$n-1$ and degree $k-1$ gives
\[
\begin{aligned}
 \sigma_{k-1}(B)^{1/(k-1)}
 &\ge \sigma_{k-1}(m(B)I_{n-1})^{1/(k-1)}
       +\sigma_{k-1}(B-m(B)I_{n-1})^{1/(k-1)}\\
 &=\binom{n-1}{k-1}^{1/(k-1)}m(B),
\end{aligned}
\]
which proves the lower bound. For the upper bound,
$I_{n-1}/(n-1)\in\mathscr S$ and \eqref{eq:cone-margin-definition} give
\begin{equation}\label{eq:margin-norm-bound}
 m(B)\le\frac{\trH B}{n-1}\le\|B\|_{\mathrm{op}}\le M.
\end{equation}
Thus \eqref{eq:margin-norm-bound} implies
\begin{equation}\label{eq:shifted-block-norm-bound}
 \|B-sI_{n-1}\|_{\mathrm{op}}
 \le\|B\|_{\mathrm{op}}+s\le2M,
 \qquad 0\le s\le m(B).
\end{equation}
For each fixed $s$, applying \eqref{eq:sigma-identity-direction} in dimension
$n-1$ and degree $k-1$ gives
\begin{equation}\label{eq:shifted-sigma-derivative}
\begin{aligned}
 \frac{\mathrm d}{\mathrm ds}\sigma_{k-1}(B-sI_{n-1})
 &= -D\sigma_{k-1} (B-sI_{n-1}) [I_{n-1} ]\\
 &=-(n-k+1)\sigma_{k-2}(B-sI_{n-1}).
\end{aligned}
\end{equation}
By \eqref{eq:cone-margin-shift}, $B-sI_{n-1}\in\mathscr C$ on this interval.
Its eigenvalues $\mu_1(s),\ldots,\mu_{n-1}(s)$ satisfy
 $|\mu_a(s)|\le2M$ by \eqref{eq:shifted-block-norm-bound}.
Definition~\ref{def:elementary-symmetric-functions} and
\eqref{eq:elementary-symmetric-functions} therefore give
\begin{equation}\label{eq:shifted-sigma-integrand-bound}
 0\le\sigma_{k-2}(B-sI_{n-1})
 \le\sum_{\substack{E\subset\{1,\ldots,n-1\}\\|E|=k-2}}
       \prod_{ a\in E}|\mu _a(s)|
 \le\binom{n-1}{k-2}(2M)^{k-2}.
\end{equation}
Integrating \eqref{eq:shifted-sigma-derivative} and using
\eqref{eq:boundary-cone} and \eqref{eq:shifted-sigma-integrand-bound}, we obtain
\begin{equation}\label{eq:margin-integrated-upper-bound}
 \begin{aligned}
 \sigma_{k-1}(B)
 &=(n-k+1)\int_0^{m(B)}\sigma_{k-2}(B-sI_{n-1})\,ds\\
 &\le(k-1)\binom{n-1}{k-1}(2M)^{k-2}m(B).
 \end{aligned}
\end{equation}
For \eqref{eq:newton-boundary}, $k=n$ gives $0\le0$, and the case
$\sigma_k(B)\le0$ is immediate. Otherwise $k<n$ and
$B\in\Gamma_k^{(n-1)}$.
For $2\le j\le k$, \eqref{eq:sigma-identity-direction} gives
\begin{equation}\label{eq:identity-direction-first-second-derivatives}
\begin{aligned}
 \left.\frac{\mathrm d}{\mathrm dt}\right|_{t=0}
 \sigma_j(B+tI_{n-1})&=(n-j)\sigma_{j-1}(B),\\
 \left.\frac{\mathrm d^2}{\mathrm dt^2}\right|_{t=0}
 \sigma_j(B+tI_{n-1})&=(n-j)(n-j+1)\sigma_{j-2}(B).
\end{aligned}
\end{equation}
By \eqref{eq:concavity-second-variation}, applied in dimension $n-1$
and degree $j$ with $Y=I_{n-1}$, and
\eqref{eq:identity-direction-first-second-derivatives}, we obtain
\begin{equation}\label{eq:concavity-newton-computation}
\begin{aligned}
 0&\ge\left.\frac{\mathrm d^2}{\mathrm dt^2}\right|_{t=0}
             \sigma_j(B+tI_{n-1})^{1/j}\\
 &=\frac{n-j}{j}\sigma_j(B)^{1/j-2}
 \left[(n-j+1)\sigma_j(B)\sigma_{j-2}(B)
       -\frac{j-1}{j}(n-j)\sigma_{j-1}(B)^2\right].
\end{aligned}
\end{equation}
Since $\sigma_\ell(B)>0$ for $0\le\ell\le k$,
\eqref{eq:concavity-newton-computation} implies
\begin{equation}\label{eq:sigma-ratio-monotonicity}
 \frac{j}{n-j}\frac{\sigma_j(B)}{\sigma_{j-1}(B)}
 \le\frac{j-1}{n-j+1}\frac{\sigma_{j-1}(B)}{\sigma_{j-2}(B)},
 \qquad 2\le j\le k.
\end{equation}
Iterating \eqref{eq:sigma-ratio-monotonicity} and using $\sigma_0(B)=1$ gives
\begin{equation}\label{eq:sigma-ratio-trace-bound}
 \frac{k}{n-k}\frac{\sigma_k(B)}{\sigma_{k-1}(B)}
 \le\cdots\le\frac{2}{n-2}\frac{\sigma_2(B)}{\sigma_1(B)}
 \le\frac{\sigma_1(B)}{n-1},
\end{equation}
which proves \eqref{eq:newton-boundary}.
\end{proof}

 \subsection{Signed distance and boundary Hessians}\label{sec:boundary-tools}
\begin{definition}[Defining function]\label{def:defining-function}
A smooth real function $\rho$ on a neighborhood of $\overline\Omega$
is a \emph{defining function} for $\Omega$ if
\begin{equation}\label{eq:defining-function}
 \rho<0\text{ in }\Omega,\qquad \rho=0\text{ on }\bd,\qquad
 |D\rho|>0\text{ on }\bd.
\end{equation}
\end{definition}

Let $\varrho$ be the signed distance to $\partial\Omega$, defined by
\[
 \varrho(q)=\begin{cases}
 -\operatorname{dist}(q,\partial\Omega),&q\in\Omega,\\
 \operatorname{dist}(q,\partial\Omega),&q\notin\Omega.
 \end{cases}
\]
It is smooth in a fixed tubular neighborhood of $\partial\Omega$
containing points on both sides of the boundary. The interior distance
is $d=-\varrho$ in $\Omega$. Set $\nu=D\varrho$ in this neighborhood.
It is the unit normal to the level sets of $\varrho$, oriented towards
increasing $\varrho$. On $\partial\Omega$, it is the outward unit normal.
At interior points of the neighborhood, $\nu=-Dd$ points towards
$\partial\Omega$.
At every point of the tubular neighborhood, we have
\begin{equation}\label{eq:distance-identities}
 |D\varrho|=1,\qquad D^2\varrho(\nu,X)=0\quad\text{for every real tangent vector }X,
 \qquad D_\nu\nu=0.
\end{equation}
Here $D^2\varrho$ is the real Hessian viewed as a symmetric bilinear form,
and $X$ is any ambient tangent vector at any point of the tubular neighborhood.
Indeed, differentiating $|D\varrho|^2=1$ gives
\[
 0=X(|D\varrho|^2)=2D^2\varrho(\nu,X),\qquad
 \langle D_\nu\nu,X\rangle_{\RR^{4n}}=D^2\varrho(\nu,X)=0.
\]

\Needspace{12\baselineskip}
\begin{definition}[Quaternionic tangent space and Levi form]
\label{def:quaternionic-levi-form}
At $q\in\bd$, keep $\nu(q)=D\varrho(q)\in\RR^{4n}$ as a real vector and set
$e_n(q):=\iota^{-1}(\nu(q))\in\HH^n$ by
\eqref{eq:real-coordinate-identification}. Define
\begin{equation}\label{eq:quaternionic-tangent-levi}
 E_q:=\{\xi\in\HH^n:e_n(q)^*\xi=0\},\qquad
 R_q:=H(\varrho)(q)|_{E_q}.
\end{equation}
The boundary is \emph{strictly quaternionic $(k-1)$-convex} if
\begin{equation}\label{eq:quaternionic-boundary-convexity}
 R_q\in\Gamma_{k-1}^{(n-1)}\qquad(q\in\bd),
\end{equation}
where $R_q$ is represented in any quaternionic orthonormal basis of $E_q$
as in \eqref{eq:subspace-matrix}.
For $k=n$, this means $R_q>0$ and is called
\emph{strict quaternionic pseudoconvexity}.
\end{definition}

At a fixed boundary point $q$, abbreviate $e_n=e_n(q)$, $\nu=\nu(q)$,
and $E=E_q=(e_n\HH)^\perp$. The real vectors
$\nu,I\nu,J\nu,K\nu$ are orthonormal and span $\iota(e_n\HH)$. Hence
\begin{equation}\label{eq:tangent-splitting}
 T_q\partial\Omega=\nu^\perp=\iota(E_q)\oplus_\RR
      \operatorname{span}_\RR\{I\nu,J\nu,K\nu\}.
\end{equation}
Thus $\dim_\RR E=4n-4$ and $\dim_\RR T_q\partial\Omega=4n-1$.

\begin{lemma}[Boundary restriction identities]\label{lem:boundary-restriction}
Let $v$ be smooth and vanish on $\partial\Omega$. Then on $\partial\Omega$,
\begin{align}
 Dv&=v_\nu\nu,\label{eq:boundary-gradient}\\
 D^2v(\tau,\tau')&=v_\nu D^2\varrho(\tau,\tau')
   \quad(\tau,\tau'\in T\partial\Omega),\label{eq:boundary-tangent}\\
 H(v)|_E&=v_\nu H(\varrho)|_E.\label{eq:boundary-restriction}
\end{align}
 For any real $u\in C^2(\overline\Omega)$ and any quaternionic
orthonormal frame $e=(e_1,\ldots,e_n)$ with the last vector $e_n$
defined above, write
$[H(u)]_e=\left(\begin{smallmatrix}B&\beta\\\beta^*&a_0\end{smallmatrix}\right)$
as in \eqref{eq:orthonormal-basis-change}. Then
\begin{equation}\label{eq:normal-entry}
 4a_0=u_{\nu\nu}+D^2u(I\nu,I\nu)+D^2u(J\nu,J\nu)+D^2u(K\nu,K\nu).
\end{equation}
\end{lemma}
\begin{proof}
Since $v=0$ and $D\varrho\ne0$ on $\partial\Omega$, the signed distance
$\varrho$ fixed above gives a local representation $v=\varrho h$ with
$h$ smooth, as in \cite[Sections~4.2--4.3]{MaZhu2026}.
On $\partial\Omega$, $\varrho=0$ and $\varrho_\nu=1$, so
$v_\nu=h$ and $Dv=v_\nu\nu$, proving \eqref{eq:boundary-gradient}.
For arbitrary real vectors $X,Y$ at a boundary point, the product rule gives
\[
 D^2v(X,Y)=v_\nu D^2\varrho(X,Y)
   +(D_X\varrho)(D_Yh)+(D_Xh)(D_Y\varrho).
\]
For $\tau,\tau'\in T\partial\Omega$, we have
$D_\tau\varrho=D_{\tau'}\varrho=0$, so the preceding identity gives
\eqref{eq:boundary-tangent}.
For $\xi\in E_q$ and $Q\in\mathscr Q$,
\eqref{eq:quaternionic-structure-orthogonality} and
\eqref{eq:tangent-splitting} give
\[
 Q^T\nu\in\operatorname{span}_{\RR}\{\nu,I\nu,J\nu,K\nu\}
 \perp\iota(E_q),\qquad
 \langle Q\iota(\xi),\nu\rangle_{\RR^{4n}}
 =\iota(\xi)^TQ^T\nu=0.
\]
Thus $Q\iota(\xi)\in T_q\partial\Omega$. By
\eqref{eq:hessian-projection}, \eqref{eq:quadratic-average}, and
\eqref{eq:boundary-tangent},
\[
 \xi^*H(v)\xi
 =\frac14\sum_Q D^2v(Q\iota(\xi),Q\iota(\xi))
 =\frac{v_\nu}{4}\sum_Q D^2\varrho(Q\iota(\xi),Q\iota(\xi))
 =v_\nu\xi^*H(\varrho)\xi.
\]
Since $\xi\in E_q$ is arbitrary, \eqref{eq:boundary-restriction} follows.
Applying the same averaging formula to $e_n$, with $\iota(e_n)=\nu$, gives
\[
 4a_0=4e_n^*H(u)e_n=\sum_QD^2u(Q\nu,Q\nu).
\]
 By the notation \eqref{eq:directional-derivatives},
 $D^2u(\nu,\nu)=u_{\nu\nu}$, proving \eqref{eq:normal-entry}.
\end{proof}

For the boundary estimates, we first express the off-diagonal
Hessian entries in terms of tangential and mixed derivatives.
For a real $C^2$ function $u$, fix $q\in\partial\Omega$ and choose a
quaternionic orthonormal basis $e_1,\ldots,e_{n-1}$ of $E_q$, completed by
$e_n=e_n(q)$. Equation \eqref{eq:orthonormal-basis-change} gives
\[
 \bigl(e_a^*H(u)(q)e_b\bigr)_{a,b=1}^n
 =\begin{pmatrix}B&\beta\\\beta^*&a_0\end{pmatrix},\qquad
 \beta=(\beta_1,\ldots,\beta_{n-1})^T,\quad
 \beta_a=e_a^*H(u)(q)e_n.
\]
For $1\le a\le n-1$ and $\omega\in\{1,\mathbf i,\mathbf j,\mathbf k\}$,
\[
 (e_a\omega)^*H(u)(q)e_n
 =\bar\omega\,e_a^*H(u)(q)e_n=\bar\omega\beta_a.
\]
Taking real parts and using \eqref{eq:hessian-projection} and the bilinear
identity obtained from \eqref{eq:quadratic-average}, we get, at $q$,
\begin{equation}\label{eq:offdiagonal-real}
 \Ree(\bar\omega\beta_a)
 =\frac14\sum_{Q\in\mathscr Q}
      D^2u\bigl(Q\iota(e_a\omega),Q\nu\bigr).
\end{equation}
Indeed, if $\beta_a=b_0+b_1\mathbf i+b_2\mathbf j+b_3\mathbf k$, then
$\Ree(\bar\omega\beta_a)$ equals $b_0,b_1,b_2,b_3$ for
$\omega=1,\mathbf i,\mathbf j,\mathbf k$, respectively.
Since $e_a\omega\in E_q$, the first vector on the right is tangent to
$\partial\Omega$ at $q$. The second is normal for $Q=I_{4n}$ and tangent
for $Q=I,J,K$. Thus the boundary values of $\beta$ involve only tangential
and mixed tangential--normal second derivatives, not $u_{\nu\nu}$.

We next construct local tangential vector fields.
Fix $p\in\partial\Omega$ and choose a quaternionic unitary matrix
$U=(e_1,\ldots,e_n)$ with $e_n=\iota^{-1}(\nu(p))$.
Use the coordinates $q=p+U\widetilde q$ from
\eqref{eq:unitary-covariance}, suppressing tildes and preserving
the real ordering in \eqref{eq:real-coordinate-identification}.
Then
\[
 p=0,\qquad \nu(0)=\partial_{x_n}=\partial_{4n-3},
 \qquad \varrho_{x_n}(0)=1,
\]
where $\varrho_{x_n}:=\partial_{x_n}\varrho$.
Thus $\varrho_{x_n}\ge1/2$ in a sufficiently small two-sided
neighborhood $\mathcal U$ of $0$. Set
\begin{equation}\label{eq:tangential-index-set}
 \mathscr I:=\{1,\ldots,4n\}\setminus\{4n-3\},\qquad
 \alpha,\alpha'\in\mathscr I.
\end{equation}
Put $\varrho_\alpha:=\partial_\alpha\varrho$ and define
 \begin{equation}\label{eq:tangent-fields}
 X_\alpha=\partial_\alpha- \frac{\varrho_\alpha}{\varrho_{x_n}}\partial_{ x_n},
 \qquad \alpha \in\mathscr I.
\end{equation}
In $\mathcal U$, we have
 $X_\alpha\varrho=\varrho_\alpha-(\varrho_\alpha/\varrho_{x_n})\varrho_{x_n}=0$.
Thus these fields are tangent to the level sets of $\varrho$ in $\mathcal U$.
On $\mathcal U\cap\partial\Omega$, they form a real basis of $T_q\partial\Omega$.
By \eqref{eq:distance-identities} and the constancy of the coordinate
vector fields,
\[
 D_\nu X_\alpha=0\qquad\text{in }\mathcal U,
 \qquad \alpha \in\mathscr I.
\]
Choose finitely many such neighborhoods covering $\partial\Omega$.
Since $\varrho_{x_n}\ge1/2$ on each neighborhood and the derivatives of $\varrho$
up to order three are bounded, there is a common constant $C$, depending
only on the boundary geometry, such that
\[
 \left\| \frac{\varrho_\alpha}{\varrho_{x_n}}\right\|_{C^2(\mathcal U)}\le C,
 \qquad \alpha \in\mathscr I,
\]
on every neighborhood $\mathcal U$ in this cover. Thus all coefficients of
$X_\alpha$ and their real coordinate derivatives of orders one and two
are bounded by a common constant.

\section{Lower-order, tangential, and mixed derivative estimates}\label{sec:lower}
We establish uniform $C^1$ bounds and reduce the $C^2$ estimate to the
double normal derivative on $\bd$. We first record the linearized
identities and estimates needed below. Comparison
and differentiation of the equation give the lower-order estimates and
the reduction to the boundary. A local barrier then bounds the
tangential and mixed derivatives, leaving the double normal estimate
to Section~\ref{sec:boundary}.

Continuity gives $H(\ubar)\in\overline{\Gamma}_k$ and
$\sigma_k(H(\ubar))\ge f(q,\ubar)>0$ on $\closure\Omega$.
Hence \eqref{eq:cone-closure-boundary} implies
$H(\ubar)(\closure\Omega)\Subset\Gamma_k$.
Set $g_0(q)=F(H(\ubar)(q))$ and
\begin{equation}\label{eq:path}
 g_t(q,z)=(1-t)g_0(q)+tg(q,z),\qquad 0\le t\le1.
\end{equation}
We consider smooth admissible solutions of
\begin{equation}\label{eq:equation}
    F(H(u))=g_t(q,u)\quad\text{in }\Omega,
 \qquad u=\varphi\quad\text{on }\bd.
\end{equation}
The function $\ubar$ solves this problem at $t=0$ and is a subsolution
for every $t$, since $g_t(q,\ubar)\le g_0(q)$. We write
$g_{t,z}(q,z):=\partial_z g_t(q,z)$ for the derivative in the second variable $z$. The
indices $r,s$ are reserved for real coordinate derivatives.

Let $h$ be the harmonic extension of $\varphi$ to $\Omega$, and set
\[
 \mathcal I=[\min_{\closure\Omega}\ubar,\max_{\closure\Omega}h].
\]
For $(q,z,t)\in\closure\Omega\times\mathcal I\times[0,1]$, we have
\[
 \sup_{0\le t\le1}\|g_t\|_{C^2(\closure\Omega\times\mathcal I)}\le C,
 \qquad 0<g_*\le g_t(q,z)\le g^*.
\]
Using the Frobenius norm in \eqref{eq:matrix-norms}, set
\[
 \delta_{\ubar}:=\operatorname{dist}
   \bigl(H(\ubar)(\overline\Omega),\partial\Gamma_k\bigr)
 =\inf_{\substack{q\in\overline\Omega\\ A\in\partial\Gamma_k}}
       |H(\ubar)(q)-A|>0.
\]
Throughout Sections~\ref{sec:lower} and~\ref{sec:boundary}, constants $C$ depend only on
\[
 \begin{gathered}
 n,k,\Omega,\quad \|\varphi\|_{C^4(\bd)},\quad
 \|\ubar\|_{C^4(\overline\Omega)},\quad\delta_{\ubar}^{-1},\\
 \|f\|_{C^2(\overline\Omega\times\mathcal I)},\qquad
 \left(\inf_{\overline\Omega\times\mathcal I}f\right)^{-1}.
 \end{gathered}
\]
They may increase from line to line, but are independent of $t$ and $u$.

\subsection{Linearization and preparatory estimates}
\label{sec:linearized-calculus}\label{sec:projection-tool}
Let $u\in C^2(\overline\Omega)$ be an admissible solution of
\eqref{eq:equation}. By continuity,
$H(u)(\overline\Omega)\subset\overline{\Gamma}_k$ and
$\sigma_k(H(u))=\binom nk g_t(q,u)^k>0$ on $\overline\Omega$.
Thus \eqref{eq:cone-closure-boundary} gives
\begin{equation}\label{eq:boundary-admissibility}
 H(u)(\overline\Omega)\Subset\Gamma_k.
\end{equation}

\begin{definition}[Linearized operator along the continuity path]
\label{def:linearized-operator}
At a solution $u$, set $F^{rs}=F^{rs}(H(u))$, with the coefficients defined
by \eqref{eq:differential-representation}, and write
\begin{equation}\label{eq:linearized-operator}
 \cL \psi:=\sum_{r,s=1}^{4n}F^{rs}\psi_{rs}-g_{t,z}(q,u)\psi,
 \qquad v:=u-\ubar.
\end{equation}
\end{definition}

By \eqref{eq:boundary-admissibility} and
Lemma~\ref{lem:concavity-ellipticity},
\[
 0<\lambda_u I_{4n}\le(F^{rs})\le\Lambda_u I_{4n}
 \quad\text{on }\overline\Omega.
\]
These constants depend on the compact set $H(u)(\overline\Omega)$.
Bounds independent of $u$ follow from \eqref{eq:global-c2} in
Section~\ref{sec:existence}.
Equations \eqref{eq:supporting-homogeneity}, \eqref{eq:concavity-tools}, and
\eqref{eq:tool} give
\begin{equation}\label{eq:trace}
 \sum_{r,s=1}^{4n}F^{rs}u_{rs}=g_t(q,u),\qquad
 \sum_{r=1}^{4n}F^{rr}=DF(H(u))[I_n]\ge1.
\end{equation}

We use the real derivatives in \eqref{eq:real-coordinate-derivatives},
with $1\le i,r,s\le4n$, and evaluate $g_{t,z}$ at $(q,u(q))$.
Equation~\eqref{eq:path} gives $g_{t,z}=t g_z(q,u)\ge0$.

\begin{lemma}[Product rules and weighted Young inequality]
\label{lem:linearized-calculus}
For smooth real functions $\psi,\chi$,
\begin{equation}\label{eq:linearized-product}
 \cL(\psi\chi)=\psi\cL\chi+\chi\cL\psi
 +2\sum_{r,s}F^{rs}\psi_r\chi_s+g_{t,z}\psi\chi.
\end{equation}
For a smooth function $Z=(Z_1,\ldots,Z_\ell)^T:\Omega\to\RR^\ell$,
where $\ell\in\{1,\ldots,4n\}$, set
\begin{equation}\label{eq:weighted-energy}
 \mathcal E(Z):=\sum_{j=1}^\ell\sum_{r,s=1}^{4n}
                   F^{rs}(Z_j)_r(Z_j)_s.
\end{equation}
Then
\begin{equation}\label{eq:linearized-square}
 \cL|Z|^2=2\mathcal E(Z)+2\sum_{j=1}^\ell Z_j\cL Z_j+g_{t,z}|Z|^2.
\end{equation}
For a smooth real vector field $X=\sum_{i=1}^{4n}X^i\partial_i$,
\begin{equation}\label{eq:vector-commutator}
 \cL(X\psi)=\sum_iX^i\cL(\psi_i)
   +2\sum_{r,s,i}F^{rs}(X^i)_r\psi_{is}
   +\sum_{r,s,i}F^{rs}(X^i)_{rs}\psi_i.
\end{equation}
For $\xi,\zeta\in\RR^{4n}$ and $\epsilon>0$,
\begin{equation}\label{eq:weighted-young}
 \begin{aligned}
 2\left|\sum_{r,s}F^{rs}\xi_r\zeta_s\right|
 &\le\epsilon\sum_{r,s}F^{rs}\xi_r\xi_s
       +\epsilon^{-1}\sum_{r,s}F^{rs}\zeta_r\zeta_s,\\
 \sum_{r,s}F^{rs}\zeta_r\zeta_s
 &\le|\zeta|^2\sum_rF^{rr}.
 \end{aligned}
\end{equation}
\end{lemma}
\begin{proof}
 The product rule and $F^{rs}=F^{sr}$ give \eqref{eq:linearized-product}.
 Applying it to $Z_j^2$ and summing gives \eqref{eq:linearized-square}.
For \eqref{eq:vector-commutator}, use
\[
 (X\psi)_{rs}=\sum_i\bigl(X^i\psi_{irs}+(X^i)_r\psi_{is}
                  +(X^i)_s\psi_{ir}+(X^i)_{rs}\psi_i\bigr).
\]
The first inequality in \eqref{eq:weighted-young} follows from
\[
 \sum_{r,s}F^{rs}
 (\sqrt\epsilon\,\xi_r\pm\epsilon^{-1/2}\zeta_r)
 (\sqrt\epsilon\,\xi_s\pm\epsilon^{-1/2}\zeta_s)\ge0.
\]
The second follows from
$\lambda_{\max}((F^{rs}))\le\sum_rF^{rr}$.
\end{proof}

\begin{lemma}[Absorption by a bounded square]\label{lem:square-absorption}
For $Z$ as above, suppose $|Z_j|\le K_Z$ and, for every $\epsilon>0$,
\[
 |\cL Z_j|\le C_\epsilon\left(1+\sum_rF^{rr}\right)
                 +\epsilon\mathcal E(Z),
 \qquad 1\le j\le\ell.
\]
Then
\begin{equation}\label{eq:square-absorption}
 \cL|Z|^2\ge\mathcal E(Z)-C\left(1+\sum_rF^{rr}\right).
\end{equation}
If the stronger bounds $|\cL Z_j|\le C(1+\sum_rF^{rr})$ hold, then
\begin{equation}\label{eq:square-strong}
 \cL|Z|^2\ge2\mathcal E(Z)-C\left(1+\sum_rF^{rr}\right).
\end{equation}
\end{lemma}
\begin{proof}
Since $g_{t,z}|Z|^2\ge0$, \eqref{eq:linearized-square} implies
\[
 \cL|Z|^2
 \ge(2-2\ell K_Z\epsilon)\mathcal E(Z)
     -2\ell K_ZC_\epsilon\left(1+\sum_rF^{rr}\right).
\]
Take $\epsilon=[2\ell(1+K_Z)]^{-1}$ to obtain \eqref{eq:square-absorption}.
Substituting the stronger bounds directly into
\eqref{eq:linearized-square} gives \eqref{eq:square-strong}.
\end{proof}

The next estimate controls the second derivative term in
\eqref{eq:vector-commutator} by derivatives in real tangential directions.

\begin{lemma}\label{lem:projection}
Fix $q\in\overline\Omega$ in the tubular neighborhood of
Subsection~\ref{sec:boundary-tools}, and let $\nu:=D\varrho(q)$ be
the value at $q$ of the normal field defined there.
Let $(\tau_\alpha)_{\alpha\in\mathscr I}$ be an orthonormal basis of $\nu^\perp$,
the tangent space to the level set of $\varrho$ through $q$. Then
\begin{equation}\label{eq:projection-est}
 \bigl|(F^{rs})D^2u\bigr|^2
 \le C_n\left[g_t(q,u)^2+
 \left(\sum_rF^{rr}\right)
 \sum_{\alpha \in\mathscr I} \sum_{r,s=1}^{4n}
 F^{rs}D^2u(\tau_\alpha,\partial_r)D^2u(\tau_\alpha,\partial_s)\right].
\end{equation}
The matrix norm on the left is the Frobenius norm in
\eqref{eq:matrix-norms}.
\end{lemma}
\begin{proof}
Set $M=(F^{rs})D^2u$, $\Pi=I_{4n}-\nu\otimes\nu$, and $\zeta=M\nu$.
In particular, $|M|^2=\sum_{r,s=1}^{4n}M_{rs}^2$.
By Lemma~\ref{lem:concavity-ellipticity}, $(F^{rs})$ commutes with
$I,J,K$ and $\Real(H(u))=\ProjR(D^2u)$. Hence
\[
 \begin{aligned}
 \ProjR(M)
 &=\frac14\sum_{Q\in\mathscr Q}Q^T(F^{rs})D^2uQ
   =(F^{rs})\ProjR(D^2u),\\
 \ProjR(M)^T
 &=\ProjR(D^2u)(F^{rs})=(F^{rs})\ProjR(D^2u)=\ProjR(M).
 \end{aligned}
\]
For a real matrix $C$, write
$\operatorname{skew}C:=\tfrac12(C-C^T)$ for its skew-symmetric part.
Since $M=M\Pi+\zeta\otimes\nu$, the symmetry of $\ProjR(M)$ gives
\begin{equation}\label{eq:projection-skew-identity}
 \operatorname{skew}\ProjR(\zeta\otimes\nu)
 =-\operatorname{skew}\ProjR(M\Pi).
\end{equation}
Choose a quaternionic orthonormal basis $e_1,\ldots,e_n$ with
$e_n=\iota^{-1}(\nu)$, and let $U=(e_1,\ldots,e_n)$.
By \eqref{eq:unitary-real-orthogonality} and
\eqref{eq:real-representation-image}, $\Real(U)$ is real orthogonal
and commutes with $I,J,K$. In the corresponding real coordinates,
\eqref{eq:real-basis-change} gives
\[
 \begin{aligned}
 \widehat\nu&:=\Real(U)^T\nu=\iota((0,\ldots,0,1)^T),\\
 \widehat\zeta&:=\Real(U)^T\zeta
    =\iota((\zeta_1,\ldots,\zeta_n)^T),
 \qquad \zeta_a:=e_a^*\iota^{-1}(\zeta)\in\HH.
 \end{aligned}
\]
Thus $\iota^{-1}(\zeta)=\sum_{a=1}^ne_a\zeta_a$.
The commutation relations and \eqref{eq:projection} give
\[
 \begin{aligned}
 \ProjR(\widehat\zeta\otimes\widehat\nu)
 &=\frac14\sum_{Q\in\mathscr Q}
       Q^T\Real(U)^T(\zeta\otimes\nu)\Real(U)Q\\
 &=\Real(U)^T\ProjR(\zeta\otimes\nu)\Real(U).
 \end{aligned}
\]
Taking skew-symmetric parts and using orthogonal invariance of the
Frobenius norm yields
\[
 |\operatorname{skew}\ProjR(\zeta\otimes\nu)|
 =|\operatorname{skew}\ProjR(\widehat\zeta\otimes\widehat\nu)|.
\]
In $4\times4$ real blocks, for $1\le a,b\le n$,
\[
 [\widehat\zeta\otimes\widehat\nu]_{ab}
 =\begin{cases}
   \iota(\zeta_a)(1,0,0,0),&b=n,\\
   0,&b<n.
  \end{cases}
\]
Since $I,J,K$ are block diagonal,
$[\ProjR(\widehat\zeta\otimes\widehat\nu)]_{ab}=0$ for $b<n$.
By \eqref{eq:quaternionic-structure-blocks} and
\eqref{eq:left-multiplication-columns},
\[
 \begin{aligned}
 4[\ProjR(\widehat\zeta\otimes\widehat\nu)]_{an}
 &=\sum_{Q_0\in\{I_4,I_0,J_0,K_0\}}
       Q_0^T\iota(\zeta_a)(1,0,0,0)Q_0\\
 &=\bigl(\iota(\zeta_a),\iota(\zeta_a\mathbf i),
                 \iota(\zeta_a\mathbf j),\iota(\zeta_a\mathbf k)\bigr)
 =L_{\zeta_a}.
 \end{aligned}
\]
Since $L_{\zeta_n}-L_{\zeta_n}^T=2L_{\operatorname{Im}\zeta_n}$,
\[
 \operatorname{skew}\ProjR(\widehat\zeta\otimes\widehat\nu)=\frac18
 \begin{pmatrix}
  0&\cdots&0&L_{\zeta_1}\\
  \vdots&\ddots&\vdots&\vdots\\
  0&\cdots&0&L_{\zeta_{n-1}}\\
  -L_{\zeta_1}^T&\cdots&-L_{\zeta_{n-1}}^T&2L_{\operatorname{Im}\zeta_n}
 \end{pmatrix}.
\]
Also, \eqref{eq:real-inner-product} and orthogonality of $\Real(U)$ give
\[
 \begin{aligned}
 \zeta\cdot\nu&=\widehat\zeta\cdot\widehat\nu=\Ree\zeta_n,\\
 |\zeta-(\zeta\cdot\nu)\nu|^2
 &=|\widehat\zeta-(\Ree\zeta_n)\widehat\nu|^2
 =\sum_{a=1}^{n-1}|\zeta_a|^2+|\operatorname{Im}\zeta_n|^2.
 \end{aligned}
\]
Using $|L_c|^2=4|c|^2$ from
\eqref{eq:real-representation-trace-norm}, we obtain
\begin{equation}\label{eq:projection-skew-bound}
 \begin{aligned}
 |\operatorname{skew}\ProjR(\zeta\otimes\nu)|^2
 &=\frac1{64}\left(2\sum_{a=1}^{n-1}|L_{\zeta_a}|^2
                 +4|L_{\operatorname{Im}\zeta_n}|^2\right)\\
 &=\frac18\sum_{a=1}^{n-1}|\zeta_a|^2
       +\frac14|\operatorname{Im}\zeta_n|^2
 \ge\frac18|\zeta-(\zeta\cdot\nu)\nu|^2.
 \end{aligned}
\end{equation}
Both $\operatorname{skew}$ and $\ProjR$ are contractions for the
Frobenius norm. Equations \eqref{eq:projection-skew-identity} and
\eqref{eq:projection-skew-bound} give
\[
 |\zeta-(\zeta\cdot\nu)\nu|
 \le\sqrt8\,|\operatorname{skew}\ProjR(M\Pi)|
 \le\sqrt8\,|M\Pi|.
\]
Also, \eqref{eq:trace} and $\trR(\zeta\otimes\nu)=\zeta\cdot\nu$ give
\[
 \begin{aligned}
 |\zeta\cdot\nu|
 &=|g_t(q,u)-\trR(M\Pi)|
 \le g_t(q,u)+\sqrt{4n}\,|M\Pi|,\\
 |M|&\le|M\Pi|+|\zeta|
 \le g_t(q,u)+(1+\sqrt8+\sqrt{4n})|M\Pi|.
 \end{aligned}
\]
Finally, each eigenvalue of $(F^{rs})$ is at most $\sum_rF^{rr}$, so
$(F^{rs})^2\le(\sum_rF^{rr})(F^{rs})$. Consequently,
\[
 \begin{aligned}
 |M\Pi|^2
 &=\sum_\alpha|M\tau_\alpha|^2
   =\sum_\alpha(D^2u\,\tau_\alpha)^T(F^{rs})^2(D^2u\,\tau_\alpha)\\
 &\le\left(\sum_rF^{rr}\right)
 \sum_{\alpha,r,s}F^{rs}
 D^2u(\tau_\alpha,\partial_r)D^2u(\tau_\alpha,\partial_s).
 \end{aligned}
\]
Squaring the preceding bound for $|M|$ proves \eqref{eq:projection-est}.
\end{proof}

If $\cL\psi\ge0$, the maximum principle gives
$\sup_\Omega\psi\le\max\{0,\sup_{\partial\Omega}\psi\}$.
Equivalently, $\cL\psi\le0$ and $\psi|_{\partial\Omega}\ge0$ imply $\psi\ge0$.
If $\psi\ge0$ on the interior side and $\psi(p)=0$ at a boundary point, then
\[
 \psi_\nu(p)=\lim_{h\downarrow0}
       \frac{\psi(p)-\psi(p-h\nu(p))}{h}\le0.
\]

\subsection{Comparison and first derivatives}
\begin{lemma}[Comparison and $C^0$ bounds]\label{lem:comparison-C0}
Every
smooth admissible solution along \eqref{eq:path} satisfies
\[
 \ubar\le u\le h\quad\text{on }\closure\Omega.
\]
\end{lemma}

\begin{proof}
Integrating the differential of $F$ along the segment joining
$H(u)$ and $H(\ubar)$ gives continuous positive definite real
coefficients
\[
 \widetilde F^{rs}(q)=\int_0^1
 F^{rs}\!\left(H(u)(q)+\theta\bigl(H(\ubar)(q)-H(u)(q)\bigr)\right)\,\dd\theta.
\]
Since $\ubar$ is a subsolution of \eqref{eq:equation} and $u$ is a solution,
\[
 \begin{aligned}
 \sum_{r,s=1}^{4n}\widetilde F^{rs}(\ubar-u)_{rs}
 &=F(H(\ubar))-F(H(u))\\
 &\ge g_t(q,\ubar)-g_t(q,u)\\
 &=\left[\int_0^1g_{t,z}(q,u+\theta(\ubar-u))\dd\theta\right](\ubar-u).
 \end{aligned}
\]
Define
\[
 c(q):=\int_0^1g_{t,z}(q,u+\theta(\ubar-u))\,\dd\theta\ge0,
 \qquad \zeta:=\ubar-u.
\]
Thus
\[
 \sum_{r,s=1}^{4n}\widetilde F^{rs}\zeta_{rs}-c(q)\zeta\ge0,
 \qquad \zeta|_{\bd}=0
 \quad\Longrightarrow\quad\zeta\le0.
\] Also,
$\Delta u=4\sigma_1(H(u))>0$, so comparison with $h$ gives $u\le h$.
\end{proof}

In particular, $u(\overline\Omega)\subset\mathcal I$.
Lemma~\ref{lem:comparison-C0} and the one-sided derivative inequality in
Subsection~\ref{sec:linearized-calculus} give
\[
 h_\nu\le u_\nu\le\ubar_\nu\quad\text{on }\bd.
\]
The common tangential boundary derivatives are those of $\varphi$, so
$\sup_{\bd}|Du|\le C$.

For a constant real unit vector $e$, the constant coefficients of $H$ and
\eqref{eq:tool} give
\[
 D_eH(u)=H(u_e),\qquad
 DF(H(u))[H(u_e)]=\sum_{r,s}F^{rs}(u_e)_{rs}.
\]
Differentiating \eqref{eq:equation} therefore yields
\begin{equation}\label{eq:first-derivative}
 \sum_{r,s}F^{rs}(u_e)_{rs}=g_{t,e}+g_{t,z}u_e,
 \qquad \cL u_e=g_{t,e}.
\end{equation}
Here $g_{t,e}$ differentiates the first variable with $z$ fixed.
The notation $g_{t,r},g_{t,ee},g_{t,ez},g_{t,zz}$ is used accordingly,
with all derivatives evaluated at $(q,u(q))$.

Set $C_{\ubar}:=\max_{\closure\Omega}\ubar$, so that
$\ubar-C_{\ubar}\le0$.
By \eqref{eq:supporting-homogeneity} and $g_{t,z}\ge0$,
\begin{equation}\label{eq:subsolution-test}
 \begin{aligned}
 \cL(\ubar-C_{\ubar})
 &=DF(H(u))[H(\ubar)]+g_{t,z}(C_{\ubar}-\ubar)\\
 &\ge F(H(\ubar))=g_0(q)\ge g_*>0.
 \end{aligned}
\end{equation}
Choose a fixed $K_1>0$ such that $K_1g_*\ge\sup|g_{t,e}|$, uniformly in
$e$ and $t$. By \eqref{eq:first-derivative},
\[
 \cL\bigl(\pm u_e+K_1(\ubar-C_{\ubar})\bigr)\ge0.
\]
Since $\ubar-C_{\ubar}\le0$, the maximum principle and the boundary
gradient bound give
\[
 \pm u_e+K_1(\ubar-C_{\ubar})
 \le\max\{0,\sup_{\bd}(\pm u_e)\}\le C.
\]
The boundedness of $\ubar-C_{\ubar}$ proves
\begin{equation}\label{eq:c1}
 \|u\|_{C^1(\closure\Omega)}\le C.
\end{equation}

\subsection{Reduction of the second derivative estimate to the boundary}
For a constant real unit vector $e$, the two sides of
\eqref{eq:equation} have second derivatives
\[
 \begin{aligned}
 D_e^2F(H(u))
 &=DF(H(u))[H(u_{ee})]+D^2F(H(u))[H(u_e),H(u_e)],\\
 D_e^2(g_t(q,u(q)))
 &=g_{t,ee}+2g_{t,ez}u_e+g_{t,zz}u_e^2+g_{t,z}u_{ee}.
 \end{aligned}
\]
Subtracting $g_{t,z}u_{ee}$ and using \eqref{eq:tool} gives
\begin{equation}\label{eq:second-linearized}
 \cL u_{ee}
 =g_{t,ee}+2g_{t,ez}u_e+g_{t,zz}u_e^2
   -D^2F(H(u))[H(u_e),H(u_e)]\ge-C.
\end{equation}
The inequality follows from \eqref{eq:concavity-second-variation} and
\eqref{eq:c1}. Increasing the constant $K_1$ chosen above if necessary,
\eqref{eq:subsolution-test} gives
$\cL\bigl(u_{ee}+K_1(\ubar-C_{\ubar})\bigr)\ge0$.
The maximum principle gives
\begin{equation}\label{eq:second-derivative-reduction}
 \sup_{\closure\Omega}u_{ee}
 \le\max\{0,\sup_{\bd}u_{ee}\}+C,
\end{equation}
where $C$ is independent of $e$. If $\lambda_i^{\RR}$ are the eigenvalues
of $D^2u$, then
\[
 0<\Delta u=\sum_i\lambda_i^{\RR}
 \quad\Longrightarrow\quad
 \lambda_i^{\RR}>-\sum_{j\ne i}\lambda_j^{\RR}
 \ge-(4n-1)\max_j\lambda_j^{\RR}.
\]
This bounds the full real Hessian once the boundary upper bound is known
\cite[Section 3, Step 2]{Zhu2016}.

\subsection{Barrier near the boundary}
Use $d=-\varrho$ and $\nu=D\varrho$ from Subsection~\ref{sec:boundary-tools}.
Choose $\delta_0>0$ so that
$\Omega_{\delta_0}:=\{q\in\Omega:0<d(q)<\delta_0\}$ lies in the
fixed tubular neighborhood. Then \eqref{eq:distance-identities} holds there.

\begin{lemma}[Barrier construction]\label{lem:collar}
For $0<\delta<\delta_1=\delta_1(\Omega,\ubar,g^*)$, there are positive constants $\eta_0,N,\varepsilon_0$ depending only on $(\Omega,\ubar,g^*)$ such that
in $\Omega_\delta=\{q\in\Omega:0<d(q)<\delta\}$ the function
$\Phi=v+\eta_0d-Nd^2$ satisfies
\begin{equation}\label{eq:collar-ineq}
 \Phi\ge v+\tfrac12\eta_0d\ge0,\qquad
 \cL \Phi\le-\varepsilon_0\left(1+\sum_{r=1}^{4n}F^{rr}\right),
\end{equation}
where $v=u-\ubar$ and $\delta_1\le \delta_0$.
\end{lemma}

\begin{proof}
Choose small fixed $\varepsilon,\eta_0>0$ so that
\[
 \varepsilon|I_n|+\eta_0\sup_{\Omega_{\delta_0}}|H(d)|
 \le\tfrac12\delta_{\ubar}.
\]
For $0\le\tau\le\eta_0$, the matrix
\[
 Y=H(\ubar)-\varepsilon I_n-\tau H(d)
\]
then satisfies
\[
 |Y-H(\ubar)|\le\tfrac12\delta_{\ubar},\qquad
 \operatorname{dist}(Y,\partial\Gamma_k)\ge\tfrac12\delta_{\ubar}.
\]
Thus these matrices stay in a fixed compact subset of $\Gamma_k$.
At each point, choose a quaternionic orthonormal frame
$U=(e_1,\ldots,e_n)$ with $e_n=\iota^{-1}(\nu)$ and set
$E_q=(e_n\HH)^\perp$. By \eqref{eq:subspace-matrix},
$Y|_{E_q}$ is represented by the leading $(n-1)\times(n-1)$ block of $U^*YU$.
Lemma~\ref{lem:block-identity} and compactness give
$\sigma_{k-1}(Y|_{E_q})\ge\gamma>0$.
The constants $\varepsilon,\eta_0,\gamma$ depend only on $(\Omega,\ubar)$.

Since $Dd=-\nu$, \eqref{eq:projected-rank-one} gives
\[
 \Proj(Dd\otimes Dd)=\Proj(\nu\otimes\nu)=\frac14e_ne_n^*,
 \qquad U^*\Proj(Dd\otimes Dd)U=\frac14E_{nn}.
\]
Thus \eqref{eq:positive-perturbation} gives
$Z:=Y+2N\Proj(Dd\otimes Dd)\in\Gamma_k$.
Applying \eqref{eq:rank-one-sigma} to $U^*YU$ yields
\[
 \begin{aligned}
 \sigma_k(Z)&=\sigma_k\left(U^*YU+\frac N2E_{nn}\right)\\
 &=\sigma_k(Y)+\frac N2
       \sigma_{k-1}(Y|_{E_q})\ge\frac{N\gamma}{2}.
 \end{aligned}
\]
Choose
\[
 N\ge\frac{2\binom nk(g^*+1)^k}{\gamma},\qquad
 \delta\le\delta_1:=\min\{\frac{\eta_0}{2N},\delta_0\}.
\]
By \eqref{eq:supporting-homogeneity},
\[
 DF(H(u))[Z]\ge F(Z)
 =\binom nk^{-1/k}\sigma_k(Z)^{1/k}
 \ge\left(\frac{N\gamma}{2\binom nk}\right)^{1/k}\ge g^*+1.
\]
At a point in $\Omega_\delta$, take $\tau=\eta_0-2Nd\in[0,\eta_0]$.
The chain rule \eqref{eq:hessian-chain} gives
\[
 H(d^2)=2dH(d)+2\Proj(Dd\otimes Dd).
\]
Hence
\[
 \begin{aligned}
 H(\Phi)
 &=H(u)-H(\ubar)+(\eta_0-2Nd)H(d)
             -2N\Proj(Dd\otimes Dd)\\
 &=H(u)-Z-\varepsilon I_n.
 \end{aligned}
\]
Since $N\delta\le\eta_0/2$,
\[
 \Phi=v+d(\eta_0-Nd)\ge v+\tfrac12\eta_0d\ge0.
\]
Using $g_{t,z}\ge0$, \eqref{eq:trace}, and \eqref{eq:tool}, we obtain
\[
 \begin{aligned}
 \cL \Phi
 &\le\sum_{r,s=1}^{4n}F^{rs}\Phi_{rs}\\
 &=g_t-DF(H(u))[Z]-\varepsilon\sum_{r=1}^{4n}F^{rr}\\
 &\le g^*-F(Z)-\varepsilon\sum_{r=1}^{4n}F^{rr}
 \le-1-\varepsilon\sum_{r=1}^{4n}F^{rr}.
 \end{aligned}
\]
Take $\varepsilon_0=\min\{1,\varepsilon\}$.
\end{proof}

\subsection{Tangential and mixed boundary derivatives}
\label{sec:tangential-estimates}\label{sec:mixed-estimates}
\noindent\textit{Tangential derivatives.}
Equations \eqref{eq:boundary-tangent} and \eqref{eq:c1} give
\[
 D^2u(\tau,\tau')
 =D^2\ubar(\tau,\tau')+v_\nu D^2\varrho(\tau,\tau'),
 \qquad |\tau|=|\tau'|=1,\quad \tau,\tau'\in T\bd.
\]
Consequently,
\[
 |D^2u(\tau,\tau')|
 \le\|D^2\ubar\|_{C^0(\bd)}
       +\|v_\nu\|_{C^0(\bd)}\|D^2\varrho\|_{C^0(\bd)}\le C.
\] On the quaternionic
tangent space $E_q=(e_n(q)\HH)^\perp$, \eqref{eq:boundary-restriction} yields
\begin{equation}\label{eq:tangential-block}
 B:=H(u)|_{E_q}=B_0+v_\nu R,\qquad
 B_0:=H(\ubar)|_{E_q},\quad R:=H(\varrho)|_{E_q}.
\end{equation}
Thus Lemma~\ref{lem:block-identity} and \eqref{eq:c1} give
\begin{equation}\label{eq:Bbounds}
 B\in\Gamma_{k-1}^{(n-1)},\qquad |B|\le C\quad\text{on }\bd.
\end{equation}
Lemma~\ref{lem:block-identity} also gives
$B_0(q)\in\Gamma_{k-1}^{(n-1)}$. By continuity and compactness of $\bd$,
\[
 \sup_{q\in\bd}|B_0(q)|\le C,\qquad
 \inf_{q\in\bd}\operatorname{dist}
    \bigl(B_0(q),\partial\Gamma_{k-1}^{(n-1)}\bigr)>0,
\]
with bounds depending only on $\Omega$ and $\ubar$.

\medskip\noindent\textit{Mixed tangential-normal derivatives.}
Fix $p\in\bd$, translate $p$ to $0$, and use the coordinates and fields
in \eqref{eq:tangent-fields}. Thus $\partial_{x_n}=\nu(0)$,
$X_\alpha(0)=\partial_\alpha$, and $D_\nu X_\alpha=0$ in the fixed
coordinate neighborhood.

Choose a fixed $\delta>0$ so small that
$\Omega\cap B_\delta(0)$ is contained in this coordinate neighborhood
and in the collar of Lemma~\ref{lem:collar}. This is possible because
$d(q)\le|q|<\delta$ there.

Set
\[
 Z_\alpha=X_\alpha v,\qquad
 Z=(Z _\alpha) _{\alpha\in\mathscr I},
\]
with any fixed ordering of $\mathscr I$ in \eqref{eq:weighted-energy}.
All sums over tangential indices below range over $\mathscr I$.
By \eqref{eq:c1}, $|Z|\le C$ and $Z=0$ on the boundary portion.
Writing $X_\alpha=\sum_iX_\alpha^i\partial_i$,
the commutator formula \eqref{eq:vector-commutator}, applied to $u$,
and \eqref{eq:first-derivative} give
\begin{equation}\label{eq:mixed-commutator}
 \begin{aligned}
 \cL(X_\alpha v)={}&\sum\limits_{i}X_\alpha ^i g_{t, i}
 -\cL(X_\alpha\ubar)\\
 &+2\sum_{r,s, i=1}^{4n}F^{rs}(X_\alpha ^i)_r u_{ is}
 +\sum_{r,s, i=1}^{4n}F^{rs}(X_\alpha ^i)_{rs}u _i.
 \end{aligned}
\end{equation}
The coefficient bounds following \eqref{eq:tangent-fields} and
\eqref{eq:c1} give
\[
 \left|\sum\limits_{i}X_\alpha ^i g_{t, i}-\cL(X_\alpha\ubar)
       +\sum_{r,s, i}F^{rs}(X_\alpha ^i)_{rs}u _i\right|
 \le C\left(1+\sum_rF^{rr}\right).
\]
Indeed, for $r\ne s$, positive definiteness gives
\[
 0<\det\begin{pmatrix}F^{rr}&F^{rs}\\F^{rs}&F^{ss}\end{pmatrix}
   =F^{rr}F^{ss}-(F^{rs})^2.
\]
Thus $|F^{rs}|\le\sqrt{F^{rr}F^{ss}}$ for all $r,s$, and
\begin{equation}\label{eq:linearized-coefficient-sum}
 \sum_{r,s}|F^{rs}|\le\left(\sum_r\sqrt{F^{rr}}\right)^2
 \le4n\sum_rF^{rr},
\end{equation}
while all the other factors in these terms have fixed bounds.
The remaining term satisfies
\[
 \left|\sum_{r,s, i}F^{rs}(X_\alpha ^i)_r u_{ is}\right|
 =\left|\sum_{r, i}(X_\alpha ^i)_r[(F^{rs})D^2u]_{ ri}\right|
 \le |DX_\alpha|\,|(F^{rs})D^2u|
 \le C|(F^{rs})D^2u|.
\]
To apply \eqref{eq:projection-est}, choose an orthonormal basis
 $(\tau_{\alpha'})_{\alpha'\in\mathscr I}$ of $\nu^\perp$ and write
 $\tau_{\alpha'}=\sum_{i=1}^{4n}\tau_{\alpha'}^i\partial_i$.
Since $\tau_{\alpha'}\cdot D\varrho=0$ and $\varrho_{x_n}\ge1/2$,
\[
 \tau _{\alpha'}^{ 4n-3}=-\sum_{\alpha \in\mathscr I}
           \frac{\varrho_\alpha}{\varrho_{x_n}}\tau _{\alpha'}^\alpha,
 \qquad \tau _{\alpha'}=\sum_{\alpha \in\mathscr I} \tau _{\alpha'}^\alpha X_\alpha,
 \qquad \sum_{\alpha \in\mathscr I} (\tau _{\alpha'}^\alpha)^2\le1.
\]
Cauchy--Schwarz for the inner product $\xi^T(F^{rs})\zeta$ gives
\[
 \begin{aligned}
 (D^2u\,\tau _{\alpha'})^T(F^{rs})(D^2u\,\tau _{\alpha'})
 &\le\left(\sum_\alpha(\tau _{\alpha'}^\alpha)^2\right)
        \sum_\alpha(D^2u\,X_\alpha)^T(F^{rs})(D^2u\,X_\alpha)\\
 &\le\sum_\alpha(D^2u\,X_\alpha)^T(F^{rs})(D^2u\,X_\alpha).
 \end{aligned}
\]
Summing over $\alpha'$ yields
\begin{equation}\label{eq:mixed-frame-comparison}
 \begin{aligned}
 &\sum_{ \alpha',r,s}F^{rs}
       D^2u(\tau _{\alpha'},\partial_r)D^2u(\tau _{\alpha'},\partial_s)\\
 &\quad\le(4n-1)\sum_{\alpha,r,s}F^{rs}
       D^2u(X_\alpha,\partial_r)D^2u(X_\alpha,\partial_s).
 \end{aligned}
\end{equation}
Differentiating $Z_\alpha=X_\alpha v$ gives
\[
 D^2u(X_\alpha,\partial_r)-(Z_\alpha)_r
 =-\sum\limits_{i}(X_\alpha ^i)_r v _i+\sum\limits_{i}X_\alpha ^i\ubar_{ ir},
 \qquad
 \sum_{\alpha,r}|D^2u(X_\alpha,\partial_r)-(Z_\alpha)_r|^2\le C.
\]
Applying \eqref{eq:weighted-young} with $\epsilon=1$ to this decomposition,
we obtain
\begin{equation}\label{eq:mixed-tangential-energy}
 \begin{aligned}
 &\sum_{\alpha,r,s}F^{rs}
       D^2u(X_\alpha,\partial_r)D^2u(X_\alpha,\partial_s)\\
 &\quad\le2\mathcal E(Z)
   +2\sum_{\alpha,r,s}F^{rs}
     \bigl[D^2u(X_\alpha,\partial_r)-(Z_\alpha)_r\bigr]
     \bigl[D^2u(X_\alpha,\partial_s)-(Z_\alpha)_s\bigr]\\
 &\quad\le2\mathcal E(Z)+C\sum_rF^{rr}.
 \end{aligned}
\end{equation}
Substituting \eqref{eq:mixed-frame-comparison} and
\eqref{eq:mixed-tangential-energy} into \eqref{eq:projection-est}, and
absorbing the factor $4n-1$ into $C_n$, gives
\[
 \begin{aligned}
 |(F^{rs})D^2u|^2
 &\le C_n\left[(g^*)^2+
   \left(\sum_rF^{rr}\right)
       \left(2\mathcal E(Z)+C\sum_rF^{rr}\right)\right]\\
 &\le C\left[1+\left(\sum_rF^{rr}\right)\mathcal E(Z)
                 +\left(\sum_rF^{rr}\right)^2\right].
 \end{aligned}
\]
Taking square roots and using \eqref{eq:mixed-commutator} yields
\begin{equation}\label{eq:mixed-linearized}
 |\cL Z_\alpha|
 \le C\left(1+\sum_rF^{rr}\right)
      +C\left(\sum_rF^{rr}\right)^{1/2}\mathcal E(Z)^{1/2}.
\end{equation}
For $\varepsilon>0$, the last term is bounded by
\[
 C\left(\sum_rF^{rr}\right)^{1/2}\mathcal E(Z)^{1/2}
 \le\varepsilon\mathcal E(Z)+\frac{C^2}{4\varepsilon}\sum_rF^{rr}.
\]
Thus
\begin{equation}\label{eq:tangential-square-estimate}
 |\cL Z_\alpha|
 \le C_\varepsilon\left(1+\sum_rF^{rr}\right)
       +\varepsilon\mathcal E(Z).
\end{equation}
Since $Z=(Z_\alpha)_{\alpha\in\mathscr I}$ is bounded, Lemma~\ref{lem:square-absorption}
applied to \eqref{eq:tangential-square-estimate} gives
\begin{equation}\label{eq:tangential-square}
 \cL|Z|^2\ge\mathcal E(Z)-C\left(1+\sum_rF^{rr}\right).
\end{equation}

For a fixed $\alpha$, consider
\[
 \Psi_\pm=K_1\Phi+K_2|q|^2\pm Z_\alpha
                 -|Z|^2
 \quad\text{in }\Omega\cap B_\delta(0).
\]
On $\bd\cap B_\delta(0)$, $\Phi=Z=0$ and $\Psi_\pm=K_2|q|^2\ge0$.
Choose $K_2$ so that
\[
 K_2\delta^2\ge\sup_{\Omega\cap B_\delta(0)}(|Z_\alpha|+|Z|^2).
\]
Then $\Phi\ge0$ gives $\Psi_\pm\ge0$ on
$\Omega\cap\partial B_\delta(0)$ as well.
Also
\[
 \cL|q|^2=2\sum_rF^{rr}-(g_{t,z})|q|^2\le2\sum_rF^{rr}.
\]
Combining \eqref{eq:collar-ineq}, \eqref{eq:tangential-square-estimate} and
\eqref{eq:tangential-square}, and choosing $\varepsilon=1/2$, gives
\[
  \cL\Psi_\pm
  \le \bigl[-K_1\varepsilon_0+C_{1/2}+C+2K_2\bigr]
       \left(1+\sum_rF^{rr}\right)-\tfrac12\mathcal E(Z).
 \]
Increase $K_1$ so that $K_1\ge(C_{1/2}+C+2K_2)/\varepsilon_0$.
Then the right-hand side is nonpositive.
The maximum principle implies $\Psi_\pm\ge0$ in $\Omega\cap B_\delta(0)$.
At $0$, $Z=0$ gives $(|Z|^2)_\nu=2\sum_\alpha Z_\alpha(Z_\alpha)_\nu=0$,
and $(|q|^2)_\nu=0$. Since $\Psi_\pm(0)=0$, the one-sided derivative
inequality in Subsection~\ref{sec:linearized-calculus} gives
\[
 0\ge(\Psi_\pm)_\nu(0)=K_1\Phi_\nu(0)\pm(Z_\alpha)_\nu(0).
\]
By \eqref{eq:tangent-fields}, $D_\nu X_\alpha=0$ and
$X_\alpha(0)=\partial_\alpha$. Equation \eqref{eq:iterated-directional-derivative} gives
\[
  ( Z_\alpha )_\nu(0)
 =D^2v(\nu,X_\alpha)(0)+Dv(0)\cdot D_\nu X _\alpha(0)
 =D^2v(\partial_\alpha,\nu)(0).
\]
The two signs in the preceding inequality therefore give
\[
 |(Z_\alpha)_\nu(0)|\le-K_1\Phi_\nu(0)
 =K_1(\eta_0-v_\nu(0))\le C.
\]
At $0$, every unit vector $\tau\in T_0\bd$ has the form
$\tau=\sum_\alpha\tau^\alpha X_\alpha(0)$ with
$\sum_\alpha(\tau^\alpha)^2=1$. Hence
\[
 |D^2u(\tau,\nu)(0)|
 \le\left(\sum_\alpha|(Z_\alpha)_\nu(0)|^2\right)^{1/2}
       +\|D^2\ubar\|_{C^0(\bd)}\le C.
\]
Since the original boundary point $p$ was arbitrary, this bound holds
on $\bd$. For every constant real unit vector $e$, the tangential and
mixed estimates give
\[
 u_{ee}\le C+C(u_{\nu\nu})^+\quad\text{on }\bd.
\]
Together with \eqref{eq:second-derivative-reduction} and $\Delta u>0$,
this proves the following estimate.

\begin{proposition}[Partial boundary $C^2$ estimate]
\label{prop:partial-boundary-C2}
Every smooth admissible solution of \eqref{eq:equation} satisfies
\[
 \sup_{q\in\bd}\sup_{\substack{\tau\in T_q\bd\\
                                  |\tau|=1}}
 \bigl(|D^2u(\tau,\tau)|+|D^2u(\tau,\nu)|\bigr)\le C.
\]
In particular,
\begin{equation}\label{eq:c2-reduction}
 \|u\|_{C^2(\closure\Omega)}
 \le C\left(1+\max_{\bd}(u_{\nu\nu})^+\right).
\end{equation}

\end{proposition}

\section{The double normal boundary estimate}\label{sec:boundary}
Proposition~\ref{prop:partial-boundary-C2} reduces the $C^2$ estimate to
an upper bound for $u_{\nu\nu}$ on $\bd$. At each $q\in\bd$, choose a
quaternionic orthonormal frame $e=(e_1(q),\ldots,e_n(q))$ with
$e_n(q)=\iota^{-1}(\nu(q))$. Equation~\eqref{eq:orthonormal-basis-change} gives
\[
 [H(u)(q)]_e=(e_a(q)^*H(u)(q)e_b(q))_{a,b=1}^n
 =\begin{pmatrix}B(q)&\beta(q)\\\beta(q)^*&a_0(q)\end{pmatrix}.
\]
By \eqref{eq:offdiagonal-real}, \eqref{eq:normal-entry}, and the
tangential and mixed estimates,
\begin{equation}\label{eq:boundary-block-bounds}
 |B|+|\beta|+|4a_0-u_{\nu\nu}|\le C\quad\text{on }\bd.
\end{equation}
Thus the block identity \eqref{eq:block} bounds $a_0$ once
$\sigma_{k-1}(B)$ has a uniform positive lower bound.
We first estimate $u_{\nu\nu}$ at a point minimizing $m(B)$, then use
the equation and \eqref{eq:margin-bounds} to obtain this lower bound
on the whole boundary.

\subsection{The boundary minimum}\label{sec:margin}
Choose a minimizer $p\in\bd$ of $m(B)$ and translate $p$ to $0$,
retaining the notation for the domain and functions.
By Lemma~\ref{lem:4.1}, choose $\Theta\in\mathscr S$ such that
\begin{equation}\label{eq:boundary-minimum}
 m_*:=\min_{q\in\bd}m(B(q))
     =m(B(0))=\pair{\Theta}{B(0)}>0.
\end{equation}
By \eqref{eq:Bbounds}, $m_*\le C$. Since
$B_0=H(\ubar)|_{E_q}$ stays in a fixed compact subset of
$\Gamma_{k-1}^{(n-1)}$, set
\[
 b_*:=\min_{q\in\bd}m(B_0(q))>0.
\]
If $m_*\ge b_*/2$, \eqref{eq:margin-bounds} already gives
$\sigma_{k-1}(B)\ge\binom{n-1}{k-1}(b_*/2)^{k-1}$ on $\bd$.
It remains to treat $m_*<b_*/2$.
At $0$, put $b(0)=\pair{\Theta}{B_0(0)}$ and
$\eta(0)=\pair{\Theta}{R(0)}$. By \eqref{eq:tangential-block} and
\eqref{eq:c1},
\[
 m_*=b(0)+v_\nu(0)\eta(0),\qquad
 b(0)\ge b_*,\qquad 0\le-v_\nu(0)\le C.
\]
Hence $v_\nu(0)<0$ and
\begin{equation}\label{eq:normal-weight-positive}
 \eta(0)=\frac{b(0)-m_*}{-v_\nu(0)}
 \ge\frac{b_*}{2C}=: \eta_*>0.
\end{equation}

Keeping $e_n=\iota^{-1}(\nu(0))$ fixed, choose the tangential
orthonormal basis $e_1,\ldots,e_{n-1}$ so that, in this basis,
$\Theta=\operatorname{diag}(\mu_1,\ldots,\mu_{n-1})$,
where $\mu_a\ge0$ and $\sum_a\mu_a=1$.
We extend this basis to a smooth quaternionic orthonormal frame of
$E_q$ near $0$ to extend $b(0)$ and $\eta(0)$ with controlled $C^2$ norms.
Here $q$ is the variable point in the fixed translated coordinates.
All frame vectors are expressed in the standard basis of $\HH^n$.
Extend $\ubar$ smoothly across $\partial\Omega$ to the fixed
tubular neighborhood where $\varrho$ is defined, and put
$e_n(q)=\iota^{-1}(D\varrho(q))$ and $E_q=(e_n(q)\HH)^\perp$.
Apply the quaternionic Gram--Schmidt
procedure to $e_n(q),e_1,\ldots,e_{n-1}$, where $e_1,\ldots,e_{n-1}$
are the vectors chosen at $0$, extended as constant vectors in $\HH^n$.
Explicitly, for $1\le a\le n-1$, set
\[
 \widetilde e_a(q)=e_a-e_n(q)(e_n(q)^*e_a)
     -\sum_{b=1}^{a-1}e_b(q)(e_b(q)^*e_a),\qquad
 e_a(q)=\frac{\widetilde e_a(q)}{|\widetilde e_a(q)|}.
\]
Since $\widetilde e_a(0)=e_a$, the fixed $C^2$ bound for $\nu$ allows
us to choose $\delta_2>0$, depending only on $\Omega$, so that
$|\widetilde e_a(q)|\ge1/2$ on $B_{\delta_2}(0)$. These formulas give
\[
 e_n(q)^*e_a(q)=0,\qquad e_a(q)^*e_b(q)=\delta_{ab},\qquad
 \sum_{a=1}^{n-1}\|e_a\|_{C^2(B _{\delta_2}(0))}\le C,
\]
with $C$ depending only on $\Omega$.

Extend $\Theta$ by $\Theta(q)e_a(q)=\mu_a e_a(q)$, keeping the weights
$\mu_a$ fixed, and define
\[
 b(q)=\pair{\Theta(q)}{H(\ubar)|_{E_q}},\qquad
 \eta(q)=\pair{\Theta(q)}{H(\varrho)|_{E_q}}.
\]
In this frame $[\Theta(q)]=\operatorname{diag}(\mu_1,\ldots,\mu_{n-1})
\in\mathscr S$. For $\psi=\ubar,\varrho$,
\eqref{eq:pairing-in-eigenbasis} gives
\[
  \pair{\Theta(q)}{H(\psi)|_{E_q}}
  =\sum_{a=1}^{n-1}\mu_a e_a(q)^*H(\psi)e_a(q) .
 \]
Since $\sum_a\mu_a=1$, the fixed $C^4$ norms of $\ubar,\varrho$ and the
$C^2$ frame bounds imply
\begin{equation}\label{eq:boundary-weight-bounds}
 \|b\|_{C^2(B _{\delta_2}(0))}+\|\eta\|_{C^2(B _{\delta_2}(0))}\le C.
\end{equation}
Choose $\delta_2\le\eta_*/[2(1+\|D\eta\|_{C^0})]$, using the fixed
bound in \eqref{eq:boundary-weight-bounds}. Then
\begin{equation}\label{eq:normal-weight-neighborhood}
 \eta(q)\ge\eta(0)-\|D\eta\|_{C^0}|q|\ge\eta_*/2
 \quad(q\in B _{\delta_2}(0)).
\end{equation}
On $B_{\delta_2}(0)\cap\overline\Omega$, define
\begin{equation}\label{eq:normal-auxiliary-functions}
 h_0=\frac{m_*-b}{\eta},\qquad
 V_i=v_i-h_0\nu_i\quad(1\le i\le4n),\qquad
 W=\nu\cdot Dv-h_0=\sum_i\nu_iV_i.
\end{equation}
Here $\nu_i$ are the real components of $\nu=D\varrho$.
The bound $m_*<b_*/2$, \eqref{eq:boundary-weight-bounds}, and
\eqref{eq:normal-weight-neighborhood} first give $|h_0|\le C$.
Differentiating $\eta h_0=m_*-b$ gives
\[
 (h_0)_r=-\frac{b_r+h_0\eta_r}{\eta},\qquad
 (h_0)_{rs}=-\frac{b_{rs}+h_0\eta_{rs}
                         +(h_0)_r\eta_s+(h_0)_s\eta_r}{\eta}.
\]
Equations \eqref{eq:boundary-weight-bounds},
\eqref{eq:normal-weight-neighborhood}, and \eqref{eq:c1} therefore give
\begin{equation}\label{eq:normal-aux-bounds}
 \|h_0\|_{C^2}\le C,\qquad |W|+|V|\le C,
 \qquad V=(V_1,\ldots,V_{4n})^T,
\end{equation}
with constants independent of $m_*$.
On $B_{\delta_2}(0)\cap\bd$, \eqref{eq:tangential-block},
\eqref{eq:boundary-minimum}, and \eqref{eq:normal-weight-neighborhood} imply
\[
 b+v_\nu\eta=\pair{\Theta(q)}{B(q)}\ge m(B(q))\ge m_*,
 \qquad W=\frac{b+v_\nu\eta-m_*}{\eta}\ge0,\qquad W(0)=0.
\]
Since $Dv=v_\nu\nu$ there,
\begin{equation}\label{eq:Vboundary}
 V_i=W\nu_i,\qquad |V|^2=W^2\quad\text{on }B _{\delta_2}(0)\cap\bd.
\end{equation}
In particular, $V(0)=0$ and $D|V|^2(0)=0$.

\subsection{The estimate at the minimizing point}
By \eqref{eq:first-derivative} and \eqref{eq:normal-aux-bounds},
\begin{equation}\label{eq:normal-vector-linearized}
 \cL V_i=g_{t,i}-\cL(\ubar_i+h_0\nu_i),\qquad
 |\cL V_i|\le C\left(1+\sum_rF^{rr}\right).
\end{equation}
The inequality follows from the fixed $C^2$ bound for
$\ubar_i+h_0\nu_i$ and \eqref{eq:linearized-coefficient-sum}.
Since $W=\sum_i\nu_iV_i$, \eqref{eq:linearized-product} gives
\begin{equation}\label{eq:LW-normal}
 \cL W=\sum_i\nu_i\cL V_i
   +\sum_{i,r,s}F^{rs}(\nu_i)_{rs}V_i
   +2\sum_{i,r,s}F^{rs}(\nu_i)_r(V_i)_s.
\end{equation}
The first two sums are bounded in absolute value by
$C(1+\sum_rF^{rr})$. By \eqref{eq:weighted-energy} and
\eqref{eq:weighted-young},
\[
 \begin{aligned}
 2\left|\sum_{i,r,s}F^{rs}(\nu_i)_r(V_i)_s\right|
 &\le\varepsilon\mathcal E(V)
      +\varepsilon^{-1}\sum_{i,r,s}F^{rs}(\nu_i)_r(\nu_i)_s\\
 &\le\varepsilon\mathcal E(V)
      +\varepsilon^{-1}|D\nu|^2\sum_rF^{rr}.
 \end{aligned}
\]
Therefore, for every $\varepsilon>0$,
\begin{equation}\label{eq:normal-barrier}
 \cL W\le\varepsilon\mathcal E(V)
      +C_\varepsilon\left(1+\sum_rF^{rr}\right).
\end{equation}
The bound for $|V|$ and \eqref{eq:normal-vector-linearized} allow us to
apply \eqref{eq:square-strong}, giving
\begin{equation}\label{eq:LQ-normal}
 \cL|V|^2\ge2\mathcal E(V)-C\left(1+\sum_rF^{rr}\right).
\end{equation}

Choose a fixed $0<\delta_3\le\min\{\delta_2,\delta\}$, where $\delta$ is the
collar width in Lemma~\ref{lem:collar}, and let $C_W$ be the bound for
$|W|$ in \eqref{eq:normal-aux-bounds}. Set
\[
 \theta_0=\frac1{2(1+C_W)},\qquad
 \Psi=W-\theta_0|V|^2+K_1\Phi+K_2|q|^2
 \quad\text{in }\Omega\cap B _{\delta_3}(0).
\]
On $\bd\cap B_{\delta_3}(0)$, \eqref{eq:Vboundary} gives
$\Psi\ge W-\theta_0W^2\ge W/2\ge0$.
Increase $K_2$ if necessary so that
 $K_2\delta_3^2\ge\sup_{\Omega\cap B_{\delta_3}(0)}(|W|+\theta_0|V|^2)$.
On $\partial B_{\delta_3}(0)\cap\Omega$ this gives
\[
 \Psi\ge-|W|-\theta_0|V|^2+K_2 \delta_3^2\ge0.
\]
Also $\cL|q|^2=2\sum_rF^{rr}-g_{t,z}|q|^2\le2\sum_rF^{rr}$.
By \eqref{eq:collar-ineq}, \eqref{eq:normal-barrier} with
$\varepsilon=\theta_0$, and \eqref{eq:LQ-normal},
\[
  \cL\Psi
  \le \bigl[-K_1\varepsilon_0+C_{\theta_0}+\theta_0 C+2K_2\bigr]
       \left(1+\sum_rF^{rr}\right)-\theta_0\mathcal E(V)\le0,
 \]
 provided $K_1\ge(C_{\theta_0}+\theta_0 C+2K_2)/\varepsilon_0$.
The maximum principle gives $\Psi\ge0$. Since $\Psi(0)=0$ and
$D|V|^2(0)=0$,
\[
 0\ge\Psi_\nu(0)=W_\nu(0)+K_1\Phi_\nu(0).
\]
By \eqref{eq:iterated-directional-derivative} and
\eqref{eq:distance-identities}, $D_\nu(D_\nu v)=v_{\nu\nu}$. Hence
\[
 \begin{aligned}
 W_\nu(0)
 &= v _{\nu\nu}(0)-(h_0)_\nu(0),\\
 \Phi_\nu(0)&=v_\nu(0)+(\eta_0-2Nd(0))d_\nu(0)=v_\nu(0)-\eta_0.
 \end{aligned}
\]
Equations \eqref{eq:c1}, \eqref{eq:normal-aux-bounds}, and
\eqref{eq:boundary-block-bounds} now give
\begin{equation}\label{eq:normal-at-p}
 u_{\nu\nu}(0)
 \le K_1(\eta_0-v_\nu(0))+(h_0)_\nu(0)+\ubar_{\nu\nu}(0)\le C,
 \qquad a_0(0)\le C_{a_0}.
\end{equation}

\subsection{The uniform boundary estimate}
Choose $C_B\ge1$ with $\|B\|_{\mathrm{op}}\le C_B$ on $\bd$.
By continuity, \eqref{eq:equation} holds on $\bd$, where $u=\varphi$.
Set
\[
 c_*:=\min_{q\in\bd,\,0\le t\le1}
              \binom nk g_t(q,\varphi(q))^k>0.
\]
By \eqref{eq:newton-positivity} and \eqref{eq:block},
\begin{equation}\label{eq:normal-equation}
 \begin{aligned}
 c_*\le\sigma_k(H(u)(0))
 &=\sigma_k(B(0))+a_0(0)\sigma_{k-1}(B(0))
       -\beta(0)^*T_{k-2}(B(0))\beta(0)\\
 &\le\sigma_k(B(0))+a_0(0)\sigma_{k-1}(B(0)).
 \end{aligned}
\end{equation}
Equation~\eqref{eq:newton-boundary} gives
\[
 \sigma_k(B(0))
 \le\frac{n-k}{k(n-1)}\sigma_1(B(0))\sigma_{k-1}(B(0))
 \le\frac{n-k}{k}C_B\sigma_{k-1}(B(0)).
\]
With $C_{a_0}$ from \eqref{eq:normal-at-p}, set
$C_0=C_{a_0}+(n-k)C_B/k$.
Then \eqref{eq:normal-equation} and the upper bound in
\eqref{eq:margin-bounds} give
\[
 c_*\le C_0\sigma_{k-1}(B(0))\le C_0C_1m_*,\qquad
 C_1=(k-1)\binom{n-1}{k-1}(2C_B)^{k-2}.
\]
Thus $m_*\ge c_*/(C_0C_1)$. Together with the case $m_*\ge b_*/2$,
this proves a uniform positive lower bound for $m_*$.
By \eqref{eq:boundary-minimum} and the lower bound in
\eqref{eq:margin-bounds},
\begin{equation}\label{eq:uniform-tangential-margin}
 \sigma_{k-1}(B(q))\ge\binom{n-1}{k-1}m(B(q))^{k-1}
 \ge\binom{n-1}{k-1}m_*^{k-1}\ge c>0
 \quad(q\in\bd).
\end{equation}
To bound the numerator in the block formula, use the eigenvalue bound
$|\lambda_a(B)|\le C_B$ from \eqref{eq:boundary-block-bounds} to obtain
\[
 |\sigma_k(B)|\le\binom{n-1}{k}C_B^k\quad(k<n),\qquad
 \|T_{k-2}(B)\|_{\mathrm{op}}
 \le\binom{n-2}{k-2}C_B^{k-2},
\]
by the definitions of the elementary symmetric functions and
\eqref{eq:newton-spectral-formula}.
When $k=n$, $\sigma_n(B)=0$.
We can now solve \eqref{eq:block} at every boundary point:
\[
 \begin{aligned}
 a_0&=\frac{\sigma_k(H(u))-\sigma_k(B)
                  +\beta^*T_{k-2}(B)\beta}{\sigma_{k-1}(B)},\\
 |a_0|&\le
 \frac{\binom nk(g^*)^k+|\sigma_k(B)|
                 +\|T_{k-2}(B)\|_{\mathrm{op}}|\beta|^2}{c}\le C.
 \end{aligned}
\]
Equation~\eqref{eq:boundary-block-bounds} now gives $|u_{\nu\nu}|\le4|a_0|+C\le C$.
Together with Proposition~\ref{prop:partial-boundary-C2}, this proves
the following result.

\begin{proposition}[Uniform boundary $C^2$ estimate]
\label{prop:boundary-C2}
Every smooth admissible solution of \eqref{eq:equation} satisfies
\[
 \sup_{\bd}|D^2u|\le C,\qquad
 \|u\|_{C^2(\overline\Omega)}\le C.
\]
\end{proposition}

\section{Proof of Theorem~\ref{thm:main} and its extension}\label{sec:existence}

The continuity method starts from the known solution $\ubar$ at $t=0$
and varies the parameter in \eqref{eq:path} until $t=1$.
We must show that solutions persist for nearby parameters and that
limits of solutions remain admissible solutions. The linearized
Dirichlet problem gives the first assertion. The uniform estimates
from Sections~\ref{sec:lower}--\ref{sec:boundary}, together with the
regularity below, give the second.

\subsection{Uniform ellipticity and the \texorpdfstring{$C^{2,\alpha}$}{C2,alpha} estimate}
\label{sec:uniform-regularity}
For every smooth admissible solution $u_t$ of \eqref{eq:equation},
Proposition~\ref{prop:boundary-C2} gives
\begin{equation}\label{eq:global-c2}
 \|u_t\|_{C^2(\closure\Omega)}\le C_1.
\end{equation}
Together with $g_t(q,u_t)\ge g_*>0$ and
\eqref{eq:cone-closure-boundary}, this yields, after enlarging $C_1$,
\begin{equation}\label{eq:uniform-admissible-set}
 H(u_t)(\closure\Omega)\subset
 \mathcal K:=\left\{A\in\overline{\Gamma}_k:|A|\le C_1,\quad
 \sigma_k(A)\ge\binom nk g_*^k\right\}\Subset\Gamma_k.
\end{equation}
By compactness and Lemma~\ref{lem:concavity-ellipticity},
\begin{equation}\label{eq:uniform-ellipticity}
 \lambda I_{4n}\le\bigl(F^{rs}(A)\bigr)_{r,s=1}^{4n}
 \le\Lambda I_{4n}\qquad(A\in\mathcal K),
 \qquad 0<\lambda\le\Lambda,
\end{equation}
with constants independent of $t$ and $u_t$.

To apply real uniformly elliptic regularity theory, define on all real
symmetric matrices
\begin{equation}\label{eq:uniformly-elliptic-extension}
 \widehat F(S):=\inf_{A\in\mathcal K}DF(A)[\Proj(S)]
 =\inf_{A\in\mathcal K}\sum_{r,s=1}^{4n}F^{rs}(A)S_{rs},
 \qquad S\in\operatorname{Sym}(4n,\RR).
\end{equation}
This infimum of linear functions is concave and is defined even when
$\Proj(S)\notin\Gamma_k$. By \eqref{eq:supporting-homogeneity} and
\eqref{eq:hessian-projection},
\begin{equation}\label{eq:extension-agreement}
 \widehat F(S)=F(\Proj(S))\quad\text{if }\Proj(S)\in\mathcal K,
 \qquad \widehat F(D^2u_t)=g_t(q,u_t).
\end{equation}
Applying \eqref{eq:uniform-ellipticity} under the infimum gives, for
$S,C\in\operatorname{Sym}(4n,\RR)$ with $C\ge0$,
\begin{equation}\label{eq:extension-ellipticity}
 \lambda\trR C\le\widehat F(S+C)-\widehat F(S)
 \le\Lambda\trR C.
\end{equation}
Thus $\widehat F$ is uniformly elliptic without requiring differentiability.

By \eqref{eq:c1} and the bounds for $g_t$,
$\|g_t(\cdot,u_t)\|_{C^{\alpha_0}(\closure\Omega)}\le C$
for any fixed $0<\alpha_0<1$.
Concavity and \eqref{eq:extension-ellipticity} give the global
$C^{2,\alpha}$ estimates for the model equation $\widehat F(D^2u)=0$
in \cite[Chapters~6 and~9]{CaffarelliCabre1995}.
Hence \cite[Theorem~1.4(b)]{SilvestreSirakov2014}, applied to
\eqref{eq:extension-agreement} with boundary data $\varphi$, gives
some $\alpha\in(0,\alpha_0)$ such that
\begin{equation}\label{eq:c2alpha}
 \|u_t\|_{C^{2,\alpha}(\closure\Omega)}\le C.
\end{equation}

\subsection{The continuity method}\label{sec:continuity-method}
\begin{proof}[Proof of Theorem~\ref{thm:main}]
Set
\[
 \mathcal T=\{t\in[0,1]:\eqref{eq:equation}\text{ has a smooth admissible solution}\}.
\]
Since $u_0=\ubar$, we have $0\in\mathcal T$.

\emph{Openness.}
Set
\[
 \mathcal X=\{v\in C^{2,\alpha}(\closure\Omega):v|_{\bd}=0\},
 \qquad \mathcal Y=C^\alpha(\closure\Omega).
\]
On the open set of $v\in\mathcal X$ satisfying
$H(\ubar+v)\in\Gamma_k$ on $\closure\Omega$, define
\[
 \mathcal F(t,v)=F(H(\ubar+v))-g_t(q,\ubar+v).
\]
At a solution, \eqref{eq:linearized-operator} gives
\[
 D_v\mathcal F(t,u-\ubar)=\cL:\mathcal X\longrightarrow\mathcal Y.
\]
By \eqref{eq:uniform-ellipticity}, $g_{t,z}\ge0$, the maximum principle,
and \cite[Theorem~6.14]{GT2001}, this map is an isomorphism.
The implicit function theorem gives admissible $C^{2,\alpha}$ solutions
for nearby parameters. Smoothness of the original equation and uniform
ellipticity near each solution's Hessian range give
$C^\infty(\closure\Omega)$ regularity by Schauder bootstrapping
\cite[Theorem~6.19 and Lemma~17.16]{GT2001}.
Thus $\mathcal T$ is open in $[0,1]$.

\emph{Closedness.}
If $t_j\in\mathcal T$ and $t_j\to t_\infty$, \eqref{eq:c2alpha} gives,
after passing to a subsequence,
\[
 u_{t_j}\longrightarrow u_\infty\quad\text{in }C^2(\closure\Omega),
 \qquad u_\infty\in C^{2,\alpha}(\closure\Omega).
\]
By \eqref{eq:equation} and \eqref{eq:uniform-admissible-set},
\[
 F(H(u_\infty))=g_{t_\infty}(q,u_\infty),\qquad
 u_\infty|_{\bd}=\varphi,\qquad
 H(u_\infty)(\closure\Omega)\subset\mathcal K\Subset\Gamma_k.
\]
The same bootstrap gives $u_\infty\in C^\infty(\closure\Omega)$, hence
$t_\infty\in\mathcal T$. Therefore $\mathcal T=[0,1]$.
Lemma~\ref{lem:comparison-C0}, applied to two solutions in both orders,
gives uniqueness and also yields $u\ge\ubar$.
\end{proof}

\subsection{Flat quaternionic quotients}\label{sec:quotient-extension}

The local lifting argument in \cite[Sections~3--5 and~8]{MaZhu2026}
gives the following extension.

\begin{theorem}[Flat quaternionic quotients]\label{thm:quotient}
Let $2\le k\le n$, let $\Gamma\subset(\HH^n,+)$ be a discrete
subgroup, and equip $X=\HH^n/\Gamma$ with the induced flat metric
and quaternionic structures. Let $\Omega\Subset X$ be a domain with
nonempty smooth boundary, and suppose that
\[
 \varphi\in C^\infty(\bd),\qquad
 f\in C^\infty(\closure\Omega\times\RR),\qquad f>0,\qquad f_z\ge0.
\]
Assume that a real-valued $\ubar\in C^\infty(\closure\Omega)$ satisfies
\[
 H(\ubar)\in\Gamma_k\quad\text{in }\Omega,\qquad
 \sigma_k(H(\ubar))\ge f(q,\ubar)\quad\text{in }\Omega,\qquad
 \ubar=\varphi\quad\text{on }\bd,
\]
where the normalized quaternionic Hessian is defined in local lifts to
$\HH^n$. Then the Dirichlet problem \eqref{eq:dirichlet} has a unique
$k$-admissible solution $u\in C^\infty(\closure\Omega)$, and $u\ge\ubar$.
\end{theorem}

\begin{proof}
Let $\pi:\HH^n\to X$ be the quotient map, set
$\widetilde\Omega=\pi^{-1}(\Omega)$, and write $\widetilde u=u\circ\pi$.
Since $\widetilde u$ is $\Gamma$-periodic,
\eqref{eq:unitary-covariance} with $U=I_n$ gives
\begin{equation}\label{eq:quotient-hessian}
 H(\widetilde u)(q+\gamma)=H(\widetilde u)(q)
 \quad(\gamma\in\Gamma),\qquad
 H(u)(\pi(q)):=H(\widetilde u)(q).
\end{equation}
Thus $H(u)$ and $k$-admissibility are independent of the lift.
The path \eqref{eq:path} is defined on $\Omega$, and its equation lifts to
\[
 F(H(\widetilde u))(q)=g_t(\pi(q),\widetilde u(q))
 \quad\text{in }\widetilde\Omega,
 \qquad \widetilde u=\varphi\circ\pi
 \quad\text{on }\partial\widetilde\Omega.
\]
The lifted subsolution is $\ubar\circ\pi$.

The harmonic extension of $\varphi$ exists on the compact manifold with
boundary $\closure\Omega$. The proof of Lemma~\ref{lem:comparison-C0}
therefore gives the same $C^0$ and boundary gradient bounds.
Constant real vector fields on $\HH^n$ descend to a global parallel
frame on $X$. For any parallel real unit field $e$, equations
\eqref{eq:first-derivative}, \eqref{eq:subsolution-test}, and
\eqref{eq:second-linearized} hold with derivatives taken in this frame.
Since $\ubar-C_{\ubar}$ is globally defined, the maximum principle on
$\Omega$ gives \eqref{eq:c1} and
\eqref{eq:second-derivative-reduction}.

For each $p\in\bd$, choose $\widetilde p\in\pi^{-1}(p)$ and $\delta_p>0$
so that $\pi$ is injective on $B_{\delta_p}(\widetilde p)$. Its restriction
\[
 \pi:B_{ \delta_p}(\widetilde p)\cap\widetilde\Omega
 \longrightarrow \pi(B_{ \delta_p}(\widetilde p))\cap\Omega
\]
is a diffeomorphism preserving the metric and quaternionic structures.
After translating $\widetilde p$ to the origin, the barrier of
Lemma~\ref{lem:collar} and the tangential and mixed derivative estimates
in Subsection~\ref{sec:tangential-estimates} apply on these bounded
local lifts. The artificial spherical boundaries are controlled by
\eqref{eq:c1}, as in the original estimates.
By \eqref{eq:subspace-matrix} and \eqref{eq:margin-continuity},
$q\mapsto m(H(u)(q)|_{E_q})$ is a continuous function on $\bd$,
independent of the chosen quaternionic orthonormal frame.
It attains its minimum on $\bd$, so the argument in
Section~\ref{sec:boundary} also applies in a lift of the minimizing point.
A finite cover of $\bd$ makes all boundary constants uniform.
Combining the resulting boundary $C^2$ bound with
\eqref{eq:second-derivative-reduction} and $\Delta u>0$ gives
\eqref{eq:global-c2} on $\closure\Omega$.

Equations \eqref{eq:uniform-admissible-set} and
\eqref{eq:uniform-ellipticity} hold on $\closure\Omega$. Applying the
argument of Subsection~\ref{sec:uniform-regularity} in finitely
many local lifts gives \eqref{eq:c2alpha}.
Linear Dirichlet Schauder theory on $\closure\Omega$ gives the same
isomorphism for $\cL$, so the openness, closedness, and uniqueness
arguments in the proof of Theorem~\ref{thm:main} apply to
\eqref{eq:path}. This proves the result.
\end{proof}

\section{Geometric consequences and examples}\label{sec:geometry}
\subsection{Formulation using defining functions}\label{sec:levi-invariance}
For any smooth defining function $\rho$ with $\rho<0$ in $\Omega$,
\eqref{eq:boundary-gradient} gives $D\rho=\rho_\nu\nu$ and
$\rho_\nu=|D\rho|>0$ on $\bd$. Equation \eqref{eq:boundary-restriction} yields
\begin{equation}\label{eq:levi-invariance}
 H(\rho)|_{E_q}=|D\rho|R_q,\qquad
 \sigma_j(H(\rho)|_{E_q})=|D\rho|^j\sigma_j(R_q)
 \quad(1\le j\le n-1).
\end{equation}
Since $|D\rho|>0$ on $\bd$, condition \eqref{eq:levi-condition} is equivalent to
$H(\rho)(q)|_{E_q}\in\Gamma_{k-1}^{(n-1)}$ for all $q\in\bd$,
for any such defining function $\rho$.

\subsection{Construction of subsolutions}\label{sec:subsolution-construction}
\begin{proof}[Proof of Corollary~\ref{cor:construction}]

\textbf{Construction of an admissible defining function.}
We first construct a defining function which is admissible in $\Omega_\delta:=\{q\in\Omega:d(q)<\delta\}$. Let $\kappa_0>0$ be chosen below and set $\rho_0=\varrho+\kappa_0\varrho^2$. On $\partial\Omega$,
\[
 H(\rho_0)=H(\varrho)+2\kappa_0\Proj(\nu\otimes\nu).
\]
 At each boundary point, choose $U=(e_1,\ldots,e_n)$ as in
\eqref{eq:unitary-columns}, with $e_n=\iota^{-1}(\nu)$.
Equations \eqref{eq:projected-rank-one} and
\eqref{eq:orthonormal-basis-change} give
\[
 \Proj(\nu\otimes\nu)=\tfrac14e_ne_n^*,\qquad
 U^*\Proj(\nu\otimes\nu)U=\tfrac14E_{nn}.
\] The polynomial identity
\eqref{eq:rank-one-sigma}, which does not require $H(\varrho)$ to be
admissible, gives, for $1\le j\le k$,
\begin{equation}\label{eq:defining-normal-completion}
 \sigma_j(H(\rho_0))
 =\sigma_j(H(\varrho))+\frac{\kappa_0}{2}\sigma_{j-1}(R_q)
 \quad\text{on }\bd.
\end{equation}

Since $R_q\in\Gamma_{k-1}^{(n-1)}$ and $\bd$ is compact,
\[
 c_R:=\min_{q\in\bd,\,1\le j\le k}\sigma_{j-1}(R_q)>0,
 \qquad C_R:=\max_{q\in\bd,\,1\le j\le k}|\sigma_j(H(\varrho)(q))|<\infty.
\]
Choose $\kappa_0>2C_R/c_R$. By \eqref{eq:defining-normal-completion}
and continuity, after reducing $\delta$,
\[
 H(\rho_0)\in\Gamma_k,\qquad 1+2\kappa_0\varrho>0\quad\text{on }\overline{\Omega_\delta}.
\]

We adapt the defining-function construction in \cite[Theorem~3.1]{Li2004}.
In the rest of this construction, all level sets are taken in
$\overline\Omega$.
Extend $\rho_0$ smoothly and negatively into $\Omega$, preserving it on
$\overline{\Omega_\delta}$. Since $\rho_0<0$ on the compact set
$\{d\ge\delta\}\Subset\Omega$, we can then choose $a>0$ such that
\[
 0<2a<-\max_{\{d\ge\delta\}}\rho_0,
 \qquad \{\rho_0\ge-2a\}\subset\Omega_\delta\cup\bd.
\]
Choose a smooth convex function $\theta$ with
\[
 \theta(s)=s\ \text{near }0,\qquad
 \theta'(s)=0\ \text{for }s\le-2a,\qquad
 \theta'(s)>0\ \text{for }-2a<s\le0.
\]
On $\{\rho_0>-2a\}$ we have $H(\rho_0)\in\Gamma_k$,
whereas $\theta'(\rho_0)=\theta''(\rho_0)=0$ elsewhere.
Thus, for $\psi=\theta(\rho_0)$, equations \eqref{eq:hessian-chain},
\eqref{eq:projected-rank-one}, and \eqref{eq:positive-perturbation} give
\[
 H(\psi)=\theta'(\rho_0)H(\rho_0)
 +\theta''(\rho_0)\Proj(D\rho_0\otimes D\rho_0)
 \in\overline{\Gamma}_k,
 \qquad H(\psi)\in\Gamma_k\ \text{on }\{\rho_0\ge-a\}.
\]
Choose $\chi_0\in C_c^\infty(\Omega)$ with $0\le\chi_0\le1$ and
 $\chi_0\equiv1$ on a neighborhood of $\{\rho_0\le-a\}$ in $\Omega$.
Set $\rho=\psi+\varepsilon\chi_0|q|^2$.
By \eqref{eq:normalization} and \eqref{eq:concavity-tools},
\[
 H(\rho)=H(\psi)+2\varepsilon I_n\in\Gamma_k
 \quad\text{on }\{\rho_0\le-a\}.
\]
On the compact set $\{\rho_0\ge-a\}$, the matrices $H(\psi)$ stay in a
compact subset of $\Gamma_k$, so $H(\rho)\in\Gamma_k$ there for
sufficiently small $\varepsilon>0$.
Decreasing $\varepsilon$ also ensures $\rho<0$ on $\operatorname{supp}\chi_0$.
Since $\rho=\rho_0$ near $\bd$, this gives a smooth defining function
with $H(\rho)\in\Gamma_k$ on $\overline\Omega$.

 \noindent\textbf{Construction of a subsolution.}
Extend $\varphi$ to a smooth real function $\widetilde\varphi$ on
$\closure\Omega$ and put $\ubar_\kappa=\widetilde\varphi+\kappa\rho$,
where $\kappa>0$ will be chosen below.
Since $H(\rho)(\closure\Omega)$ is compact in $\Gamma_k$,
for all sufficiently large $\kappa$ the matrices
\[
 \kappa^{-1}H(\ubar_\kappa)
 =H(\rho)+\kappa^{-1}H(\widetilde\varphi)
\]
stay in a fixed compact set $\mathcal K_\rho\Subset\Gamma_k$. Thus
$c_1:=\min_{\mathcal K_\rho}F>0$, and positivity and homogeneity of $F$ give
\[
 H(\ubar_\kappa)\in\Gamma_k,\qquad
 F(H(\ubar_\kappa))\ge c_1\kappa\quad\text{on }\closure\Omega.
\]
Set $\varphi_{\max}=\max_{\bd}\varphi$. Then
\[
 \Delta\ubar_\kappa>0\ \Longrightarrow\ \ubar_\kappa\le \varphi_{\max},
 \qquad
 g(q,\ubar_\kappa)\le g(q,\varphi_{\max})
 \le\max_{\overline\Omega}g(q,\varphi_{\max}).
\]
Choose $\kappa$ still larger so that $c_1\kappa$ exceeds the last bound.
Then $F(H(\ubar_\kappa))\ge g(q,\ubar_\kappa)$ and
$\ubar_\kappa=\varphi$ on $\bd$, so Theorem~\ref{thm:main} applies.
\end{proof}

\begin{proposition}\label{prop:quotient-defining-function}
Let $\Omega\Subset X=\HH^n/\Gamma$ have smooth boundary satisfying
\eqref{eq:levi-condition}. Suppose that a real
 $\chi\in C^\infty(\overline\Omega)$ satisfies
$H(\chi)\in\Gamma_k$ on $\overline\Omega$.
 Then $\Omega$ has a smooth defining function $\rho$ with
$H(\rho)\in\Gamma_k$ on $\overline\Omega$.
Consequently, the conclusion of Corollary~\ref{cor:construction}
holds on $\Omega$.
\end{proposition}
\begin{proof}
Use $\rho_0,\psi,\chi_0$ from the preceding construction, with all level
sets taken in $\overline\Omega$. Subtract a constant from $\chi$ so that
$\chi\le0$ and put $\rho=\psi+\varepsilon\chi_0\chi$.
 Since $\chi_0\equiv1$ near $\{\rho_0\le-a\}$,
\eqref{eq:concavity-tools} gives
\[
 H(\rho)=H(\psi)+\varepsilon H(\chi)\in\Gamma_k
 \quad\text{on }\{\rho_0\le-a\}.
\]
 On the compact set $\{\rho_0\ge-a\}$, $H(\psi)$ stays in a compact
subset of $\Gamma_k$, so the same conclusion holds for sufficiently
small $\varepsilon>0$. Also $\rho<0$ in $\Omega$ and $\rho=\rho_0$
near $\bd$. The subsolution construction above and
Theorem~\ref{thm:quotient} complete the proof.
\end{proof}
The hypothesis on $\chi$ imposes no boundary values or inequality involving $f$.
The periodic domain in \eqref{eq:periodic-quartic-domain} satisfies these
hypotheses for $2\le k<n$, although it is not strictly quaternionic pseudoconvex.

\subsection{Constant boundary values}
\begin{proof}[Proof of Corollary~\ref{cor:necessity}]
Sufficiency follows from Corollary~\ref{cor:construction}, so it remains
to prove necessity.
Let $u\in C^2(\overline\Omega)$ solve the equation with boundary value
$c_0$. By \eqref{eq:boundary-admissibility} at $t=1$,
$H(u)\in\Gamma_k$ on $\overline\Omega$. Equation \eqref{eq:normalization},
the strong maximum principle, and Hopf's lemma give
\[
 \Delta u=4\sigma_1(H(u))>0,\qquad u<c_0\ \text{in }\Omega,
 \qquad u_\nu>0\ \text{on }\bd.
\]
By \eqref{eq:boundary-restriction} applied to $u-c_0$ and
Lemma~\ref{lem:block-identity},
\[
 u_\nu R_q=H(u)|_{E_q}\in\Gamma_{k-1}^{(n-1)}.
\]
Dividing by $u_\nu>0$ proves \eqref{eq:levi-condition}. The same argument
applies to an admissible subsolution with constant boundary value.

For the last assertion, \eqref{eq:levi-condition} immediately implies
\eqref{eq:scalar-levi-condition}, so only the converse remains.
Assume that $\bd$ is connected and \eqref{eq:scalar-levi-condition}
holds. We use the connectedness argument of
\cite[Remark~1.2]{CNS1985-III}.
Fix $c\in\Omega$ and choose $p\in\bd$ such that
$|p-c|=\max_{q\in\closure\Omega}|q-c|$.
Then $\Omega\subset B_{|p-c|}(c)$, and tangency at $p$ gives
$\iota(p-c)=|p-c|\nu(p)$.
For a unit $\tau\in T_p\bd$, choose a smooth curve
$\gamma$ in $\iota(\bd)$ with $\gamma(0)=\iota(p)$ and
$\gamma'(0)=\tau$. Differentiating
$\varrho\circ\iota^{-1}\circ\gamma=0$ twice gives
\begin{equation}\label{eq:boundary-curve-acceleration}
 0=\left.\frac{\mathrm d^2}{\mathrm d\theta^2}\right|_{\theta=0}
       (\varrho\circ\iota^{-1}\circ\gamma)(\theta)
 =D^2\varrho_p(\tau,\tau)+\nu(p)\cdot\gamma''(0).
\end{equation}
Since $|\gamma(\theta)-\iota(c)|^2$ has a maximum at $\theta=0$,
\eqref{eq:boundary-curve-acceleration} yields
\begin{equation}\label{eq:enclosing-ball-curvature}
 \begin{aligned}
 0&\ge\left.\frac{\mathrm d^2}{\mathrm d\theta^2}\right|_{\theta=0}
      |\gamma(\theta)-\iota(c)|^2\\
  &=2|\gamma'(0)|^2
     +2\bigl(\gamma(0)-\iota(c)\bigr)\cdot\gamma''(0)\\
  &=2+2|p-c|\nu(p)\cdot\gamma''(0)
   =2-2|p-c|D^2\varrho_p(\tau,\tau).
 \end{aligned}
\end{equation}
For every unit $\xi\in E_p$, \eqref{eq:quaternionic-structure-orthogonality}
and \eqref{eq:tangent-splitting} show that
$Q\iota(\xi)$ is a unit tangent vector at $p$ for each
$Q\in\mathscr Q$.
Applying \eqref{eq:enclosing-ball-curvature} to these vectors and using
\eqref{eq:hessian-projection} and \eqref{eq:quadratic-average}, we obtain
\begin{equation}\label{eq:positive-levi-contact-point}
 \xi^*R_p\xi
 =\frac14\sum_{Q\in\mathscr Q}
        D^2\varrho_p(Q\iota(\xi),Q\iota(\xi))\ge |p-c|^{-1}>0.
\end{equation}
Equation \eqref{eq:cone-closure-boundary}, applied to $\mathcal H_{n-1}$
and $\sigma_{k-1}$, states that
\[
 \partial\Gamma_{k-1}^{(n-1)}
 =\{A\in\overline{\Gamma}_{k-1}^{(n-1)}:\sigma_{k-1}(A)=0\}.
\]
Since \eqref{eq:scalar-levi-condition} gives $\sigma_{k-1}(R_q)>0$
at every $q\in\bd$, none of the forms $R_q$ lies on this cone boundary.
Consequently,
\[
 \mathcal V:=\{q\in\bd:R_q\in\Gamma_{k-1}^{(n-1)}\}
 =\{q\in\bd:R_q\in\overline{\Gamma}_{k-1}^{(n-1)}\}.
\]
By continuity, the first description makes $\mathcal V$ open in $\bd$
and the second makes it closed.
Equation \eqref{eq:positive-levi-contact-point} gives $p\in\mathcal V$.
Since $\bd$ is connected, $\mathcal V=\bd$,
which proves \eqref{eq:levi-condition}.
\end{proof}

\subsection{Examples}
For $2\le k<n$, the two examples below satisfy
\eqref{eq:levi-condition} but are not strictly quaternionic pseudoconvex.
The first is a bounded domain in $\HH^n$. The second adapts the
construction to a periodic quotient and uses
 Proposition~ \ref{prop:quotient-defining-function} to construct subsolutions.
Both admit smooth solutions for arbitrary smooth boundary values under
the stated assumptions on $f$.

\medskip
\Needspace{6\baselineskip}
\noindent\textbf{Example 1. A strictly quaternionic $(k-1)$-convex domain.}
\par\nobreak\smallskip\nobreak
\noindent
Let $2\le k<n$ and
$0<\varepsilon<(n-k)/(k-1)$. Put
\begin{equation}\label{eq:quartic-domain}
 \rho_\varepsilon=\sum_{a=1}^{n-1}|q_a|^2-\varepsilon|q_n|^2+|q_n|^4-1,
 \qquad \Omega_\varepsilon=\{\rho_\varepsilon<0\}.
\end{equation}
The negative quadratic term also appears in
\cite[Example~3.3]{AhagCzyzHed2018}, where intersecting a quadratic
sublevel set with a polydisc gives a non-pseudoconvex complex domain
admitting a negative $m$-subharmonic exhaustion.
Here the quartic term produces a smooth bounded domain, and we verify
the strict boundary condition \eqref{eq:levi-condition} despite a
negative quaternionic Levi eigenvalue.

By \eqref{eq:quartic-domain}, $\rho_\varepsilon(q)\to+\infty$ as
$|q|\to\infty$, and on
$\{\rho_\varepsilon=0\}$,
\begin{equation}\label{eq:quartic-domain-regularity}
 D\rho_\varepsilon(q)\cdot\iota(q)
 =2\sum_{a<n}|q_a|^2-2\varepsilon|q_n|^2+4|q_n|^4
 =2+2|q_n|^4>0.
\end{equation}
Equation \eqref{eq:quartic-domain-regularity} gives $D\rho_\varepsilon\ne0$
on $\{\rho_\varepsilon=0\}$, so the implicit function theorem shows that
the boundary is smooth. Since $\rho_\varepsilon(0)=-1$ and every
crossing along a ray is from negative to positive values, each ray meets
the boundary exactly once. Hence the domain is bounded and star-shaped,
with connected boundary.

Fix $q\in\partial\Omega_\varepsilon$ and let
$\epsilon_n=(0,\ldots,0,1)^T\in\HH^n$ be the last standard basis vector
from Definition~\ref{def:real-coordinates}. Direct computation from
\eqref{eq:hessian-definition} gives
\begin{equation}\label{eq:quartic-hessian}
 \begin{aligned}
 H(\rho_\varepsilon)(q)
 &=\operatorname{diag}(2,\ldots,2,6|q_n|^2-2\varepsilon)\\
 &=2I_n+(6|q_n|^2-2\varepsilon-2)\epsilon_n\epsilon_n^*.
 \end{aligned}
\end{equation}
Choose a quaternionic unitary matrix $U=(e_1,\ldots,e_n)$ whose columns
are expressed in the standard basis of $\HH^n$, with
$e_n=\iota^{-1}(\nu(q))$. Then $e_1,\ldots,e_{n-1}$ form an orthonormal
basis of $E_q$. Set
\[
 V=(e_1,\ldots,e_{n-1})\in\HH^{n\times(n-1)},\qquad
 V^*V=I_{n-1},\qquad \zeta=V^*\epsilon_n\in\HH^{n-1}.
\]
Equation \eqref{eq:orthonormal-basis-change} and unitarity give
\begin{equation}\label{eq:quartic-coordinate-decomposition}
 \begin{gathered}
 U^*\epsilon_n=\begin{pmatrix}\zeta\\e_n^*\epsilon_n\end{pmatrix},
 \qquad \epsilon_n=V\zeta+e_n(e_n^*\epsilon_n),\\
 |\zeta|^2+|e_n^*\epsilon_n|^2=1,\qquad |\zeta|\le1.
 \end{gathered}
\end{equation}
Thus $\zeta\in\HH^{n-1}$ is the coordinate column in the basis
$e_1,\ldots,e_{n-1}$ of the component $V\zeta\in E_q\subset\HH^n$.
By \eqref{eq:subspace-matrix} and \eqref{eq:quartic-hessian}, the restriction
$H(\rho_\varepsilon)(q)|_{E_q}$ is represented in this basis by
\begin{equation}\label{eq:quartic-tangential-block}
 \begin{aligned}
 B:=V^*H(\rho_\varepsilon)(q)V
 &=2V^*V+(6|q_n|^2-2\varepsilon-2)
       (V^*\epsilon_n)(V^*\epsilon_n)^*\\
 &=2I_{n-1}+(6|q_n|^2-2\varepsilon-2)\zeta\zeta^*.
 \end{aligned}
\end{equation}
 For $\zeta\ne0$, let
$\zeta^\perp=\{\xi\in\HH^{n-1}:\zeta^*\xi=0\}$.
 On $\HH^{n-1}=\zeta\HH\oplus\zeta^\perp$, we have
\[
 (\zeta\zeta^*)\zeta=\zeta |\zeta|^2,\qquad
 (\zeta\zeta^*) \xi=0\quad( \xi\in\zeta^\perp) .
\]
 Thus \eqref{eq:quartic-tangential-block} has eigenvalues
$2$ on $\zeta^\perp$ and $\mu$ on $\zeta\HH$,
 where the bound uses \eqref{eq:quartic-coordinate-decomposition}:
\begin{equation}\label{eq:quartic-tangential-eigenvalue}
 \begin{aligned}
 \mu&=2+(6|q_n|^2-2\varepsilon-2)|\zeta|^2\\
 &=2(1-|\zeta|^2)+(6|q_n|^2-2\varepsilon)|\zeta|^2
 \ge-2\varepsilon.
 \end{aligned}
\end{equation}
If $\zeta=0$, then $\zeta^\perp=\HH^{n-1}$, $B=2I_{n-1}$ and $\mu=2$.
In either case, the eigenvalues of $B$ are
$\underbrace{2,\ldots,2}_{n-2},\mu$.
Using \eqref{eq:quartic-tangential-eigenvalue}, for $1\le j\le k-1$ we obtain
\begin{equation}\label{eq:quartic-tangential-sigma}
 \begin{aligned}
 \sigma_j(B)
 &=2^j\binom{n-2}{j}+2^{j-1}\mu\binom{n-2}{j-1}\\
 &\ge2^j\left[\binom{n-2}{j}
                   -\varepsilon\binom{n-2}{j-1}\right]>0,
 \end{aligned}
\end{equation}
because $\binom{n-2}{j}/\binom{n-2}{j-1}
=(n-1-j)/j\ge(n-k)/(k-1)>\varepsilon$.
Together with \eqref{eq:levi-invariance},
\eqref{eq:quartic-tangential-sigma} proves strict quaternionic
$(k-1)$-convexity.

At $p=(1,0,\ldots,0)$, $|D\rho_\varepsilon|=2$ and
\eqref{eq:quartic-hessian} gives
\begin{equation}\label{eq:quartic-negative-tangential-eigenvalue}
 H(\rho_\varepsilon)|_{E_p}
   =\operatorname{diag}(2,\ldots,2,-2\varepsilon).
\end{equation}
By \eqref{eq:quartic-negative-tangential-eigenvalue}, this domain is neither
real convex nor strictly quaternionic $(n-1)$-convex, while
Corollary~\ref{cor:construction} gives smooth Dirichlet solvability for
every smooth $\varphi$ and positive smooth $f(q,z)$ with $f_z\ge0$.

\medskip
\Needspace{6\baselineskip}
\noindent\textbf{Example 2. A periodic domain.}
\par\nobreak\smallskip\nobreak
\noindent
Let $\Gamma_x\subset\RR^n$ be a full lattice and set
\[
 X=(\RR^n/\Gamma_x)_x\times\RR^{2n}_{y,z}
       \times(\RR^n/\mathbb Z^n)_w.
\]
For the same $k$ and $\varepsilon$, replace $|q_a|^2$ in
\eqref{eq:quartic-domain} by $y_a^2+z_a^2$ and define
$\Omega_\varepsilon^{\mathrm{per}}=\{\rho_\varepsilon^{\mathrm{per}}<0\}$ by
\begin{equation}\label{eq:periodic-quartic-domain}
 \rho_\varepsilon^{\mathrm{per}}
 =\sum_{a=1}^{n-1}(y_a^2+z_a^2)
   +(y_n^2+z_n^2)^2-\varepsilon(y_n^2+z_n^2)-1.
\end{equation}
As $|y|^2+|z|^2\to\infty$, we have
$\rho_\varepsilon^{\mathrm{per}}(y,z)\to\infty$. On
$\{\rho_\varepsilon^{\mathrm{per}}=0\}$,
\begin{equation}\label{eq:periodic-quartic-regularity}
 \begin{aligned}
 \sum_{a=1}^n(y_a\partial_{y_a}+z_a\partial_{z_a})
       \rho_\varepsilon^{\mathrm{per}}
 &=2\sum_{a<n}(y_a^2+z_a^2)
       -2\varepsilon(y_n^2+z_n^2)+4(y_n^2+z_n^2)^2\\
 &=2+2(y_n^2+z_n^2)^2>0.
 \end{aligned}
\end{equation}
Equation \eqref{eq:periodic-quartic-regularity} gives
$D\rho_\varepsilon^{\mathrm{per}}\ne0$ on this level set. As after
\eqref{eq:quartic-domain-regularity}, \eqref{eq:periodic-quartic-domain}
with $\rho_\varepsilon^{\mathrm{per}}(0,0)=-1$ yields a bounded star-shaped
$(y,z)$ sublevel set with smooth connected boundary. Compactness of
the $x,w$ tori gives $\Omega_\varepsilon^{\mathrm{per}}\Subset X$.

Each nonconstant summand in \eqref{eq:periodic-quartic-domain} depends
on a single pair $(y_a,z_a)$, so mixed real derivatives involving two
distinct pairs vanish. Since the function is independent of $x,w$,
\eqref{eq:hessian-definition} and \eqref{eq:hessian-component-expansion} give
\[
 H(\rho_\varepsilon^{\mathrm{per}})_{ab}=0\quad(a\ne b),\qquad
 H(\rho_\varepsilon^{\mathrm{per}})_{aa}
 =\tfrac14(\partial_{y_a}^2+\partial_{z_a}^2)
       \rho_\varepsilon^{\mathrm{per}}.
\]
Using
\[
 (\partial_{y_n}^2+\partial_{z_n}^2)(y_n^2+z_n^2)^2
 =16(y_n^2+z_n^2),
\]
we obtain
\begin{equation}\label{eq:periodic-quartic-hessian}
 H(\rho_\varepsilon^{\mathrm{per}})
 =\operatorname{diag}(1,\ldots,1,4(y_n^2+z_n^2)-\varepsilon).
\end{equation}
At $q\in\partial\Omega_\varepsilon^{\mathrm{per}}$, use the same basis
change as in \eqref{eq:quartic-tangential-block}. For
\eqref{eq:periodic-quartic-hessian}, the tangential eigenvalues are now
$1,\ldots,1,\mu_{\mathrm{per}}$, where $|\zeta|\le1$ and
\[
 \mu_{\mathrm{per}}
 =1-|\zeta|^2+\bigl(4(y_n^2+z_n^2)-\varepsilon\bigr)|\zeta|^2
 \ge-\varepsilon.
\]
Thus the expansion in \eqref{eq:quartic-tangential-sigma} becomes,
for $1\le j\le k-1$,
\begin{equation}\label{eq:periodic-quartic-tangential-sigma}
 \begin{aligned}
 \sigma_j\bigl(H(\rho_\varepsilon^{\mathrm{per}})(q)|_{E_q}\bigr)
 &=\binom{n-2}{j}+\mu_{\mathrm{per}}\binom{n-2}{j-1}\\
 &\ge\binom{n-2}{j}-\varepsilon\binom{n-2}{j-1}>0.
 \end{aligned}
\end{equation}
Equations \eqref{eq:periodic-quartic-tangential-sigma} and
\eqref{eq:levi-invariance} give strict quaternionic $(k-1)$-convexity.
At the boundary point $p$ with $y_1=1$ and all other
coordinates zero, $E_p=\{\xi:\xi_1=0\}$ and
$|D\rho_\varepsilon^{\mathrm{per}}(p)|=2$, so
\begin{equation}\label{eq:periodic-quartic-negative-levi}
 R_p=\tfrac12\operatorname{diag}(1,\ldots,1,-\varepsilon).
\end{equation}
By \eqref{eq:periodic-quartic-negative-levi}, strict quaternionic
pseudoconvexity fails, and \eqref{eq:levi-invariance} excludes a strictly
quaternionic plurisubharmonic defining function.

For comparison, the periodic annuli in
\cite[Section~7]{MaZhu2026}, with the $w$ torus added, have the form
\[
 \Omega_{\mathrm{ann}}
 =\{q\in X: \delta_-^2<|y|^2+|z|^2< \delta_+^2\},\qquad 0< \delta_-< \delta_+.
\]
Near the inner boundary, $\delta_-^2-|y|^2-|z|^2$ is a defining function.
Its Hessian is $-I_n$ and its gradient has length $2\delta_-$ on this
boundary. Hence \eqref{eq:levi-invariance} gives
\begin{equation}\label{eq:annular-negative-levi}
 R_q=- \frac{1}{2\delta_-}I_{n-1},\qquad
 \sigma_1(R_q)=- \frac{n-1}{2\delta_-}<0
 \quad\text{on }\{|y|^2+|z|^2= \delta_-^2\}.
\end{equation}
Equation \eqref{eq:annular-negative-levi} excludes
\eqref{eq:levi-condition} for every $k\ge2$.
Nevertheless, $u_0=|y|^2+|z|^2$ satisfies
$H(u_0)=I_n$ and $\sigma_k(H(u_0))=\binom nk$.
It is therefore a subsolution when $f(q,u_0)\le\binom nk$, with
boundary values $\delta_-^2$ and $\delta_+^2$ on the inner and outer components,
respectively.

For the present domain, choose
$C\ge\max_{\overline{\Omega_\varepsilon^{\mathrm{per}}}}(|y|^2+|z|^2)$.
Set
\begin{equation}\label{eq:periodic-global-admissible-function}
 \chi=|y|^2+|z|^2-C\le0
 \quad\text{on }\overline{\Omega_\varepsilon^{\mathrm{per}}},
 \qquad H(\chi)=I_n\in\Gamma_k.
\end{equation}
This function is global on $X$, being independent of $x,w$.
Equations \eqref{eq:periodic-quartic-tangential-sigma} and
\eqref{eq:periodic-global-admissible-function} permit the construction
in Subsection~\ref{sec:subsolution-construction} and
 Proposition~ \ref{prop:quotient-defining-function}, giving a smooth
$k$-admissible defining function and a subsolution for each smooth boundary value.
Theorem~\ref{thm:quotient} therefore gives smooth Dirichlet solvability
for every smooth boundary value and positive smooth $f(q,z)$ with
$f_z\ge0$.

\subsection*{Acknowledgements} Xi-Nan Ma and Yang Zhou were supported by the National Key Research and Development Program of China (Grant No. 2025YFA1017600) and NSFC No. 12141105. 

\subsection*{Statement of usage of artificial intelligence} The central mathematical ideas and overall
approach are due to the authors. OpenAI ChatGPT was used to assist with proofreading
and improving the clarity of the exposition. The authors take full responsibility for the
paper’s content and the accuracy of its references.

\end{document}